\documentclass[letter,11pt]{article}
\usepackage{graphicx}
\usepackage{xcolor}
\usepackage{amsmath}
\usepackage{amsfonts}
\usepackage{cite}
\usepackage{ascmac}
\usepackage{comment}
\usepackage{amsthm}
\usepackage{hyperref}
\usepackage{appendix}

\usepackage{ulem}

\newtheorem{thm}{Theorem}

\newtheorem{la}[thm]{Lemma}

\newtheorem{cor}[thm]{Corollary}
\newtheorem{prop}[thm]{Proposition}

\newtheorem{conj}{Conjecture}

\newcommand{\white}{\qquad{\framebox{\rule{0pt}{4pt}}}}
\newcommand{\owari}{\hfill\white}

\theoremstyle{definition}
\newtheorem{rem}{Remark}

\def\NP{{$\mathcal{NP}$}}

\theoremstyle{definition}
\newtheorem{alg}{Algorithm}
\newtheorem{prob}{Problem}

\DeclareMathOperator{\Od}{O}

\DeclareMathOperator{\BT}{BT}
\DeclareMathOperator{\Peb}{Peb}
\DeclareMathOperator{\col}{col}

\DeclareMathOperator{\Alg}{Alg}
\DeclareMathOperator{\Cen}{Centroid}
\DeclareMathOperator{\CenE}{CentroidEdge}
\DeclareMathOperator{\CenV}{CentroidVertex}
\DeclareMathOperator{\Joi}{Joint}
\DeclareMathOperator{\JoiV}{JointVertex}
\DeclareMathOperator{\JoiE}{JointEdge}
\DeclareMathOperator{\NS}{NodeSet}
\DeclareMathOperator{\PacCon}{PackedConfig}
\DeclareMathOperator{\Pac}{Packing}
\DeclareMathOperator{\CutOff}{CutOff}
\DeclareMathOperator{\Extraction}{Extraction}
\DeclareMathOperator{\Main}{Main}

\DeclareMathOperator{\SP}{SP}
\DeclareMathOperator{\D}{D}
\DeclareMathOperator{\CL}{CL}

\newcommand{\Ord}{\mathrm{O}}
\newcommand{\ord}{\mathrm{o}}

\newcommand{\OPT}{\mathrm{OPT}}

\newcommand{\qedgo}{\vspace{0.35cm}}

\newcommand{\tume}{\hspace*{0.2em}}

\def\NP{{$\mathcal{NP}$}}

\def\deg{\textrm{deg}}
\newcommand{\maru}[1]{{\ooalign{\hfill$\small#1$\hfill\crcr\scalebox{1.3}{$\bigcirc$}}}}

\newenvironment{alg-enumerate}{%
  \begin{enumerate}%
 }{\end{enumerate}%
}
\newenvironment{(R)-enumerate}{%
  \begin{enumerate}%
 }{\end{enumerate}%
}

\newenvironment{(A)-enumerate}{%
  \begin{enumerate}%
 }{\end{enumerate}%
}

\newenvironment{Proofof}[1]{\expandafter
\paragraph{{\it Proof of #1.}}}{\hspace{5mm}\owari\qedgo}

  \newenvironment{namelist}[1]{%
    \begin{list}{}
      {
      \settowidth{\labelwidth}{#1}
      \setlength{\leftmargin}{1.1\labelwidth}
      \setlength{\labelsep}{5pt}
      \setlength{\parsep}{0pt}}
  }{%
  \end{list}}

\def\CAL{\mathcal}
    
 \def\cF{{\CAL F}}

  \def\cS{{\CAL S}}

\begin{document}

\thispagestyle{empty}

\begin{center}
%{\textcolor{blue}{\Large The Structure of Cubic Complexity in Pebble Motion}}
 {\Large Structural Origins of Cubic Complexity in Pebble Motion}
% {\Large Characterizing Cubic Complexity in Pebble Motion}
%    {\Large On the order of the shortest solution sequences for the pebble motion problems}
\end{center}

%%%%%%%%%%%%%%%%%%%%%%%%%%%%%%%%%%%%%%%%%%%%%%%%%%%%%%%%%%%%%%%%%%%%%%%%%%%%%%%%%
%%%%%%%%%%%%%%%%%%%%%%%%%%%%%%%%%%%%%%%%%%%%%%%%%%%%%%%%%%%%%%%%%%%%%%%%%%%%%%%%%
%%%%%%%%%%%%%%%%%%%%%%%%%%%%%%%%%%%%%%%%%%%%%%%%%%%%%%%%%%%%%%%%%%%%%%%%%%%%%%%%%
%%%%%%%%%%%%%%%%%%%%%%%%%%%%%%%%%%%%%%%%%%%%%%%%%%%%%%%%%%%%%%%%%%%%%%%%%%%%%%%%%
%%%%%%%%%%%%%%%%%%%%%%%%%%%%%%%%%%%%%%%%%%%%%%%%%%%%%%%%%%%%%%%%%%%%%%%%%%%%%%%%%
%%%%%%%%%%%%%%%%%%%%%%%%%%%%%%%%%%%%%%%%%%%%%%%%%%%%%%%%%%%%%%%%%%%%%%%%%%%%%%%%%
%%%%%%%%%%%%%%%%%%%%%%%%%%%%%%%%%%%%%%%%%%%%%%%%%%%%%%%%%%%%%%%%%%%%%%%%%%%%%%%%%
%%%%%%%%%%%%%%%%%%%%%%%%%%%%%%%%%%%%%%%%%%%%%%%%%%%%%%%%%%%%%%%%%%%%%%%%%%%%%%%%%
%%%%%%%%%%%%%%%%%%%%%%%%%%%%%%%%%%%%%%%%%%%%%%%%%%%%%%%%%%%%%%%%%%%%%%%%%%%%%%%%%
%%%%%%%%%%%%%%%%%%%%%%%%%%%%%%%%%%%%%%%%%%%%%%%%%%%%%%%%%%%%%%%%%%%%%%%%%%%%%%%%%
%\maketitle 
\begin{center}
{\large Tomoki Nakamigawa
\footnote{This work was supported by JSPS KAKENHI Grant Number JP16K05260.\\
Email: {\texttt{nakamigwt@gmail.com}}}
%Email: {\texttt{nakami@info.shonan-it.ac.jp}}}
}
\end{center}
%\maketitle 
\begin{center}
{
Department of Information Science \\%
Shonan Institute of Technology \\%
1-1-25 Tsujido-Nishikaigan, Fujisawa 251-8511, Japan
}
\end{center}
\begin{center}
{\large Tadashi Sakuma
\footnote{This work was supported by JSPS KAKENHI Grant Number JP16K05260, JP26400185, JP18K03388.\\
Email: {\texttt{sakuma@sci.kj.yamagata-u.ac.jp}}}
}
\end{center}
\begin{center}
{
Faculty of Science \\%
Yamagata University \\%
1-4-12 Kojirakawa, Yamagata 990-8560, Japan
}
\end{center}
%%%%%%%%%%%%%%%%%%%%%%%%%%%%%%%%%%%%%%%%%%%%%%%%%%%%%%%%%%%%%%%%%%%%%%%%%%%%%%%%%
%%%%%%%%%%%%%%%%%%%%%%%%%%%%%%%%%%%%%%%%%%%%%%%%%%%%%%%%%%%%%%%%%%%%%%%%%%%%%%%%%
%%%%%%%%%%%%%%%%%%%%%%%%%%%%%%%%%%%%%%%%%%%%%%%%%%%%%%%%%%%%%%%%%%%%%%%%%%%%%%%%%
%%%%%%%%%%%%%%%%%%%%%%%%%%%%%%%%%%%%%%%%%%%%%%%%%%%%%%%%%%%%%%%%%%%%%%%%%%%%%%%%%
%%%%%%%%%%%%%%%%%%%%%%%%%%%%%%%%%%%%%%%%%%%%%%%%%%%%%%%%%%%%%%%%%%%%%%%%%%%%%%%%%
%%%%%%%%%%%%%%%%%%%%%%%%%%%%%%%%%%%%%%%%%%%%%%%%%%%%%%%%%%%%%%%%%%%%%%%%%%%%%%%%%
%%%%%%%%%%%%%%%%%%%%%%%%%%%%%%%%%%%%%%%%%%%%%%%%%%%%%%%%%%%%%%%%%%%%%%%%%%%%%%%%%
%%%%%%%%%%%%%%%%%%%%%%%%%%%%%%%%%%%%%%%%%%%%%%%%%%%%%%%%%%%%%%%%%%%%%%%%%%%%%%%%%
%%%%%%%%%%%%%%%%%%%%%%%%%%%%%%%%%%%%%%%%%%%%%%%%%%%%%%%%%%%%%%%%%%%%%%%%%%%%%%%%%
%%%%%%%%%%%%%%%%%%%%%%%%%%%%%%%%%%%%%%%%%%%%%%%%%%%%%%%%%%%%%%%%%%%%%%%%%%%%%%%%%

\vspace{10pt}

\begin{abstract}
The pebble motion problem (PMP) asks whether one configuration of
labeled pebbles on a graph can be transformed into another by moving
pebbles to adjacent unoccupied vertices. It is closely related to
multi-agent path finding (MAPF).

Kornhauser, Miller, and Spirakis (FOCS 1984) established the tight
$\Theta(N^3)$ worst-case shortest solution length on $N$-vertex graphs
and suggested that cubic complexity was exceptional. Its structural
origin remains unclear. In particular, whether the cubic worst-case
bound is tight for pebble motion on trees (PMT) remained a major open
question for over four decades.

We establish the first $\ord(N^3)$ upper bound for PMT: an
output-sensitive algorithm constructs a solution sequence of length
$\Ord(N^2\log N)$ for every feasible instance on an $N$-vertex tree,
matching the known $\Omega(N^2)$ lower bound up to a logarithmic factor.

Kornhauser, Miller, and Spirakis also raised the
problem of efficiently approximating the number
of moves required on general board graphs.
In this direction, we extend our tree method
to a hybrid algorithm yielding structure-dependent
bounds on solution length.
In particular, when $2\leq N-n=\Ord(1)$,
every instance on a feasible puzzle admits
a solution sequence of length
$\Ord(n\,\CL(G)\,\D(G))$,
where $n$ is the number of pebbles, $\D(G)$
the diameter, and $\CL(G)$ the maximum length
of a shortest cycle through a vertex.
The above order matches the classical
$\Ord(N^{3/2})$ bound for square grids,
which they cited as an example of subcubic complexity.

Cubic complexity requires $\Ord(1)$ blanks.
Under this condition, we develop an efficient reconfiguration algorithm
across cut vertices. We also prove that subcubic reconfiguration
complexity is preserved under arbitrary ear extensions to biconnected
graphs. Using these results, we show that cubic complexity requires
a unique biconnected component $X$ of order $\Theta(N)$ whose own
reconfiguration cost is $\Theta(N^3)$. Moreover, for each cut vertex
of $G$ belonging to $X$, the union of all structures attached to $X$
through that vertex has order $\Ord(1)$.
We show that feasible puzzles on board graphs obtained by attaching
connected graphs of order $\Ord(1)$ along $\Ord(1)$ individual edges
of a cycle of order $\Theta(N)$ exhibit cubic complexity, whereas two
cycles of order $\Theta(N)$ sharing a single vertex have worst-case
shortest solution length $\Theta(N^2\log N)$.
Our structural results further imply that $X$ contains no cycle of
unbounded sublinear order and that intersecting cycles of order
$\Theta(N)$ in $X$ share at least $N-n$ vertices.
We conjecture that the construction with bounded-size attachments
along $\Ord(1)$ cycle edges characterizes the cubic regime.
\end{abstract}

\newpage

{\small
\tableofcontents
}

\section{Introduction}\label{intro}

\subsection{Definition of the Pebble Motion Problem}

Let $G$ be a finite simple undirected graph.
We denote its vertex and edge sets by $V(G)$ and $E(G)$, and write
$N:=|V(G)|$ for the number of vertices.
Let $P=\{1,\ldots,n\}$ be a set of labeled pebbles with $n\leq N-1$.
Let $\D(G)$ denote the diameter of $G$.
For each vertex $v\in V(G)$, let $\CL(G,v)$ denote the length of a
shortest cycle containing $v$ if such a cycle exists, and set
$\CL(G,v):=1$ otherwise. Define $\CL(G):=\max_{v\in V(G)}\CL(G,v)$.

In the context of the Pebble Motion Problem, the graph $G$ serves as the
\textit{board graph} on which the pebbles move.
When the board graph $G$ is a tree, we refer to it as the
\textit{board tree}.
A \textit{configuration} of $P$ on $G$ is a function
$f:V(G)\to\{0,1,\ldots,n\}$ such that $|f^{-1}(i)|=1$ for each
$1\leq i\leq n$.
Thus $f^{-1}(i)$ is the vertex occupied by pebble $i$, while
$f^{-1}(0)$ is the set of unoccupied vertices.
We also refer to an unoccupied vertex as a \textit{blank}.
A \textit{move} shifts a pebble from its current vertex to an adjacent
blank vertex.

An instance of the \textit{Pebble Motion Problem} (PMP) consists of a
board graph $G$, a set of pebbles $P$, and an initial and a target
configuration.
A sequence of moves transforming the initial configuration into the
target configuration is called a \textit{solution sequence}, and its
\textit{length} is the number of moves it contains.
An instance is called \textit{feasible} if a solution sequence exists;
otherwise it is called \textit{infeasible}.
The classical ``15-puzzle,'' popularized by Sam Loyd
\cite{loyd1959mathematical}, is a standard example in which the board
graph is a $4\times4$ grid
(cf.~\cite{aa887662-c3e2-30b0-a475-bcd8c655206a,
WOS:000083999100001}).

The pebble motion problem is a fundamental model of reconfiguration
on graphs.
It is closely related to \textit{Multi-Agent Path Finding} (MAPF)
and \textit{Coordinated Motion Planning} (CMP).
These models share the basic setting of coordinating multiple agents
on a graph under exclusion constraints, although they differ in their
movement rules, optimization objectives, and scheduling models.
Classical graph theory and graph algorithms often seek to understand a
discrete structure through static invariants of the underlying graph.
Motion and reconfiguration problems introduce a complementary viewpoint:
the structure of a graph is also reflected in the collective dynamics
that it supports when multiple distinguishable objects move under local
constraints.
The pebble motion problem is one of the simplest models in which this
interaction between static graph structure and global discrete dynamics
can be studied exactly.

\subsection{Background}

The algorithmic theory of the pebble motion problem was initiated by
Kornhauser, Miller, and Spirakis in their seminal FOCS~1984 paper
\cite{715921}.
They showed that, on general graphs, every feasible instance admits
a solution sequence of length $\Ord(N^3)$ and also exhibited families
of instances requiring $\Omega(N^3)$ moves.
Thus the worst-case length of a shortest solution sequence is
$\Theta(N^3)$.
They already suggested, however, that this cubic behavior might be
exceptional: ``... it seems that only a small fraction of the graph
puzzles actually require $\Ord(N^3)$ moves.''

Kornhauser, Miller, and Spirakis~\cite{715921} also considered the
special case in which the board graph is a tree.
We refer to this case as the \textit{Pebble Motion Problem on Trees}
(PMT). 
They stated that the worst-case length of a shortest solution sequence
for PMT is $\Ord(N^2)$. 
Their argument, however, relies on the assertion that, on a tree, an
arbitrary pebble can be moved to an arbitrary vertex using at most
$\Ord(N)$ moves. 
This assertion is false in general, and consequently their argument
does not establish the claimed $\Ord(N^2)$ bound; see
Remark~\ref{counterexample} in the Appendix for a concrete
counterexample. 
Auletta, Monti, Parente, and Persiano~\cite{WOS:000077614700003}
later remarked that they were unable to verify this quadratic claim.
Consequently, the only general worst-case upper bound available for
trees remained the $\Ord(N^3)$ bound inherited from the algorithm of
Kornhauser, Miller, and Spirakis~\cite{715921}.
Even the recent algorithm of Ardizzoni, Saccani, Consolini, Locatelli,
and Nebel~\cite{WOS:001177663400002} has worst-case complexity
$\Ord(N^3)$ on trees.

On the other hand, the best known lower bound for PMT is
$\Omega(N^2)$, obtained from instances in which $\Theta(N)$ pebbles
must each traverse a path of length $\Theta(N)$.
Let $\OPT(N)$ denote the maximum length of a shortest solution sequence
over all feasible instances on $N$-vertex trees.
Before the present work, therefore, the extremal complexity of PMT
remained within the gap $\Omega(N^2)\leq \OPT(N)\leq \Ord(N^3)$.
Thus the asymptotic complexity of pebble motion on trees remained
unresolved for more than four decades.

Trees are also structurally restrictive from the viewpoint of
feasibility.
Because every edge of a tree is a bridge, opportunities for exchanging
pebbles through alternative routes are severely limited.
In particular, if the number of blanks is too small relative to the
length of a corridor, that is, a maximal path consisting of consecutive
degree-two vertices, some target configurations may become unreachable.
This type of obstruction does not occur in $2$-connected graphs other
than cycles when at least two blanks are available.
The tree case therefore provides a natural testing ground for
understanding how restricted graph structure affects both feasibility
and the length of shortest solution sequences.

The broader feasibility theory of the Pebble Motion Problem also has a long
history.
Wilson~\cite{zbMATH03446929} characterized the case of a single blank vertex,
and Kornhauser, Miller, and Spirakis~\cite{715921} gave a
polynomial-time feasibility test for undirected graphs with at least
two blank vertices.
The directed setting is substantially harder.
Nebel~\cite{WOS:001155044500001} proved that reachability is \NP-complete on
general directed graphs, while polynomial-time algorithms remain
available on strongly connected digraphs.

A closely related recent direction concerns the complexity of
graph-based multi-agent motion. 
Coordinated motion of labeled robots on grids, including approximation and hardness results
under parallel movement rules, was studied in \cite{WOS:000546869700004}. 
Eiben, Ganian, and Kanj~\cite{eiben_et_al:LIPIcs.SoCG.2023.28} studied parameterized
Coordinated Motion Planning on rectangular grids.
Deligkas, Eiben, Ganian, Kanj, and Ramanujan~\cite{deligkas_et_al:LIPIcs.ICALP.2024.53} extended
this line to general graphs, including the objective of minimizing
total travel distance.
Subsequent work has considered graph-based coordinated motion under
a prescribed move budget~\cite{deligkas_et_al:JCSS.2026} and parameterized
Multi-Agent Path Finding specifically on trees~\cite{deligkas_et_al:AAMAS.2025}.
Trees are also algorithmically relevant beyond the tree-restricted
problem itself: several approaches to MAPF solve or approximate more
general instances by reducing or reformulating them as pebble-motion
problems on suitably defined trees~\cite{Khorshid2011APA,Krontiris2013FromFT,WOS:000948128106003,WOS:001177663400002},
and pebble-motion algorithms have themselves been used as fast MAPF
solvers~\cite{WOS:001225902400001}.

For the sequential PMP movement model considered here, solution length
coincides exactly with the total travel distance of the pebbles; this
correspondence is recorded formally in Section~\ref{pre}.

Exact optimization, however, remains computationally difficult even
on highly structured board graphs.
Ratner and Warmuth~\cite{zbMATH04155895} proved that finding a shortest solution
sequence is \NP-hard even for square-grid instances with a single
blank vertex.
This is particularly instructive because square grids are also a
basic example on which substantially subcubic solution-length bounds
are possible.
More recently, Geft~\cite{Gef26} proved that exact PMT optimization is
\NP-hard even on subdivided stars of depth~3 and, separately, on
trees of maximum degree~3, while Koyfman, Atzmon, Shperberg, and
Felner~\cite{KASF26} independently established related hardness
results for MAPF optimization on trees.

These hardness results make worst-case solution-length analysis a
natural complementary objective.
When an exact optimum cannot in general be computed efficiently, it
remains meaningful to ask how much motion is intrinsically sufficient
or necessary on a given structural class of board graphs, and to
construct solution sequences attaining the corresponding asymptotic
bound.
For PMP this distinction is particularly concrete: the quantity being
bounded is not merely the running time of an algorithm but the total
motion performed by the pebbles themselves.
Accordingly, we study the extremal order of the shortest solution
length and the graph-theoretic structures that determine that order.

Beyond the existence of general-graph instances requiring
$\Theta(N^3)$ moves, the structural origin of this cubic behavior has
also remained unclear. In particular, no asymptotic separation between
trees and general graphs had been established, despite the substantial
structural difference between the corresponding motion mechanisms.
This leads to two closely related questions: can the classical cubic
worst-case complexity already occur on trees, and if not, what
additional graph structure is necessary to support $\Theta(N^3)$
behavior? The present paper answers the first question by establishing
the first subcubic upper bound for trees, and then uses the resulting
algorithmic framework to address the second by identifying strong
structural restrictions on graphs that can exhibit cubic complexity.

\subsection{Our Contributions}

In this paper, we investigate the structural origin of the classical
$\Theta(N^3)$ worst-case complexity of the pebble motion problem.
We first show that this cubic behavior cannot occur on trees.
Our first main result is the following.

\noindent\textbf{Theorem\tume\ref{main} (informal).}\textit{
Every feasible instance of the pebble motion problem on an $N$-vertex
tree admits a solution sequence of length $\Ord(N^2\log N)$,
constructible in time linear in its length.
This is the first $\ord(N^3)$ upper bound for PMT.
Together with the known $\Omega(N^2)$ lower bound, it determines the
extremal complexity on trees up to a logarithmic factor.}

We next extend the tree-based method to general graphs.
Theorem\tume\ref{refinement2-general-constant-blanc} shows that, when
$2\leq N-n=\Ord(1)$, every instance on a feasible puzzle admits a solution sequence
of length
$\Ord\bigl(n\,\CL(G)\,\D(G)\bigr)$.
This gives a structural refinement of the classical cubic bound in
the direction suggested by Kornhauser, Miller, and
Spirakis~\cite{715921}, and matches their
$\Ord(N^{3/2})$ bound for square grids.
More generally, Theorem\tume\ref{refinement3} implies that cubic
behavior is possible only when the number of blanks is bounded by a
constant.

We then develop additional algorithms for this constant-blank regime
and derive the following structural conclusion.

\noindent\textbf{Theorem\tume\ref{thm-Theta-n3-unique-large-component} (informal).}\textit{
If the worst-case shortest solution length on an $N$-vertex graph is
$\Theta(N^3)$, then the graph has a unique biconnected component $X$
of order $\Theta(N)$.
All structures attached to $X$ through any one of its cut vertices
have bounded total size.
Moreover, $X$ contains no cycle of unbounded sublinear order, and any
two intersecting cycles of order $\Theta(N)$ have at least as many
common vertices as there are blanks.}

Finally, Theorem\tume\ref{thm-Theta-n3} shows that cubic complexity
does occur when a cycle of order $\Theta(N)$ forms the backbone of the
graph and the total number of additional vertices in attachments along
the cycle is $\Ord(1)$.
In contrast, Theorem\tume\ref{thm-two-cycles-Theta-n2logn} shows that
two cycles of order $\Theta(N)$ meeting at a single vertex have
worst-case shortest solution length exactly $\Theta(N^2\log N)$.
Motivated by these results, we formulate
Conjecture\tume\ref{conj-Theta-n3}, which asserts that the former
construction essentially characterizes the cubic regime.

\subsection{High-level Technical Overview of the Tree Algorithm}

Throughout this overview, let $b:=N-n$ denote the number of blank
vertices.

Our algorithm has two layers, corresponding to two distinct tasks:
efficient routing inside a local region and global management of the
blank vertices that make such local routing repeatedly possible.

The \textbf{inner layer} is the routing algorithm established in
Theorem\tume\ref{sbs3}.
Whenever a connected region containing sufficiently many blank
vertices is available, it provides an efficient local routing
primitive.
Intuitively, the blank vertices form the local working space in which
the pebbles can be rearranged.
Internally, this primitive uses a divide-and-conquer organization:
the pebbles are partitioned into smaller groups, the corresponding
subproblems are solved recursively, and the resulting configurations
are merged.
This local procedure alone, however, does not yield the main tree
algorithm, because the required workspace is neither globally
available nor automatically preserved after each routing task.

The main algorithmic contribution in the tree case is the
\textbf{outer layer}, which is fundamentally different from
divide-and-conquer.
The current unresolved part of the tree is viewed as consisting of a
central region together with peripheral subtrees.
At each stage, the algorithm selects one such peripheral target
subtree, moves the required pebbles into their target positions there,
and then declares that subtree permanently completed.
The completed subtree is never used again, while the remaining tree,
together with the relocated blank workspace, becomes the next
instance.
Thus the outer layer repeatedly reduces a single global instance,
rather than decomposing it into independent subproblems.

The principal difficulty is that Theorem\tume\ref{sbs3} requires a
connected region containing sufficiently many blank vertices. The
outer algorithm therefore maintains the \textbf{Workspace Invariant}:
throughout the execution, there exists a connected workspace
containing sufficiently many blank vertices to invoke
Theorem\tume\ref{sbs3} on the current local subproblem. After a target
subtree has been completed, the algorithm moves this workspace into
the remaining tree without disturbing any previously completed
region. A centroid decomposition is used to choose the successive
regions so that their sizes decrease in a controlled way and the
workspace can be relocated from one region to the next. At the same
time, the reduction satisfies the \textbf{Monotone-Progress Property}:
once a target subtree has been completed, its vertices are never
modified again. Hence every stage permanently decreases the unresolved
part of the tree, until the remaining instance is small enough for a
final application of Theorem\tume\ref{sbs3}.

The cost analysis reflects this two-level organization. Each
invocation of Theorem\tume\ref{sbs3} operates on a local region of
size $\Theta(b)$ and hence requires $\Ord(b^2\log b)$ moves. The
centroid-based reduction produces $\Ord(N/b)$ target subtrees, and for
each target subtree the workspace passes through only $\Ord(N/b)$
temporary regions of size $\Theta(b)$. Thus the total number of moves
is $\Ord((N/b)(N/b)b^2\log b)=\Ord(N^2\log b)$. The resulting
solution sequence can furthermore be generated in $\Ord(N+L)$ time,
where $L$ denotes its length. The essential point is that the local
routing primitive becomes globally useful precisely because the outer
algorithm preserves and relocates the blank workspace while completed
regions are fixed permanently.

\subsection{High-level Technical Overview of the General-Graph Structural Results}

The tree algorithm also provides the starting point for our analysis
of general board graphs.
Given a board graph $G$, we first construct a breadth-first-search
spanning tree $T_G$ and attempt to route the pebbles by applying the
tree procedure to $T_G$.
The difficulty is that a feasible instance on $G$ need not remain
feasible when restricted to $T_G$: a long corridor in the spanning
tree may contain too few blanks to permit the required rearrangement.
When this obstruction occurs, we use edges outside $T_G$ to route
pebbles around an appropriate cycle of $G$, thereby restoring the
local feasibility conditions needed by the tree procedure.

This combination of tree routing and cycle-based repairs yields, when
$2\leq N-n=\Ord(1)$ and $(G,P)$ is a feasible puzzle, a solution
sequence of length $\Ord\bigl(n\,\CL(G)\,\D(G)\bigr)$.
Thus geometric information about the board graph directly translates
into an upper bound on the motion complexity.
This provides a general structural refinement of the classical cubic
bound in the direction proposed by Kornhauser, Miller, and
Spirakis~\cite{715921} in their concluding open problems, where they
observed that only a small fraction of graph puzzles appear to require
cubic-length solutions and pointed to square grids as a principal
subcubic example.
Indeed, on a $\sqrt N\times\sqrt N$ square grid,
$\CL(G)=4$ and $\D(G)=\Theta(\sqrt N)$, giving
$\Ord(N^{3/2})$ moves and matching their classical bound.
%%%%%%%%%%%%%%%%%%%%%%%%%%%%%%%%%%%%%%%%%%%%%%%%%%%%%%%%%%%%%%%%%%%
Related coordinated-motion results on grids obtain comparable
distance-scale bounds under parallel movement rules
\cite{WOS:000546869700004}, whereas our bound applies under the
sequential PMP rule and extends beyond grids through the parameter
$\CL(G)$.
%%%%%%%%%%%%%%%%%%%%%%%%%%%%%%%%%%%%%%%%%%%%%%%%%%%%%%%%%%%%%%%%%%%
More generally, the same analysis shows that cubic complexity is
impossible when the number of blanks is unbounded.
Hence, in studying the $\Theta(N^3)$ regime, we may restrict attention
to a constant number of blanks.

The remainder of the structural analysis exploits the decomposition
of $G$ into its biconnected components and the connected components of
the forest formed by all bridges of $G$, augmented by the cut vertices
through which they attach to the rest of the graph; we call the latter
\textit{bridge-tree components}. 
We refer to the collection of these biconnected components and
bridge-tree components as the 
\textit{biconnected--bridge-tree decomposition} of $G$.
The main issue is no longer merely how to route a prescribed pebble
locally, but how efficiently pebbles can be transferred between these
components through their common cut vertices.
For a component under consideration, we call such a cut vertex through
which it is connected to the rest of the graph an
\textit{interface vertex}.
For this purpose, we introduce a procedure that repeatedly brings
pebbles to an interface vertex in a controlled order.
Its key property is that every relevant pebble is presented at the
interface exactly once while the blanks remain concentrated in a
connected local workspace.

This presentation order can then be used as an address system.
If a component can present its pebbles at a cut vertex in a prescribed
order, selected pebbles can be replaced there by pebbles arriving from
an adjacent component, after which the presentation procedure is run
in reverse.
In this way, local presentation sequences can be
\textit{transcribed} into global reconfiguration sequences without
reconstructing the entire configuration from scratch.
This mechanism allows us to compare the reconfiguration cost of a
large component with that of the whole graph.

Applying this principle repeatedly rules out intermediate-scale
structures: a component of order $\omega(1)$ but $\ord(N)$ would
already permit a subcubic global reconfiguration, and suitable ear
extensions of a subcubic component also remain subcubic. Consequently,
if cubic complexity occurs, every component in the
biconnected--bridge-tree decomposition has either bounded order or
order $\Theta(N)$; moreover, there is a unique component of order
$\Theta(N)$, and it must be biconnected. All bridge-tree components
are therefore bounded, and even the total structure attached through
any single cut vertex of this giant component is bounded. 
Finally, the
same reconfiguration mechanisms exclude cycles of order
$\omega(1)\cap\ord(N)$ inside the giant component and show that any two
cycles of order $\Theta(N)$ are either disjoint or overlap in at least
as many vertices as there are blanks.
%Finally, the
%same reconfiguration mechanisms exclude cycles of order
%$\omega(1)\cap\ord(N)$ inside the giant component and show that any two
%cycles of order $\Theta(N)$ must overlap in at least as many vertices
%as there are blanks. 
Thus cubic complexity forces essentially all
unbounded structure into one giant biconnected component whose
internal cycle structure is itself highly restricted. The sharpness
results in Section~\ref{lower bound} show that these restrictions are
substantive: a giant cycle with only $\Ord(1)$ vertices in local
attachments can still exhibit $\Theta(N^3)$ complexity, whereas two
cycles of order $\Theta(N)$ meeting at a single cut vertex have
complexity only $\Theta(N^2\log N)$.

\subsection{Organization}

Section~\ref{pre} reviews preliminary definitions and classical
feasibility results, and Section~\ref{results} states the main results
of the paper. Section~\ref{base-cases} establishes the local routing
primitive used throughout the tree algorithm;
Section~\ref{splitting-tree} develops the centroid-based decomposition
needed to organize the local subproblems; and Section~\ref{MainAlg}
presents the global workspace-preserving reduction.
Section~\ref{Main1Proof} proves Theorem\tume\ref{main}, while
Section~\ref{Main2Proof} proves Theorem\tume\ref{main2} for general
board graphs. Section~\ref{refinement} presents refinements of these
two results. Specifically, Subsection~\ref{refinement-tree} proves
Theorem\tume\ref{refinement1}, refining Theorem\tume\ref{main}, and
Subsection~\ref{refinement-general} establishes
Theorem\tume\ref{refinement2}, refining Theorem\tume\ref{main2};
Theorem\tume\ref{refinement2-general-constant-blanc} is an immediate
consequence of Theorem\tume\ref{refinement2}. The same subsection also
establishes Theorem\tume\ref{refinement3} and
Theorem\tume\ref{Restriction}, which give the first structural
restrictions on the cubic regime.

Section~\ref{sec:pebble-presentation} develops the Pebble Presentation
Algorithm and the biconnected--bridge-tree decomposition, and
Section~\ref{sec:structural-restrictions} introduces transcription
across cut vertices and derives further structural restrictions on
graphs exhibiting cubic complexity. Section~\ref{lower bound}
establishes sharp sufficient and separating examples, formulates our
conjectured characterization, and combines the preceding results into
our strongest structural description of the cubic regime.
Section~\ref{final} concludes the paper and discusses the remaining
open problems. Proofs of auxiliary results deferred from
Sections~\ref{base-cases}--\ref{MainAlg} are collected in the Appendix;
the main results are proved in the sections indicated above.

%%%%%%%%%%%%%%%%%%%%%%%%%%%%%%%%%%%%%%%%%%%%%%%%%%%%%%%%%%%%%%%%%%%%%%%%%%%%%%%%
%%%%%%%%%%%%%%%%%%%%%%%%%%%%%%%%%%%%%%%%%%%%%%%%%%%%%%%%%%%%%%%%%%%%%%%%%%%%%%%%
%%%%%%%%%%%%%%%%%%%%%%%%%%%%%%%%%%%%%%%%%%%%%%%%%%%%%%%%%%%%%%%%%%%%%%%%%%%%%%%%

\section{Preliminaries}\label{pre}

In this section, we collect the definitions and classical feasibility
results that will be used throughout the paper.
The basic terminology for the Pebble Motion Problem, including board
graphs, configurations, blanks, moves, and solution sequences, is as
defined in Section~\ref{intro}.

\subsection{Feasibility and isthmuses}\label{sec:feasibility}

We first recall the structural characterization of feasibility due to
Kornhauser, Miller, and Spirakis~\cite{715921}.

A path $I=v_1\cdots v_d$ of a connected graph $G$ is called an
\textit{isthmus} if
\begin{enumerate}
\item every edge of $I$ is a bridge of $G$;
\item every vertex of $I$ is a cut vertex of $G$; and
\item $\deg_G(v_i)=2$ for $1<i<d$.
\end{enumerate}
An isthmus with $d$ vertices is called a $d$-isthmus.
For a connected graph $G$, let
\[
k(G):=\max\{|V(I)|: I \text{ is an isthmus of }G\}
\]
if $G$ contains an isthmus, and set $k(G):=1$ otherwise.
We call $k(G)$ the \textit{maximum isthmus size} of $G$.
Thus $k(G)=1$ whenever $G$ has no isthmus of size at least two;
in particular, this is the case for every $2$-edge-connected graph.
When the graph under consideration is clear from the context, we simply
write $k:=k(G)$.

Given a board graph $G$ and a pebble set $P$, let
$n:=|P|$ and $N:=|V(G)|$.
We say that the pair $(G,P)$ is \textit{transitive} if, for every
configuration $f$ of $P$ on $G$, every pebble can be moved from its
current vertex to any prescribed vertex of $G$ by a finite sequence
of moves.

The following classical theorem characterizes the case in which all
configurations are mutually reachable when at least two blanks are
present.

\medskip
\noindent
\textbf{Theorem A}
(Kornhauser, Miller, and Spirakis~\cite{715921}).
\textit{
Let $2\leq q\leq N-1$ and assume $n+q=N$.
Let $G$ be a connected graph with $N$ vertices that is not a cycle,
and let $P$ be a set of $n$ pebbles placed on the vertices of $G$.
Then the following are equivalent:
\begin{enumerate}
\item $G$ has no $q$-isthmus;
\item the pair $(G,P)$ is transitive;
\item every pair of configurations of $P$ on $G$ can be transformed
      into each other by a finite sequence of moves.
\end{enumerate}
}
\medskip

We call the puzzle $(G,P)$ \textit{feasible} if every pair of
configurations of $P$ on $G$ is mutually reachable.  Under the
assumptions of Theorem~A, this is equivalent to saying that $(G,P)$ is
transitive.

Thus, when $q=N-n\geq 2$, the absence of a $q$-isthmus precisely
characterizes unrestricted mutual reachability, except for the cycle
case excluded in Theorem~A.
In particular, an isthmus whose size is at least the number of blanks
constitutes a barrier across which the pebbles cannot in general be
rearranged freely.
When such an isthmus is present, the board naturally separates into
smaller regions within which the induced pebble-motion problems can
be treated independently.

Throughout the general-graph analysis in this paper, unless explicitly
stated otherwise, we restrict attention to connected board graphs $G$
that are not cycles, to instances with at least two blanks, i.e.,
$N-n\geq 2$, and to feasible puzzles $(G,P)$.  Recall that a
\textit{feasible instance} is one whose target configuration is reachable
from its initial configuration, whereas a \textit{feasible puzzle} is one
in which every two configurations are mutually reachable.  Thus every
instance on a feasible puzzle is feasible, but the converse need not
hold.  In particular, pure cycles and the single-blank regime are
excluded from the general-graph results of this paper.

This leads to the formulation of the general-graph problem used
throughout the paper.

\begin{prob}[\textsc{Pebble Motion Problem on General Graphs}]
\label{PMPonGeneralGraphs}
\hfill
\begin{namelist}{INSTANCE}
\item[{\it Instance}] \rm
A connected graph $G$ with $N$ vertices serving as the board graph,
where $G$ is not a cycle;
a set of pebbles $P:=\{1,\ldots,n\}$ with $n\leq N-2$ such that
the puzzle $(G,P)$ is feasible;
and an initial configuration $f$ and a target configuration $g$ of
$P$ on $G$.

\item[{\it Task}] \rm
Compute a sequence of moves that transforms $f$ into $g$, where each
move shifts a pebble from its current vertex to an adjacent
unoccupied vertex of $G$.
\end{namelist}
\end{prob}

By Theorem~A, under the standing assumptions above, the feasibility of
the puzzle $(G,P)$ is equivalent to the absence of an $(N-n)$-isthmus
in $G$.  In particular, every pair of configurations $f$ and $g$ in
the problem above is mutually reachable.

\subsection{The Pebble Motion Problem on Trees}

We now specialize the preceding characterization to trees.
Let $G$ be a tree and write $k:=k(G)$ for its maximum isthmus size.

This leads to the formulation of the tree problem used throughout the
paper.

%We now specialize the preceding characterization to trees.
%Let $G$ be a tree and write $k:=k(G)$ for its maximum isthmus size.
%Since $k\geq 1$, the condition $N-n>k$ implies $N-n\geq 2$.
%Hence, since a tree is not a cycle, Theorem~A implies that whenever
%$N-n>k$, the pair $(G,P)$ is transitive.
%Equivalently, every configuration of the pebbles can be transformed
%into every other configuration.
%
%This leads to the formulation of the tree problem used throughout the
%paper.

\begin{prob}[\textsc{Pebble Motion Problem on Trees}]
\label{PMPonTree}
\hfill
\begin{namelist}{INSTANCE}
\item[{\it Instance}] \rm
A tree $G$ with $N$ vertices serving as the board graph
(i.e., the board tree);
a set of pebbles $P:=\{1,\ldots,n\}$ with $n\leq N-2$;
and an initial configuration $f$ and a target configuration $g$ of
$P$ on $G$, where $n+k(G)<N$.

\item[{\it Task}] \rm
Compute a sequence of moves that transforms $f$ into $g$, where each
move shifts a pebble from its current vertex to an adjacent
unoccupied vertex of $G$.
\end{namelist}
\end{prob}

%\subsection{The Pebble Motion Problem on Trees}
%
%We now specialize the preceding characterization to trees.
%Let $G$ be a tree, and let $k$ denote the number of vertices in a
%longest isthmus of $G$.
%Since a tree is not a cycle, Theorem~A implies that whenever
%$N-n>k$, the pair $(G,P)$ is transitive.
%Equivalently, every configuration of the pebbles can be transformed
%into every other configuration.
%
%This leads to the formulation of the tree problem used throughout the
%paper.
%
%\begin{prob}[\textsc{Pebble Motion Problem on Trees}]
%\label{PMPonTree}
%\hfill
%\begin{namelist}{INSTANCE}
%\item[{\it Instance}] \rm
%A tree $G$ with $N$ vertices serving as the board graph
%(i.e., the board tree);
%an integer $k$, defined as the number of vertices contained in a
%longest isthmus of $G$;
%a set of pebbles $P:=\{1,\ldots,n\}$ with $n\leq N-1$;
%and an initial configuration $f$ and a target configuration $g$ of
%$P$ on $G$, where $n+k<N$.
%
%\item[{\it Task}] \rm
%Compute a sequence of moves that transforms $f$ into $g$, where each
%move shifts a pebble from its current vertex to an adjacent
%unoccupied vertex of $G$.
%\end{namelist}
%\end{prob}

The condition $n+k(G)<N$ is equivalent to $N-n>k(G)$.
Since $k(G)\geq 1$, this implies $N-n\geq 2$.
Hence, since a tree is not a cycle, Theorem~A guarantees that every
pair of configurations in the instance class above is mutually
reachable.
%The condition $n+k(G)<N$ is equivalent to $N-n>k(G)$.
%Since $k(G)\geq 1$, it also implies $N-n\geq 2$.
%Therefore, by Theorem~A, every pair of configurations in the instance
%class above is mutually reachable.
%The inequality $n+k<N$ is equivalent to $N-n>k$ and therefore
%guarantees, by Theorem~A, that every pair of configurations in the
%instance class above is mutually reachable.
This allows the remainder of the paper to focus on the length and
construction of solution sequences rather than on reachability itself.
More generally, if an instance on a tree is feasible but the puzzle
itself is not feasible, the isthmuses obstructing transitivity decompose
the motion into independent feasible subpuzzles in the sense of the
classical decomposition of Kornhauser, Miller, and
Spirakis~\cite{715921}. Applying the same bounds to these subpuzzles
independently gives the same asymptotic bounds for every feasible
instance on a tree: if the subpuzzles contain $n_i$ pebbles and have
orders $N_i$, then $\sum_i n_iN_i\leq nN$ and
$\sum_i n_i^2\leq n^2$, while distances within each subtree agree with
their distances in the original tree.

\subsection{Solution length and total travel distance}

For a PMP solution, let $\pi_i$ denote the trajectory of pebble $i$.
Since each move sends exactly one pebble across one edge, a solution
sequence of length $L$ satisfies
$L=\sum_{i=1}^{n}\operatorname{length}(\pi_i)$.
Thus, under the sequential PMP movement model, solution length is
exactly the total travel distance of the pebbles.

\section{Results}\label{results}

We now state the main results of the paper.  We first give the
algorithmic upper bounds that separate trees from the general
$\Theta(N^3)$ regime.  We then describe the presentation and
transcription mechanisms developed in
Sections~\ref{sec:pebble-presentation} and
\ref{sec:structural-restrictions}, and the structural restrictions
they impose on cubic instances. Finally, we state the strongest 
structural result and the sharp examples from Section~\ref{lower bound}.

\subsection{Algorithmic upper bounds}

We begin with the \textsc{Pebble Motion Problem on Trees} defined in
Section~\ref{pre}.  
Recall that $n$ denotes the number of pebbles, $N$ the number of
vertices of the board tree, and $k=k(G)$ its maximum isthmus size.
%Recall that $n$ denotes the number of pebbles, $N$
%the number of vertices of the board tree, and $k$ the number of
%vertices in a longest isthmus.

\setcounter{thm}{0}
\begin{thm}\label{main}
A solution sequence to the \textsc{Pebble Motion Problem on Trees}
of length
\[
\Ord\!\left(
nN+n^2\log\bigl(1+\min\{n,k\}\bigr)
\right)
\]
can be computed within the same asymptotic time bound.
\end{thm}

Since $n\leq N$, Theorem\tume\ref{main} implies
$\Omega(N^2)\leq\OPT(N)\leq\Ord(N^2\log N)$, where $\OPT(N)$ denotes
the maximum, over all feasible instances on $N$-vertex trees, of the
length of a shortest solution sequence.  In particular, cubic
worst-case behavior cannot occur on trees.  A refined bound keeping
track of the diameter and the maximum isthmus size is the following.

\begin{thm}\label{refinement1}
A solution sequence of length
\[
\Ord\!\left(
n\D(G)+
\min\left\{
kn\D(G),\,
n^2\log\bigl(1+\min\{n,k\}\bigr)
\right\}
\right)
\]
to the \textsc{Pebble Motion Problem on Trees} can be computed in time
\[
\Ord\!\left(
N+n\D(G)+
\min\left\{
kn\D(G),\,
n^2\log\bigl(1+\min\{n,k\}\bigr)
\right\}
\right).
\]
\end{thm}

Theorem\tume\ref{refinement1} implies, in particular, that if
$k=\Ord(a^{N/n})$ for some positive constant $a$, then the shortest
solution length is in $\Ord(nN)$.

The main result of Ardizzoni, Saccani, Consolini, Locatelli, and
Nebel~\cite{WOS:001177663400002} gives solution sequences of length
$\Ord(N^2+nkN)$ for the \textsc{Pebble Motion Problem on Trees}.
In contrast, Theorem\tume\ref{refinement1} above constitutes a genuine
improvement over their result.

We do not know whether the logarithmic factor in
Theorem\tume\ref{main} is necessary.

\begin{conj}\label{conjecture}
There exist instances of the \textsc{Pebble Motion Problem on Trees}
for which the length of the shortest solution sequence is
$\Theta(N^2\log N)$.
\end{conj}

We next extend the tree algorithm to general board graphs.  Throughout
the remainder of this subsection, we work under the standing
assumptions for general board graphs stated in
Section~\ref{sec:feasibility}.  Recall that $k:=k(G)$ denotes the
maximum isthmus size of $G$; in particular, $k=1$ when $G$ is
$2$-edge-connected.

\begin{thm}\label{main2}
Let $(G,P)$ be a feasible puzzle, and let $f$ and $g$ be
configurations of $P$ on $G$.  Suppose that a shortest cycle through each
vertex $v\in V(G)$, or $\emptyset$ if no such cycle exists, is provided.
Let $\theta$ be any integer satisfying $k\leq\theta\leq N-n-1$.  Then a
solution sequence from $f$ to $g$ of length
$\displaystyle
\Ord\!\left(n\D(G)+\frac{n^2\min\{n,\CL(G)\}}{\theta}
\allowbreak+n^2\log\bigl(1+\min\{n,\theta\}\bigr)\right)$
can be computed in time
$\displaystyle
\Ord\!\left(|E(G)|+n\D(G)+\frac{n^2\min\{n,\CL(G)\}}{\theta}
\allowbreak+n^2\log\bigl(1+\min\{n,\theta\}\bigr)\right).$
\end{thm}

For convenient reference, we record two direct specializations of
Theorem\tume\ref{main2}.  Taking $\theta=N-n-1$ and using
$N-n-1\geq (N-n)/2$ gives a solution sequence of length
\[
\Ord\!\left(
 n\D(G)
 + \frac{n^2\min\{n,\CL(G)\}}{N-n}
 + n^2\log\bigl(1+\min\{n,N-n\}\bigr)
\right),
\]
with computational time $\Ord(|E(G)|)$ plus the same bound.
Alternatively, taking $\theta=k$ gives a solution sequence of length
\[
\Ord\!\left(
 n\D(G)
 + \frac{n^2\min\{n,\CL(G)\}}{k}
 + n^2\log\bigl(1+\min\{n,k\}\bigr)
\right),
\]
again with computational time $\Ord(|E(G)|)$ plus the same bound.

Theorem\tume\ref{main2} is further refined in
Subsection~\ref{refinement-general}; the full refinement is stated there
as Theorem\tume\ref{refinement2}. Since its general form is somewhat
involved, we record here the following particularly important
constant-blank consequence.

\begin{thm}\label{refinement2-general-constant-blanc}
Let $(G,P)$ be a feasible puzzle, and let $f$ and $g$ be
configurations of $P$ on $G$.  Assume that $2\leq N-n=\Ord(1)$, and that
a shortest cycle through each vertex $v\in V(G)$, or $\emptyset$ if no
such cycle exists, is provided.  Then a solution sequence from $f$ to $g$
of length $\Ord(n\CL(G)\D(G))$ can be computed in time
$\Ord(|E(G)|+n\CL(G)\D(G))$.
\end{thm}

For a $\sqrt N\times\sqrt N$ square grid,
$\CL(G)=4$ and $\D(G)=\Theta(\sqrt N)$, so
Theorem\tume\ref{refinement2-general-constant-blanc} gives
$\Ord(N^{3/2})$ moves when $2\leq N-n=\Ord(1)$, matching the
classical square-grid bound of
Kornhauser, Miller, and Spirakis~\cite{715921}.

To find, for every vertex $v$ of a given graph $G$, a shortest cycle
passing through $v$ (or $\emptyset$ if no such cycle exists) requires a
total computational effort of $\Ord(|E(G)||V(G)|)$. However, if it
suffices to find, for each vertex $v$, some cycle passing through $v$,
not necessarily a shortest one, then all such cycles can be obtained
in $\Ord(|E(G)|)$ time. Thus the requirement in
Theorem\tume\ref{main2} that shortest cycles be supplied in advance can
be removed, at the cost of replacing $\min\{n,\CL(G)\}$ by the trivial
upper bound $n$. This gives the following immediate corollary of
Theorem\tume\ref{main2}.

\begin{cor}\label{main2-cor}
Let $(G,P)$ be a feasible puzzle, let $f$ and $g$ be configurations of
$P$ on $G$, and let $\theta$ be any integer satisfying
$k\leq\theta\leq N-n-1$.  Then a solution sequence from $f$ to $g$ can be
computed in time
$\displaystyle
\Ord\!\left(|E(G)|+n\D(G)+\frac{n^3}{\theta}
\allowbreak+n^2\log\bigl(1+\min\{n,\theta\}\bigr)\right),$
and its length is bounded by
$\displaystyle
\Ord\!\left(n\D(G)+\frac{n^3}{\theta}
\allowbreak+n^2\log\bigl(1+\min\{n,\theta\}\bigr)\right).$
\end{cor}

In particular, taking $\theta=N-n-1$ yields the bound
\[
\Ord\!\left(
 n\D(G)
 + \frac{n^3}{N-n}
 + n^2\log\bigl(1+\min\{n,N-n\}\bigr)
\right)
\]
on the solution length, and $\Ord(|E(G)|)$ plus the same bound on the
computational time.

Theorem\tume\ref{refinement3} further gives, for every pair of
configurations on a feasible puzzle,
\[
\Ord\!\left(
 n\D(G)+\frac{n}{N-n}\sum_{v\in V(G,N-n-2)}\CL(G,v)
 +n^2\log\bigl(1+\min\{n,N-n\}\bigr)
\right)
\]
as an upper bound on the length of an optimal solution sequence.
Here $V(G,q)$ denotes the set of degree-two vertices of $G$ that occur
in a contiguous sequence of at least $q$ degree-two vertices.
In particular, cubic complexity is possible only in the
constant-blank regime and requires linearly many vertices in
$V(G,N-n-2)$ whose shortest cycles have linear length; the precise
statement is given in Theorem\tume\ref{Restriction} at the end of
Section~\ref{refinement}.

\subsection{Presentation and structural restrictions}

The principal tool for analyzing the constant-blank regime is the
\textit{Pebble Presentation Algorithm} developed in
Section~\ref{sec:pebble-presentation}.
Theorem\tume\ref{thm:pebble-presentation}, the
\textit{Pebble Presentation Theorem}, states that if
$2\leq N-n=\Ord(1)$, then for every feasible puzzle $(G,P)$ on a
connected $N$-vertex board graph, every configuration $\sigma$ of $P$ on
$G$, and every vertex $u\in V(G)$, there exists a move sequence of length
$\Ord(N^2)$ in which every pebble is presented at
$u$ exactly once; moreover, such a sequence can be constructed in
$\Ord(N^2)$ time.  This presentation order provides an address system
that can be transcribed across cut vertices.

The first major consequence is
Theorem\tume\ref{thm:intermediate-cut-component}.  If the board graph
$G$ of a feasible puzzle is written as $G=X\cup Y$, where $X$ and $Y$
are connected and $V(X)\cap V(Y)=\{u\}$ for a cut vertex $u$, then
$|V(X)|=\omega(1)\cap\ord(N)$ implies that the worst-case shortest
solution length on $G$ is $\ord(N^3)$.  Thus a side of intermediate
order separated by a single cut vertex is incompatible with cubic
complexity.

Applying this transcription mechanism to the
biconnected--bridge-tree decomposition $\mathcal{B}(G)$ determines the
possible scales of its components.  Here and below, the
reconfiguration cost of a component is understood in the sense defined
in Section~\ref{sec:pebble-presentation}.
Theorem\tume\ref{thm:no-intermediate-decomposition-component} states
that, for a family with worst-case shortest solution length
$\Theta(N^3)$, there are constants $K<\infty$ and $\delta>0$ such
that every component of $\mathcal{B}(G)$ has either at most $K$
vertices or at least $\delta N$ vertices, for all sufficiently large
$N$.  In particular, no component can have order
$\omega(1)\cap\ord(N)$.
Theorem\tume\ref{thm:linear-component-cubic} strengthens this
conclusion: every component $X\in\mathcal{B}(G)$ of order
$\Theta(N)$ must itself have reconfiguration cost
$\Theta(|V(X)|^3)=\Theta(N^3)$.  Since a tree of order $\Theta(N)$ has
subcubic reconfiguration cost by Theorem\tume\ref{main},
Corollary\tume\ref{cor:bridge-tree-components-bounded} follows: every
bridge-tree component of $\mathcal{B}(G)$ has order $\Ord(1)$.

Theorem\tume\ref{thm:linear-ear-closure} establishes a strong closure
principle: once a suitable subcubic core is present, arbitrary successive
ear extensions cannot create cubic complexity, provided that the resulting
graph is biconnected; neither the number nor the individual lengths of the
ears are restricted.  This result will be used in
Section~\ref{lower bound} to exclude additional internal structures from
the cubic regime.
%The final closure result of Section~\ref{sec:structural-restrictions}
%concerns ear extensions.
%Theorem\tume\ref{thm:linear-ear-closure} shows that, under its stated
%hypotheses, subcubic reconfiguration complexity is preserved under an
%arbitrary sequence of ear additions leading to a biconnected graph,
%with no restriction on the number or the individual lengths of the ears.
%This result will be used in Section~\ref{lower bound} to exclude
%additional internal structures from the cubic regime.

Combining these reductions with the cycle analysis of
Section~\ref{lower bound} yields our strongest structural result.

\noindent
\textbf{Theorem~\ref{thm-Theta-n3-unique-large-component}
(structural restrictions for cubic complexity).}
For every family of feasible puzzles on $N$-vertex board graphs whose
worst-case shortest solution length is $\Theta(N^3)$, the decomposition
$\mathcal{B}(G)$ has the following properties for all sufficiently
large $N$:
\begin{enumerate}
\item
there exists a biconnected component $X$ such that
$|V(X)|=\Theta(N)$ and the reconfiguration cost of $X$ is
$\Theta(N^3)$;
\item
$X$ contains no cycle $C$ with
$|V(C)|=\omega(1)\cap\ord(N)$;
\item
if $C_1,C_2\subseteq X$ are cycles of order $\Theta(N)$, then
$|V(C_1)\cap V(C_2)|\geq N-n$ or
$V(C_1)\cap V(C_2)=\emptyset$;
\item
every component of $\mathcal{B}(G)$ other than $X$ has order
$\Ord(1)$;
\item
for every cut vertex $v\in V(X)$, the union of all components of
$\mathcal{B}(G)$ attached to $X$ through $v$ has order
$\Ord(1)$.
\end{enumerate}

In particular, $X$ is the unique component of order $\Theta(N)$ in
$\mathcal{B}(G)$, so all unbounded structure responsible for cubic
complexity is confined to a single biconnected core.  Moreover, this
core contains no cycle of unbounded sublinear order, and any two
intersecting cycles of order $\Theta(N)$ must share at least $N-n$
vertices.

\subsection{Sharpness for the cubic regime}

For brevity, we call a family of feasible puzzles a
$\Theta(N^3)$-puzzle family if its worst-case shortest solution length
is $\Theta(N^3)$; when the family is clear, we suppress the family
index and simply write $\Theta(N^3)$-puzzle.  Section~\ref{lower bound}
defines an \textit{edge-amalgamation} onto a cycle by replacing each cycle
edge with a connected graph attached through a designated pair of vertices.

Section~\ref{lower bound} shows that the preceding restrictions are
substantive and identifies sharp separating examples.
Theorem\tume\ref{thm-Theta-n3} proves that cubic complexity genuinely
occurs in a nearly cyclic structure: if $G$ is an edge-amalgamation
of a sequence of connected graphs onto a cycle $C_0$, $(G,P)$ is
feasible, and $|V(G)|-|V(C_0)|=\Ord(1)$, then $(G,P)$ is a
$\Theta(N^3)$-puzzle.  Thus a cycle containing all but
$\Ord(1)$ vertices can retain the full cubic complexity.

Theorem\tume\ref{thm-cycle-o-n3} shows a complementary threshold
within the same family.  If $|V(C_0)|=\Theta(N)$, the total number of
vertices outside $C_0$ is $\omega(1)$, and every individual
attachment contributes only $\Ord(1)$ vertices outside $C_0$, then
the resulting feasible puzzle is an $\ord(N^3)$-puzzle.

A different separation occurs when two large regions meet at a single
vertex.  Theorem\tume\ref{thm-two-cycles-o-n3} considers board graphs
$H_1,H_2$ of feasible puzzles with $N_i=|V(H_i)|$ and
$N_1\leq N_2$, identified at one vertex, and gives the upper bound
$\Ord(N^3/N_1+N^2\log N_1)$ for the worst-case shortest solution
length.  Corollary\tume\ref{cor-two-cycles-o-n3} therefore gives an
$\ord(N^3)$ bound whenever $N_1=\omega(1)$.  In the special case
where $H_1$ and $H_2$ are cycles of order $\Theta(N)$,
Theorem\tume\ref{thm-two-cycles-Theta-n2logn} determines the exact
order $\Theta(N^2\log N)$.  Thus even two linear-size cycles meeting
at a single cut vertex do not produce cubic complexity.

Together, these results motivate Conjecture\tume\ref{conj-Theta-n3}:
if $(G,P)$ is a $\Theta(N^3)$-puzzle, then $G$ is an
edge-amalgamation onto a cycle $C_0$ with
$|V(G)|-|V(C_0)|=\Ord(1)$.

\section{Solving the Problem on Trees with Ample Blank Space}\label{base-cases}

%%%nakami

This section describes the base cases, where the size of the board $N$ is significantly larger than the number of pebbles $n$, which are the cornerstone of the proof of Theorem\tume\ref{main}.

%In the base case where $n \ll N$, we will demonstrate that for given configurations $f, g \in (G, P)$, $f$ can be transformed into $g$ with $\Ord (nN + n^2 \log n)$ moves, and that such a sequence of moves can also be computed with the same time complexity $\Ord (nN + n^2 \log n)$. 

The main result of this section is the following theorem.

\begin{thm}\label{sbs3}
Let $(G,P)$ be a feasible puzzle on a tree, and let $f$ and $g$ be
configurations of $P$ on $G$.
If $N\ge 3n$, then a solution sequence from $f$ to $g$ of length
\[
\Ord\bigl(n\D(G)+n^2\log n\bigr)
\]
can be computed in time
\[
\Ord\bigl(N+n\D(G)+n^2\log n\bigr).
\]
\end{thm}

We will show the outline of the proof of Theorem\tume\ref{sbs3}.
Let $\cF(G, P)$ be the set of all configurations on $(G, P)$.
\begin{itemize}
\item Under the condition $n \ll N$, there exists a subtree $H$ of $G$ such that the puzzle $(H, P)$ is feasible and $(H, P)$ can be divided into a family of sub-puzzles on $H$.  
\item For given configurations $f, g \in \cF(G, P)$, choose a pair of configurations $f_H, g_H \in \cF(H, P)$ such that $f$ and $g$ can be transformed to $f_H$ and $g_H$ respectively with $\Ord (n\D(G))$ moves.
\item In order to transfer $f_H$ to $g_H$ with $\Ord (n\D(G) + n^2 \log n)$ moves, we apply binary comparisons, where the puzzle $(H, P)$ is divided into sub-puzzles, in which each sub-puzzle has a half number of original pebbles, based on a divide-conquer strategy. 
\item The algorithm works similarly to the merge sorting algorithm:
 it splits the puzzle into two halves until each sub-puzzle has one pebble, and merges the pebbles to build a target configuration.  
\end{itemize}

%We propose an algorithm to move unlabeled pebbles from a set of vertices to another. 
%For a vertex $v \in V(G)$, let us denote the set of neighbors of $v$ by $N(v)$.
Firstly, we state a basic lemma concerning unlabeled pebble motion problem on trees.

\begin{la}\label{sbs0}
Let $S_1, S_2 \subset V(G)$ with $|S_1| = |S_2| = n$.
Then there exists a sequence of moves of length $\Ord (n\D(G))$ which transfers all the pebbles from $S_1$ to $S_2$, and it can be computed in time $\Ord (N + n\D(G))$.
\end{la}

Although Lemma\tume\ref{sbs0} can be derived from the proof of Theorem 2.1 of  \cite{WOS:000254460900008}, we will show a proof in the Appendix for the convenience of the reader.
Lemma\tume\ref{sbs0} is used in subsequent lemmas concerning labeled pebbles.

%Let $T$ be a tree with $N$ vertices.
%Let $P$ be a set of $n$ pebbles.
Let us denote $N + n\D(G)$ by $\alpha(G, P)$.
For a configuration $\pi \in \cF(G,P)$, the {\it support} of $\pi$ is defined as $\pi^{-1}(P)$, which is the subset of vertices of $G$ occupied by $P$.
We denote the support of $\pi$ by ${\rm sup}(\pi)$.   
For two configurations $\pi, \pi^{\ast} \in \cF(G,P)$, we say that $\pi$ and $\pi^{\ast}$ are {\it strongly related}, if there exists an algorithm which transfers $\pi$ to $\pi^{\ast}$ with a time complexity of $\Ord (\alpha(G,P))$.
For any configuration $\pi$ and for any $S \subset V(G)$ with $|S| = n$, by applying Lemma\tume\ref{sbs0}, we have a configuration $\pi^{\ast}$ such that ${\rm sup}(\pi^{\ast})$ is $S$ and $\pi^{\ast}$ is strongly related to $\pi$. 

To prove Theorem\tume\ref{sbs3}, the next two cases are particularly important.
The first key case is where $G$ has a rivet, a vertex of degree at least $3$, dividing $G$ into two parts with a suitable size.

\begin{la}\label{sbs1}
Let $G$ be a tree having a vertex $o$ with $\deg_G o \ge 3$ and let $G$ have subtrees $G_1$ and $G_2$ such that $V(G_1) \cap V(G_2) = \{ o \}$, and $|V(G_i)| \ge n + 1$ for $i=1, 2$. 
%Then the pebble motion problem on $G$ is computted in time $\Ord (Nn + n^2 \log n)$.
Then a solution sequence of length $\Ord (n\D(G) + n^2 \log n)$ for the \textsc{Pebble Motion Problem on Trees} can be computed in time $\Ord (N + n\D(G) + n^2 \log n)$.
\end{la}

The second key case for the proof of Theorem\tume\ref{sbs3} is where $G$ has a long isthmus with respect to the number of pebbles.

\begin{la}\label{sbs2}
Let $G$ be a tree with an isthmus $I$ of size at least $n$.
Let $G_1$ and $G_2$ be the two subtrees induced by $V(G) \setminus V_0$, where $V_0$ is a set of inner vertices of $I$.
%If $|V(G_i)| \ge n/2$ for $i = 1,2$, then the pebble motion problem on $G$ can be computed in time $\Ord (Nn + n^2 \log n)$.
If $|V(G_i)| \ge n/2 + 1$ for $i = 1,2$, then a solution sequence of length $\Ord (n\D(G) + n^2 \log n)$ for the \textsc{Pebble Motion Problem on Trees} can be computed in time $\Ord (N + n\D(G) + n^2 \log n)$.
\end{la}

Combining Lemma\tume\ref{sbs0}, Lemma\tume\ref{sbs1} and Lemma\tume\ref{sbs2}, we can prove Theorem\tume\ref{sbs3}. 

Before closing the section, we note the following lemma which states that starting from a given configuration, we can exchange the positions of two specified pebbles in time $\Ord (N + n\D(G))$.

Let $G$ be a tree.
We assume that the puzzle $(G, P)$ is feasible.
For a configuration $\pi \in \cF(G, P)$ and two pebbles $p, q \in P$, a configuration $\pi^\ast = \pi(p\leftrightarrow q)$ is defined as $\pi^\ast \circ \pi^{-1} (p) = q$, $\pi^\ast \circ \pi^{-1} (q) = p$ and $\pi^\ast \circ \pi^{-1} (r) = r$ for all $r \in P \setminus \{ p, q \}$.

\begin{la}\label{sbs4}
If the puzzle $(G, P)$ satisfies $N \ge 3n$, then for any $\pi \in \cF(G, P)$ and $p, q \in P$, there exists a sequence of moves of length $\Ord (n\D(G))$ which transforms $\pi$ to $\pi(p\leftrightarrow q)$, and it can be calculated in time $\Ord (N + n\D(G))$.
\end{la}

Lemma\tume\ref{sbs4} will be used in Section\tume\ref{refinement-tree}. 
 
%%%

\section{Decomposing a Tree into Smaller Subtrees}\label{splitting-tree}

To prove the main theorem, we need to decompose an instance puzzle 
into smaller sub-puzzles that satisfy the conditions of the base cases 
in the previous section. In this section we develop several structural
tools for this purpose.

The following three lemmas play complementary roles in our algorithm.
Lemma\tume\ref{centroid-splitting} provides a balanced way to split a tree
into two subtrees while sharing at most one vertex.
Lemma\tume\ref{feasible-subtree} allows us to construct subtrees of
prescribed size while preserving the bound on the maximum isthmus.
Finally, Lemma\tume\ref{bunkatsu} repeatedly applies 
Lemma\tume\ref{centroid-splitting} to obtain a global decomposition of a
tree into many smaller subtrees of bounded size, which will serve as the
domains of the subproblems solved by Theorem\tume\ref{sbs3}.

We first recall a standard property of tree centroids.
Every tree $T$ has either a vertex $v$ (called a {\it centroid vertex} of $T$)
such that the order of each connected component of $T-v$ does not exceed $(|V(T)|-1)/2$,
or an edge $e$ (called a {\it centroid edge} of $T$) such that each of the two connected
components of $T-e$ has order exactly $|V(T)|/2$.
\begin{la}[\hspace{-0.3pt}\cite{MR396308}]\label{centroid}
A centroid vertex or a centroid edge of a tree can be found in linear time.
\end{la}
Using this fact, we derive the following lemma.
\begin{la}\label{centroid-splitting}
Let $T$ be an arbitrary tree with $N$ vertices. 
Then one can find, in $\Ord(N)$ time, a pair of subtrees $(T_1, T_2)$ of $T$ such that all of 
the following three conditions hold:
\begin{enumerate}
    \item $|V(T_1)\cap V(T_2)| \leq 1$;
    \item $V(T_1) \cup V(T_2) = V(T)$;
    \item $(N+2)/3 \leq |V(T_1)| \leq |V(T_2)| \leq (N+1)/2$.
\end{enumerate}
\end{la}
This Lemma\tume\ref{centroid-splitting} is used when recursively splitting a given tree 
into two subtrees as balanced as possible, while sharing at most one vertex.

%%%%%%%%%%%%%%%%%%%%%%%%%%%%%%%%%%%%%%%%%%%%%%%%
%%%%%%%%%%%%%%%%%%%%%%%%%%%%%%%%%%%%%%%%%%%%%%%%

Consider a tree $X$ (the board tree) serving as the input board graph of a PMT. 
If the order $N$ of a tree $X$ (the board tree) is considerably larger than the number of pebbles $n = |P|$ 
to be placed on it, solving the PMT on an appropriate subtree $Y$ of $X$ can be more efficient---both 
in computation time and in the number of moves---than solving it directly on $X$. Here, $Y$ is chosen 
so that its order is at least $n+k+1$, where $k:=k(X)$ is the maximum
isthmus size of $X$, and the maximum isthmus size of $Y$ is at most $k$.
The solution of the PMT on $(Y, P)$ can 
then be transferred back to the original instance $(X, P)$. 
The following Lemma\tume\ref{feasible-subtree} is used to find such a suitable $Y$ within $X$. 
In our proposed algorithm, situations repeatedly arise in which one must find a subtree $Z$ of $X$ 
that contains a specified subtree $W$, has a prescribed order, and does not contain any isthmus 
larger than the maximum isthmus of $X$, and then solve the PMT on $Z$. In these cases as well, 
this lemma is used extensively. 
\begin{la}\label{feasible-subtree}
Let $T$ be a tree such that the size of a maximum isthmus of $T$ is at most $k$. 
For any subtree $T^{\prime}$ of $T$ such that the size of a maximum isthmus of $T^{\prime}$ is also at most $k$, the following statements hold:
\begin{itemize}
\item For any natural number $d$ satisfying $|V(T^{\prime})| \leq d \leq |V(T)|$, 
we can find, in $\Ord(d)$, a $d$-vertex subtree $T^{\prime\prime}$ of $T$ such that the size of a maximum isthmus 
of $T^{\prime\prime}$ is at most $k$ and that $T^{\prime\prime}$ contains $T^{\prime}$ 
as a subgraph.
\end{itemize}
\end{la}

%%PMTの入力board graphである木（board tree）$X$の位数$N$が、$X$の上に配置するpebble達$P$の個数$n=|P|$に比べて小さいときには、
%%$X$の上のPMTを直接扱うよりも、$X$の最大isthmusのサイズ$k$と$n$の和でバウンドされた$n+k+1$の位数を持つ$X$の適切なsub tree$Y$で、
%%再び$Y$の最大isthmusのサイズも$k$を超えないものを見つけて、$(Y,P)$というインスタンスに対するPMTを解き、それを$(X,P)$のインスタンスに
%%戻してやる方が計算効率も上がるし、moveの手数も抑えることが出来る。以下の補題\ref{feasible-subtree}は、$X$からそのようなうまい$Y$を
%%見つける際に用いられる。また、我々が提案するアルゴリズムにおいては、指定された$X$の部分木$Z$を含み、かつ指定された位数を持ち、
%%更に、$X$の最大isthmusのサイズを越える最大isthmusを持たないような部分木$W$を見つけて、$W$の上のPMTを解かなければならない
%%シチュエーションが繰り返される。そうした際にも、この補題は頻繁に利用される。
%\begin{la}\label{feasible-subtree}
%Let $T$ be a tree such that the size of a maximum isthmus of $T$ is at most $k$. 
%For any subtree $T^{\prime}$ of $T$ such that the size of a maximum isthmus of $T^{\prime}$ is also at most $k$, the following statements hold.
%%最大isthmus長が$k$以下である任意の木$T$の任意の部分木$T^{\prime}$について、以下の言明が成り立つ。
%\begin{itemize}
%\item{For any natural number $d$ satisfying $|V(T^{\prime})|\leq d \leq |V(T)|$, 
%we can find, in $\Ord(d)$, a $d$-vertex subtree $T^{\prime\prime}$ of $T$ such that the size of a maximum isthmus 
%of $T^{\prime\prime}$ is at most $k$ and that $T^{\prime\prime}$ contains $T^{\prime}$ 
%as a subgraph.}
%%$T^{\prime}$の最大isthmus長も$k$以下であれば、$|V(T^{\prime})|\leq d \leq |V(T)|$を満たす任意の自然数$d$に対して、最大isthmus長が$k$以下であり$T^{\prime}$を部分グラフとして含む$T$の位数$d$の部分木が線形時間の手間で見つかる。}
%\end{itemize}
%\end{la}

%%%%%%%%%%%%%%%%%%%%%%%%%%%%%%%%%%%%%%%%%%%%%%%%
%%%%%%%%%%%%%%%%%%%%%%%%%%%%%%%%%%%%%%%%%%%%%%%%
%%%%%%%%%%%%%%%%%%%%%%%%%%%%%%%%%%%%%%%%%%%%%%%%
%%%%%%%%%%%%%%%%%%%%%%%%%%%%%%%%%%%%%%%%%%%%%%%%

By repeatedly applying Lemma\tume\ref{centroid-splitting}, we obtain the
following structural decomposition.
\begin{la}\label{bunkatsu}
For any tree $T$ with $n$ vertices and any positive integer $k\le n$,
there exists an algorithm that constructs a family $\mathcal{F}$ of at most
$4n/k$ subtrees such that
\[
\bigcup_{H\in\mathcal{F}}V(H)=V(T),
\]
any two subtrees in $\mathcal{F}$ share at most one vertex, and every
$H\in\mathcal{F}$ satisfies
\[
\frac{k}{3}<|V(H)|\le k.
\]
\end{la}

\noindent
Lemma\tume\ref{bunkatsu} provides a decomposition of the board tree
into subtrees of sizes to which Theorem\tume\ref{sbs3} applies.
Our algorithm repeatedly solves PMT subproblems associated with these 
subtrees.
In the algorithm of Lemma\tume\ref{bunkatsu}, the subtrees covering $T$ are
stored and maintained in a binary tree data structure $\BT$.
Each node $x$ of $\BT$ corresponds to a subtree $X(x)$ of $T$ that
arises during the recursive decomposition.
When a subtree is split during the recursive decomposition,
the resulting pieces may share either a vertex or an edge.
We refer to this shared object as the \textit{joint} of the subtree.
Accordingly, the joint is represented by either a vertex
(a \textit{joint vertex}) or an edge (a \textit{joint edge}).
For each node $x$ of $\BT$, let $p(x)$ denote its parent node and
$c_i(x)$ $(i=1,2)$ denote its child nodes.
Each node $x$ stores the following three objects as its attributes:
\begin{itemize}
\item $X(x)$,
\item $\Cen(X(x)) := (\CenV(X(x)), \CenE(X(x)))$,
\item $\Joi(X(x)) := (\JoiV(X(x)), \JoiE(X(x)))$.
\end{itemize}
First, $X(x)$ is a subtree of $T$ determined by the inductive
procedure described below.
For a subtree $X(x)$, the pair
$\Cen(X(x)) = (\CenV(X(x)),\CenE(X(x)))$ stores its centroid
information.
If $X(x)$ has a centroid vertex, it is stored in $\CenV(X(x))$;
otherwise $\CenV(X(x)):=\emptyset$.
Similarly, if $X(x)$ has a centroid edge, it is stored in
$\CenE(X(x))$; otherwise $\CenE(X(x)):=\emptyset$.
If $x$ is a leaf node of $\BT$, we define
$\Cen(X(x)):=(\emptyset,\emptyset)$ to avoid unnecessary computation.
We now describe how the objects associated with the nodes of $\BT$
are defined.
\begin{alg-enumerate}
\item
For the root node $r$ of $\BT$, let
$X(r):=T$ and $\Joi(X(r)):=\Cen(X(r))$.
\item
Suppose that $X(x)$ is stored at a non-leaf node $x$ of $\BT$.
The objects $X(c_i(x))$ and $\Joi(X(c_i(x)))$ $(i=1,2)$ to be stored
in the child nodes $c_i(x)$ are defined inductively by the following
procedure:
\begin{alg-enumerate}
\item
Using Lemma\tume\ref{centroid-splitting}, find two subtrees
$X_1$ and $X_2$ of $X(x)$ satisfying
$V(X_1) \cap V(X_2)=\{\CenV(X(x))\}$ and
$V(X_1) \cup V(X_2)=V(X(x))$.
If $X(x)$ has a centroid edge, $X_1$ and $X_2$ are the two subtrees
obtained by removing that edge from $X(x)$.
\item
Of the two subtrees $X_1$ and $X_2$ of $X(x)$,
the one that contains $\JoiV(X(x))$ or one of the endpoints
of $\JoiE(X(x))$ is stored at the second child node $c_2(x)$ of $x$.
The other subtree is stored at the first child node $c_1(x)$.
\label{No2Child}
\item
$\Joi(X(c_2(x))):=\Joi(X(x))$.
\item
$\Joi(X(c_1(x))):=\Cen(X(x))$.
\end{alg-enumerate}
\end{alg-enumerate}

\noindent\textit{Re-rooting.}
The assignment of the first and second children in Steps~2-b--2-d depends on the joint of the parent and hence on the current root. Whenever suppression replaces the root of $\BT$ (or of a copy $\BT2$) by a surviving child $x$, we reset
\[
    \Joi(X(x)):=\Cen(X(x)).
\]
We then reorient the surviving binary tree rooted at $x$ from the root downward. At each non-leaf node, the two stored child subtrees themselves are left unchanged, but the child containing the current joint is designated as the second child and the other as the first child, interchanging the two child pointers when necessary; their joints are then reassigned according to Steps~2-c and~2-d. Thus re-rooting changes neither the stored subtrees $X(\cdot)$ nor the underlying decomposition, but only the first/second-child designation and the joint data.
If both children satisfy the criterion in Step~2-b, as happens at the root where $\Joi(X(r))=\Cen(X(r))$, either designation may be used. This choice does not affect Lemma~\ref{children}, which holds for both children of the root.

Strictly speaking, symbols such as $p^{\BT}(x)$, $c^{\BT}_i(x)$,
and $X^{\BT}(x)$ should be used to indicate that these functions are
defined on the nodes of the binary tree $\BT$.
However, to avoid notational clutter, we omit such superscripts
throughout the paper whenever no confusion arises, even when
multiple binary trees (such as $\BT_1$, $\BT_2$, etc.) are considered
simultaneously.

Then, the following lemmas hold. 
%このとき、以下の事実が成り立ちます。
\begin{la}\label{children}
For any child node $c_i(r) (i=1,2)$ of the root node $r$ of this binary tree $\BT$, the following statement holds:
%二分木$\BT$のrootノード$r$のどの子ノード$c_i(r) (i=1,2)$についても、以下の言明が成り立つ：
\begin{itemize}
\item{Regarding the leaf node $\ell(i)$ that can be reached by successively selecting only the first child of each descendant of the node $c_i(r)$, even if we delete the part $X(\ell(i))-\JoiV(X(\ell(i)))$ from the graph $X(r)$, the resulting graph $X(r)-(X(\ell(i))-\JoiV(X(\ell(i))))$ is connected. 
%$X(\ell(i))-\JoiV(X(\ell(i)))$を$X(r)$から削除して得られるグラフ$X(r)-(X(\ell(i))-\JoiV(X(\ell(i))))$は連結である。
}
\end{itemize}
\end{la}

Whenever the root is replaced during $\Alg.\CutOff$ or $\Alg.\Extraction$, the re-rooting operation above restores $\Joi(X(r))=\Cen(X(r))$ for the current root $r$ and reassigns the first and second children consistently with the current joints. Hence Lemma\tume\ref{children} applies to the current root of $\BT$ and of $\BT2$ at every stage of these algorithms.

\begin{la}\label{BT}
The effort required to create the binary tree $\BT$ from a given tree $T$ with $n$-vertices using the algorithm from Lemma\tume\ref{bunkatsu} is $\Od(n\log n)$. 
%与えられた$n$頂点の木$T$からLemma~\ref{bunkatsu}のアルゴリズムを用いて二分木リスト$\BT$を作るのに必要な手間は$\Od(n\log n)$である。
\end{la}

It is worth noting that Lemma\tume\ref{children} is essential for ensuring the correctness of the Main Algorithm, 
which will be presented in the next section. 
Lemma\tume\ref{BT}, on the other hand, provides a complexity bound for constructing the binary tree database 
used in the algorithm's complexity analysis.

\section{The Global Tree Algorithm}\label{MainAlg}

In this section, we present the main algorithm for solving the \textsc{Pebble Motion Problem on Trees} in $\Ord(nN+n^2\log(1+\min\{n,k\}))$ 
without any additional constraints on the instances.

%%%%%%%%%%%%%%%%%%%%%%%%%%%%%%%%%%%%%%%%%%%%%%%%%%%%%%%%%%
%%%%%%%%%%%%%%%%%%%%%%%%%%%%%%%%%%%%%%%%%%%%%%%%%%%%%%%%%%
%%%%%%%%%%%%%%%%%%%%%%%%%%%%%%%%%%%%%%%%%%%%%%%%%%%%%%%%%%
%%%%%%%%%%%%%%%%%%%%%%%%%%%%%%%%%%%%%%%%%%%%%%%%%%%%%%%%%%
%%%%%%%%%%%%%%%%%%%%%%%%%%%%%%%%%%%%%%%%%%%%%%%%%%%%%%%%%%
%%%%%%%%%%%%%%%%%%%%%%%%%%%%%%%%%%%%%%%%%%%%%%%%%%%%%%%%%%
%%%%%%%%%%%%%%%%%%%%%%%%%%%%%%%%%%%%%%%%%%%%%%%%%%%%%%%%%%

In the previous section, we established the recursive framework
outlined in the Introduction.
Its key idea is that Theorem~\ref{sbs3} provides an efficient
algorithm whenever the computation is performed inside a connected
workspace containing sufficiently many empty vertices.
The purpose of this section is to explain how such workspaces are
constructed, relocated, and preserved throughout the recursion.

Our base algorithm works only when a connected workspace containing
sufficiently many empty vertices is available.
The main contribution of the recursive algorithm therefore lies not
only in recursively decomposing the tree into smaller subproblems,
but also in continually creating, relocating, and preserving such
workspaces while progressively shrinking the remaining search space.
The procedures introduced in this section are precisely the mechanisms
that maintain the workspace invariant introduced in the Introduction.

\subsection{Roadmap of the Main Algorithm}

The recursive algorithm consists of four procedures organized into
two conceptual layers.

The preprocessing layer constructs a compact target configuration
before the recursive algorithm begins.

$\Alg.\PacCon$
(\textsc{PacCon})
transforms the target configuration into a compact configuration on
the reduced subtree $T$ while recording the corresponding sequence
of moves.

$\Alg.\Pac(T,\BT,r)$
(\textsc{Pac})
is a recursive subroutine used internally by
\textsc{PacCon}.
It recursively packs pebbles according to the decomposition stored
in $\BT$, thereby producing the compact configuration required by
the recursive algorithm.

The recursive layer repeatedly completes one target subtree after
another while maintaining the workspace invariant introduced in the
Introduction.

$\Alg.\CutOff$
(\textsc{CutOff})
is the top-level recursive procedure.
It repeatedly reduces the remaining instance until it becomes
sufficiently small for direct solution by
Theorem~\ref{sbs3}.

$\Alg.\Extraction(T,\BT,g',t)$
(\textsc{Extraction})
is an auxiliary procedure called internally by
\textsc{CutOff}.
It assembles the workspace required for each recursive step by
gathering the necessary pebbles into a temporary subtree where
Theorem~\ref{sbs3} can be applied.

Conceptually,
\textsc{CutOff} drives the overall recursion.
Whenever a workspace must be assembled,
it invokes \textsc{Extraction},
which in turn relies on
\textsc{PacCon} and \textsc{Pac}
to organize temporary packings.

The procedures are presented in the next subsection according to
their logical dependencies rather than their execution order.
This organization allows each procedure to be introduced only after
all lower-level mechanisms on which it depends have already been
established.

Together, these four procedures constitute the recursive framework
that maintains the workspace invariant introduced in the
Introduction and underlies the proof of
Theorem~\ref{main}.

%%%%%%%%%%%%%%%%%%%%%%%%%%%%%%%%%%%%%%%%%%%%%%%%%%%%%%%%%%
%%%%%%%%%%%%%%%%%%%%%%%%%%%%%%%%%%%%%%%%%%%%%%%%%%%%%%%%%%
%%%%%%%%%%%%%%%%%%%%%%%%%%%%%%%%%%%%%%%%%%%%%%%%%%%%%%%%%%
%%%%%%%%%%%%%%%%%%%%%%%%%%%%%%%%%%%%%%%%%%%%%%%%%%%%%%%%%%
%%%%%%%%%%%%%%%%%%%%%%%%%%%%%%%%%%%%%%%%%%%%%%%%%%%%%%%%%%
%%%%%%%%%%%%%%%%%%%%%%%%%%%%%%%%%%%%%%%%%%%%%%%%%%%%%%%%%%
%%%%%%%%%%%%%%%%%%%%%%%%%%%%%%%%%%%%%%%%%%%%%%%%%%%%%%%%%%

\subsection{Recursive Procedures}

We now describe the four procedures introduced above.
Following the dependency structure explained in the roadmap,
they are presented in an order that follows their logical
dependencies rather than their execution order.

%%%%%%%%%%%%%%%%%%%%%%%%%%%%%%%%%%%%%%%%%%%%%%%%
%%%%%%%%%%%%%%%%%%%%%%%%%%%%%%%%%%%%%%%%%%%%%%%%
First, we state the following lemma, which handles the case where the
maximum isthmus size is bounded by a constant.

\begin{la}\label{constant}
Let $k:=k(G)$ be the maximum isthmus size of the board tree $G$.
If $k$ is bounded by a positive constant $M$, i.e.,
$k\leq M\in\Ord(1)$, then a solution sequence of length
$\Ord(n\D(G))\,(\subset \Ord(nN))$ to the
\textsc{Pebble Motion Problem on Trees} can be computed in time
$\Ord(N+n\D(G))\,(\subset \Ord(nN))$.
\end{la}

Therefore, in the remainder of this paper, we assume $k \gg 1$ for all subsequent discussions.

Then, the outline of the remaining part of the main algorithm is as follows:
%メインアルゴリズムの概要は次の通りである。
\begin{enumerate}
\item{For the input $(G,P)$, use Lemma\tume\ref{feasible-subtree} to find a subtree $T$ of $G$ such that it has order $n+k+1$ and the size of its maximum isthmus is $k$ or less. Then, move all pebbles in the initial configuration $f$ on $G$ onto $T$ with $\Ord(N+n\D(G))$ effort and $\Ord(n\D(G))$ moves, and call this configuration, which is strongly related to $f$, $f^{\prime}$.}
%, while sequentially recording the sequence of pebble movements that occur during this process.}
%\item{入力$(G,P)$のboard tree $G$に対して、位数が$n+k+1$であり、かつ、その最大isthmus長が$k$以下であるような$G$の部分木$T$をLemma\tume\ref{feasible-subtree}を用いて見つけ、$G$上の初期配置$f$におけるすべての碁石を$T$の上に$\Ord(Nn)$の手間で移動させる（そしてこの作業で生じたPPPの移動の系列を逐次記録しておく）。}
\item{Similarly, for the goal configuration $g$ on $G$, convert it in advance to a configuration $g^{\prime}$ on $T$ that is strongly related to $g$ and matches the operations from Step\tume\ref{kiriotoshi}.
%, while also sequentially recording the sequence of pebble movements that occur during this process. 
%同様にして、$G$上のゴール配置$g$についても、$g$とstrongly relatedであり、かつステップ\ref{kiriotoshi}の操作にマッチするような$T$上の配置$g^{\prime}$に予め変換しておく。
}
\item{By using the algorithm from Lemma\tume\ref{bunkatsu}, decompose the tree $T$ into subtrees of sufficiently small order with respect to $k$.
%Lemma\tume\ref{bunkatsu}の手法を用いて、$T$を、その最大isthmus長$k$に対して十分小さな部分木達に分解する。
}\label{bunkai}
\item{
Among the subtrees obtained from the procedure in Step\tume\ref{bunkai}, there is a tree $T^{\prime}$ such that removing all of $T^{\prime}$ except at most one vertex from $T$ keeps the remaining graph connected. Therefore, starting from the configuration $f^{\prime}$, repeat the following sequence of operations: find one such tree $T^{\prime}$, place pebbles on it without empty space according to the configuration $g^{\prime}$, and then cut off that portion from the main body. 
%上記\ref{bunkai}の手続きで得られた部分木達の中には、その高々$1$点を残して他のすべてを$T$から削除しても、残されたグラフが連結に保たれるような木が常に存在するので、そのような木$T^{\prime}$上に、配置$g^{\prime}$の通りに碁石を隙間なく配置してその部分をグラフから切り落とす操作を繰り返す。
}\label{kiriotoshi}
\end{enumerate}
In Step\tume\ref{kiriotoshi} of the outline above, in order to place pebbles according to the configuration $g^{\prime}$ on $T^{\prime}$, it is first necessary to gather the set of pebbles $g^{\prime}(V(T^{\prime}))\setminus\{0\}$ that are scattered throughout the entire graph $T$ near the vicinity of $T^{\prime}$. To achieve this, the following method will be used.
%なお、上記のステップ\tume\ref{kiriotoshi}において、「$T^{\prime}$上に配置$g^{\prime}$の通りに碁石を配置する」ためには、まず$T$全体に散在している$(g^{\prime})^{-1}(V(T^{\prime}))$という碁石の集合を、$T^{\prime}$付近に集めてくる必要がある。そのためには、下記の手法を用いる。
\begin{itemize}
\item{Using the method of Lemma\tume\ref{bunkatsu} to decompose $T$, there will always be a subtree $T^{\prime\prime}$, aside from $T^{\prime}$ such that removing all but at most one vertex from $T$ keeps the remaining graph connected. On such a subtree $T^{\prime\prime}$, tightly pack any pebbles that do not belong to $g^{\prime}(V(T^{\prime}))$ and then "temporarily" cut them off from the whole. Once the pebbles are arranged on $T^{\prime}$ according to the configuration $g^{\prime}$, all parts that were "temporarily" removed in this operation will be restored.
%Lemma\tume\ref{bunkatsu}の手法を用いて分解すると、$T$には、自身の高々$1$点以外のすべてを$T$から削除しても残されたグラフが連結に保たれるような部分木が、$T^{\prime}$以外にも、常に存在する。そうした部分木$T^{\prime\prime}$の上に$(g^{\prime})^{-1}(V(T^{\prime}))$以外の碁石を隙間なく詰め込んで、$T$から一旦切り落とす（その頂点集合を灰色に着色する）という操作を繰り返す。配置$g^{\prime}$の通りに$T^{\prime}$上に碁石が配置されたら、この操作でテンポラリに削除した部分はすべて復元する。
}
\end{itemize}

%%%%%%%%%%%%%%%%%%%%%%%%%%%%%%%%%%%%%%%%%%%%%%%%
%%%%%%%%%%%%%%%%%%%%%%%%%%%%%%%%%%%%%%%%%%%%%%%%
%%%%%%%%%%%%%%%%%%%%%%%%%%%%%%%%%%%%%%%%%%%%%%%%

In the following sections of this paper, we will assume that each vertex $v$ of the input tree $G$ (and consequently all of its subgraphs) is associated with the following two attributes.
%以降の議論においては、木$G,T$の各頂点$v$に次の2つの属性を持たせる。
\begin{namelist}{********}
\item[$\Peb(v)$]{A function that returns the pebble $x (\in \{1,\ldots,n\})$ or a blank $0$ placed on vertex $v$.}
%%頂点$v$に乗っている碁石$x (\in \{1,\ldots,n\})$もしくは空白$0$を返す関数
\item[$\col(v)$]{A function that returns the color of vertex $v$. The possible colors are white, gray, or black.
%頂点$v$の色を返す関数．頂点の色は、白色、灰色、黒色のいずれか．
Regarding the meaning of the color, white indicates a vertex that has been completely removed from the graph, gray signifies a vertex that is temporarily removed during the execution of the algorithm but will be restored later, and black represents a vertex that remains in the graph.
%なお、treeの各頂点$v$の色の意味についてであるが、白色はグラフから完全に削除された頂点であることを、灰色はアルゴリズムの遂行上一旦はグラフから削除されるがのちに復元される頂点であることを、黒色はグラフに残された頂点であることを示している。
}
\end{namelist}
In other words, the graphs are treated not just as simple graphs, but
also as objects carrying information about the current pebble
configuration and the current status of their vertices.

At every stage of the algorithm, let
\[
    T_{\mathrm{black}}
    :=
    T\bigl[
        \{v\in V(T)\mid \col(v)=\mathrm{black}\}
     \bigr]
\]
denote the subtree induced by the vertices that currently remain
active.

The graphs $G$ and $T$ themselves remain fixed throughout the
algorithm.  The subtrees $X(x)$ stored at the nodes $x$ of $\BT$
are retained as decomposition data and, after nodes of $\BT$ have
been deleted, are not required to coincide with the current
remaining subtree.  The latter is represented by
$T_{\mathrm{black}}$.  The only changes made to the stored
subtrees $X(x)$ are the explicit updates
\[
    X(w):=X(w)\cup uv
\]
prescribed below for leaf nodes $w$.
%In other words, the graphs are treated not just as simple graphs, but also as objects that also contain information about the configurations of pebbles on them.
%%すなわち、$G$や$T$は単なるグラフであるのみならず、その上の碁石の配置の情報も併せ持つオブジェクトとして扱います。
%In this algorithm, the shapes such as the graphs $G$ and $T$, as well as the graphs $X(x)$ stored in each node $x$ of $\BT$ created by the operations in Step\tume\ref{bunkai} of the outline above, do not change; however, the values of the functions $\Peb(v)$ and $\col(v)$ at each vertex $v$ of the graphs continuously change.
%%このアルゴリズムにおいては、$G$や$T$あるいは$\BT$の各ノード$x$に格納された$X(x)$などの形状は変化しないが、treeの各頂点$v$の持つ関数$\Peb(v), \col(v)$の値は刻々と変化していく。

%%%%%%%%%%%%%%%%%%%%%%%%%%%%%%%%%%%%%%%%%%%%%%%%
%%%%%%%%%%%%%%%%%%%%%%%%%%%%%%%%%%%%%%%%%%%%%%%%
%%%%%%%%%%%%%%%%%%%%%%%%%%%%%%%%%%%%%%%%%%%%%%%%

According to Theorem\tume{A}, if the size of the pebble set $P$ is $n$ and the
size of the maximum isthmus of a tree $T$ is at most $k$, then the
pebble motion problem on $(T,P)$ is feasible whenever $|V(T)| \geq n+k+1$.
Therefore, given an instance $(G,P)$, our algorithm first uses
Lemma\tume\ref{feasible-subtree} to find a subtree $T$ of $G$ whose order
is $n+k+1$ and whose maximum isthmus size is at most $k$.
Next, we transform the goal configuration $g$ on $G$ into a
``suitable'' configuration $g'$ on $T$ that is strongly related to $g$.
By Lemma\tume\ref{sbs0}, the support of $g$ can first be transferred
onto $T$ in $\Ord(n\D(G))\subseteq\Ord(Nn)$ moves.
The procedure for constructing $g'$ is given by the algorithm
$\Alg.\PacCon(G,T,\BT,g)$ below.
%According to Theorem A, if the size of the pebble set $P$ is $n$ and the size of the maximum isthmus of the graph $T$ is at most $k$, then if the order of $T$ is at most $n + k + 1$, it is guaranteed that the pebble motion problem on $(T, P)$ is feasible.
%%the size of the set of pebbles $P$が$n$である場合、Theorem Aによれば、
%%最大isthmus長が$k$以下でありさえすれば、board tree $T$の位数は、
%%高々$n+k+1$あれば、$(T,P)$の上のpebble motion problemはfeasibleとなる。
%Therefore, the algorithm presented in this paper first uses Lemma\tume\ref{feasible-subtree} to find a subtree $T$ of the tree $G$ from the input $(G, P)$ such that its order is $n + k + 1$ and its size of a maximum isthmus is at most $k$.
%%よって、本論文に示すアルゴリズムでは、まず入力$(G,P)$のboard tree $G$に対して、
%%位数が$n+k+1$であり、かつ、その最大isthmus長が$k$以下であるような$G$の
%%部分木$T$をLemma\tume\ref{feasible-subtree}を用いて見つける。
%Furthermore, we rearrange the goal configuration $g$ on $G$ into a "suitable" configuration $g^{\prime}$ on $T$ that is strongly related to $g$ with a time complexity and number of moves of at most $\Ord(Nn)$. The algorithm used for this purpose is $\Alg.\PacCon(\BT, g)$ shown below.
%%さらに、$(G,P)$のgoal configuration $g$を$\Ord(Nn)$の手数で（$g$とstrongly relatedな）$T$上の「適切な」配置$g^{\prime}$に変換する。
%%そのために利用するアルゴリズムが、以下に示す$\Alg.$詰込配置$(\BT, g)$である。
\begin{alg}[$\Alg.\PacCon(G, T, \BT, g)$]\hfill
\begin{namelist}{OUTPUT}
\item[\sl Input]
A tree $G$ and its subtree $T$, the binary tree $\BT$ storing the
subtrees covering $T$ obtained by Lemma\tume\ref{bunkatsu},
and a configuration $g$ on $(G,P)$.

\item[\sl Output]
A configuration $g'$ on $(T,P)$ and a sequence of moves
$\phi[r]\cdot\psi$ transforming $g$ into $g'$.
\end{namelist}

\begin{alg-enumerate}

\item
Set $\Peb := g$.
For each $u\in V(T)$ set $\col(u):=\mathrm{black}$,
and for each $v\in V(G)\setminus V(T)$ set $\col(v):=\mathrm{white}$.

\item
Rearrange the configuration $g$ on $(G,P)$ into a configuration on
$(T,P)$ that is strongly related to $g$, and record the performed
sequence of moves as $\psi$.

\item
Call $\Alg.\Pac(T,\BT,r)$ for the root node $r$ of $\BT$ 
obtaining the move sequence $\phi[r]$.

\item
Return the move sequence $\phi[r]\cdot\psi$ and the configuration
$g' := \phi[r]\cdot\psi(g)$ on $(T,P)$.

\end{alg-enumerate}
\end{alg}
%\begin{alg}[$\Alg.\PacCon(G, T, \BT, g)$]\hfill
%\begin{namelist}{OUTPUT}
%\item[\sl Input]{A tree $G$ and its subtree $T$, the binary tree $\BT$ that stores the set of subtrees covering the tree $T$ obtained from the algorithm stated in Lemma\tume\ref{bunkatsu}, and a configuration $g$ on the pair $(G,P)$.}
%\item[\sl Output]{A configuration $g^{\prime} $ on $(T,P)$ and a sequence of moves $\phi[r]\cdot\psi$ to rearrange $g$ into $g^{\prime}$.}
%\end{namelist}
%\begin{alg-enumerate}
%\item{Set $\Peb:=g$. Set $\col(u):={\textrm{black}}$ for each $u\in V(T)$. Set $\col(v):={\textrm{white}}$ for each $v\in V(G) \setminus V(T)$.}
%\item{Rearrange the configuration $g$ on $(G,P)$ to some configuration on $(T,P)$ which is strongly related to $g$, and remember the sequence of moves performed here as $\psi$.}
%%配置$g$のすべての碁石を$T$上に移動させると共に、ここで行った一連のmoveの系列を$\psi$として記憶せよ}
%\item{Call $\Alg.\Pac(T, \BT, r)$ for the root node $r$ of $\BT$. }
%%$\BT$のrootノード$r$に対して、サブルーチン$\Alg.$詰込$(\BT,r)$を呼ぶ}
%\item{Return the sequence of moves $\phi[r]\cdot\psi$ and the configuration $g^{\prime} := \phi[r] \cdot \psi(g)$ on $(T,P)$.}
%%moveの系列$\phi[r]\cdot\psi$を出力して停止する}
%\end{alg-enumerate}
%\end{alg}
%%%
The procedure $\Alg.\Pac$ recursively packs pebbles toward the first
child subtrees of $\BT$ so that the configuration becomes compact
with respect to the decomposition structure.
\begin{alg}[$\Alg.\Pac(T,\BT,x)$]\hfill
\begin{namelist}{OUTPUT}

\item[\sl Input]
A tree $T$, the binary tree $\BT$ storing the set of subtrees covering
$T$ obtained from Lemma\tume\ref{bunkatsu},
and a node $x$ of $\BT$.

\item[\sl Output]
A sequence of moves $\phi[x]$.

\end{namelist}

\begin{alg-enumerate}

\item
If $x$ is a leaf node of $\BT$, return the empty sequence.

If there is no pebble on $X(x)$, or
if every vertex of $X(x)-\JoiV(X(x))$ is occupied,
return the empty sequence.

\item
Move pebbles into vertices of 
$X(c_1(x))-\JoiV(X(c_1(x)))$ 
until either
\begin{enumerate}
\item all vertices of $X(c_1(x))-\JoiV(X(c_1(x)))$ are occupied, or  
\item no pebble remains in $X(x)$ that can be moved there.
\end{enumerate}
Record the performed moves as $\phi[x]$.

\item
Compute $\phi[c_1(x)] := \Alg.\Pac(T,\BT,c_1(x))$
and $\phi[c_2(x)] := \Alg.\Pac(T,\BT,c_2(x))$.

\item
Update
\[
\phi[x] := \phi[c_2(x)]\cdot\phi[c_1(x)]\cdot\phi[x]
\]
and return $\phi[x]$.

\end{alg-enumerate}
\end{alg}

%\begin{alg}[$\Alg.\Pac(T, \BT, x)$]\hfill
%\begin{namelist}{OUTPUT}
%\item[\sl Input]{A tree $T$, the binary tree $\BT$ that stores the set of subtrees covering the tree $T$ obtained from the algorithm stated in Lemma\tume\ref{bunkatsu},
%and a node $x$ of $\BT$.}
%\item[\sl Output]{A sequence of moves $\phi[x]$.}
%\end{namelist}
%\begin{alg-enumerate}
%\item{If the node $x$ is a leaf node of $\BT$, else if there is no pebble on $X(x)$, else if there is no unoccupied vertex in $X(x)-\JoiV(X(x))$, then return a zero-length sequence of moves (i.e., an identity transformation) as $\phi[x]$.}
%%ノード$x$が$\BT$のleafノードであるか、$X(x)-\JoiV(X(x))$上に空白がないか、あるいは$X(x)$上に碁石が一つも乗っていない場合には、$\phi[x]$をmoveの長さ$0$の系列（何もしない恒等写像）として返して停止する
%\item{Fill the pebbles in the current configuration into the portion of the tree $X (c_1(x))$ except for $\JoiV(X(c_1(x)))$ % (i.e. $X(c_1(x))-\JoiV(X(c_1(x)))$) 
%as much as possible, and record the sequence of moves performed here as $\phi[x]$.}
%%現在の碁石の配置に対して、碁石を$\BT$のノード$x$の第1子$c_1(x)$に格納されている木$X(c_1(x))$の$\JoiV(X(c_1(x)))$以外の部分($X(c_1(x))-\JoiV(X(c_1(x)))$)に詰められるだけ詰め込むと共に、ここで行った一連のmoveの系列を$\phi[x]$として記憶する
%\item{Call $\Alg.\Pac(T, \BT, c_1(x))$ and $\Alg.\Pac(T, \BT, c_2(x))$.}
%%\item{$\Alg.$詰込$(\BT, c_1(x))$と$\Alg.$詰込$(\BT, c_2(x))$を呼ぶ}
%\item{Update $\phi[x]:=\phi[c_2(x)]\cdot\phi[c_1(x)]\cdot\phi[x]$ and return $\phi[x]$.}
%%\item{$\phi[x]:=\phi[c_2(x)]\cdot\phi[c_1(x)]\cdot\phi[x]$に更新して、$\phi[x]$を返して停止する}
%\end{alg-enumerate}
%\end{alg}
\noindent
In the main algorithm, the initial configuration $f$ is first rearranged to the configuration 
$g^{\prime} = \phi[r] \cdot \psi(g)$, which is the output of $\Alg.\PacCon(\BT, g)$.  
Then, by applying the inverse sequence of moves $\psi^{-1}\cdot\phi[r]^{-1}$ from $g^{\prime}$ to $g$, 
we obtain $g$ as $\psi^{-1} \cdot \phi[r]^{-1}(g^{\prime})$. 
%メインアルゴリズムでは、まず初期配置$f$を$g^{\prime}:=\phi[r]\cdot\psi(g)$に配置換えしてから、
%それを$\psi^{-1}\cdot\phi[r]^{-1}$というmoveの逆系列に掛けることにより、$g=\psi^{-1}\cdot\phi[r]^{-1}(g^{\prime})$を実現する。

%%%%%%%%%%%%%%%%%%%%%%%%%%%%%%%%%%%%%%%%%%%%%%%%
%%%%%%%%%%%%%%%%%%%%%%%%%%%%%%%%%%%%%%%%%%%%%%%%
We now turn to the top-level recursive procedure
\textsc{CutOff}, which repeatedly completes one target subtree
and shrinks the remaining instance.
After reducing the instance on $(G,P)$ to the one on $(T,P)$, we will use the following algorithm $\Alg.\CutOff(G,T,\BT,f,g^{\prime})$ to perform the following steps repeatedly: 1) find a subtree at the tip of tree $T$ that, when removed, does not cause the remaining parts to become disconnected (we will refer to this subtree as the {\it target} subtree), 2) fill the target subtree with pebbles according to the configuration $g^{\prime}$, and 3) cut that portion off from the whole. 
%$(T,P)$上のpebble motion problemに落とし込んだ後は、下記の$\Alg.$切り落とし$(\BT,f,h)$というアルゴリズムを用いて、$g^{\prime}$の配置通りに、$T$の端にあるsubtree（taget subtreeと呼ぶ）にpebbleを詰め込んで切り落とすという作業を繰り返す。
\begin{alg}[$\Alg.\CutOff(G,T,\BT,f,g^{\prime})$]\hfill
\begin{namelist}{OUTPUT}
\item[\sl Input]{An instance tree $G$ of the \textsc{Pebble Motion Problem on Trees}, 
its subtree $T$ of order $n+k+1$ whose maximum isthmus size is at most $k$, 
the binary tree $\BT$ storing the family of subtrees covering $T$ obtained in Lemma~\ref{bunkatsu}, 
a configuration $f$ on $(G,P)$, and an output configuration $g'$ on $(T,P)$ of $\Alg.\PacCon(G,T,\BT,g)$.}

\item[\sl Output]{A sequence of moves $\sigma$ that rearranges $f$ into $g'$.}
\end{namelist}

\begin{alg-enumerate}

\item
Set $\Peb := f$.  
Set $\col(u) := \mathrm{black}$ for each $u \in V(T)$, and  
$\col(v) := \mathrm{white}$ for each $v \in V(G)\setminus V(T)$.

\item
Rearrange the configuration $f$ on $G$ to a configuration on $T$ 
that is strongly related to $f$, and store the sequence of moves as $\sigma$.

\item
If
\[
    |V(T_{\mathrm{black}})|
    \leq
    \frac{3(k+1)}{2},
\]
then apply Theorem\tume\ref{sbs3} to the current puzzle on
$T_{\mathrm{black}}$ so that the remaining black vertices agree
with the configuration $g^{\prime}$.
Let $\sigma_0$ be the resulting sequence of moves.
Update
\[
    \sigma:=\sigma_0\cdot\sigma
\]
and return $\sigma$.
(\#\# Note that $\sigma(f)=g^{\prime}$.)
\label{initial-check}

\item
Set $t:=c_1(r)$.  
While $t$ is not a leaf node of $\BT$, update $t:=c_1(t)$.

\item
Call $\Alg.\Extraction(T,\BT,g',t)$,
obtaining $\sigma_t$.

\item
Update $\sigma := \sigma_t\cdot\sigma$.

\item
Set $\col(v):=\mathrm{white}$ for each 
$v\in V(X(t)-\JoiV(X(t)))$.  
\label{t-eliminate}

\item
Set $x:=p(t)$ and $y:=c_2(x)$.  
(\#\# Note that $t=c_1(x)$.)
\item
If $T_{\mathrm{black}}$ contains an isthmus of size at least
$k+1$, then choose an edge $uv$ with
\[
    u\in V(X(t)-\JoiV(X(t)))
    \qquad\text{and}\qquad
    v\in V(X(y)).
\]
Change the color of $u$ (which was set to white in
Step\tume\maru{\ref{t-eliminate}}) back to black,
find a leaf node $w$ of $\BT$ such that $v\in V(X(w))$, and update
\[
    X(w):=X(w)\cup uv.
\]
\label{hige}

\item
Delete the leaf node $t$ from $\BT$ and suppress its parent $x$ by
replacing $x$ with its remaining child $y$.
Thus $y$ inherits the position of $x$ in $\BT$ and becomes the root
if $x$ was the root.
Relabel this surviving node as $x$.
\label{del-t}

\item
If $x$ has become the root of $\BT$, re-root the surviving binary tree at $x$ as described in the paragraph \textit{Re-rooting} in Section~\ref{splitting-tree}.

\item
Go to Step~\ref{initial-check}.

\end{alg-enumerate}
\end{alg}

\noindent
The following procedure $\Alg.\Extraction(T, \BT,g^{\prime},t)$ is 
the only auxiliary routine called by \textsc{CutOff}.
Its role is to assemble the workspace required for the
current target subtree while preserving the workspace
invariant.
In the algorithm $\Alg.\CutOff(G,T,\BT,f,g^{\prime})$, it is necessary to gather the set of pebbles $g^{\prime}(V(X(t)))\setminus\{0\}$, which are scattered throughout the entire $T$, onto the target subtree $X(t)$. For this purpose, the following $\Alg.\Extraction(T,\BT,g^{\prime},t)$ is provided.
%アルゴリズム$\Alg.\CutOff(G,T,\BT,f,h)$においては、$T$全体に散在している$(g^{\prime})^{-1}(V(X(t)))$という碁石の集合を、
%target subtreeである $X(t)$ 上に集めてくる必要がある。そのためのアルゴリズムとして、$\Alg.\Extraction(\BT,h,t)$が用意されている。
\begin{alg}[$\Alg.\Extraction(T, \BT,g^{\prime},t)$]\label{alg:extraction}\hfill
\begin{namelist}{OUTPUT}
\item[\sl Input]{A tree $T$ of order $n + k + 1$ with maximum isthmus size at most $k$,
a binary tree $\BT$ storing subtrees of $T$, an output configuration $g^{\prime}$ on $(T,P)$ of $\Alg.\PacCon(G,T,\BT,g)$,
and a leaf node $t$ of $\BT$.}

\item[\sl Output]{A sequence of moves $\sigma_t$.}
\end{namelist}

\begin{alg-enumerate}

\item Copy the binary tree $\BT$ and denote it by $\BT2$.

\item
If
\[
    |V(T_{\mathrm{black}})|
    \leq
    \frac{3(k+1)}{2},
\]
then apply Theorem\tume\ref{sbs3} to the current puzzle on
$T_{\mathrm{black}}$ so as to place the pebbles
\[
    g^{\prime}(V(X(t)))\setminus\{0\}
\]
on $X(t)$ according to $g^{\prime}$.
After that, change the color of every gray vertex of $T$ back to
black, and return the sequence of moves $\sigma_t$ obtained so far.
\label{initial-check2}

\item
Let $r$ be the current root of $\BT2$.
Among the two children of $r$, choose the child whose descendant subtree
does not contain the target leaf $t$.
Starting from this child, repeatedly follow the first child until a
leaf node $\ell$ is reached.
In particular, $\ell\neq t$.

\item By Lemma\tume\ref{feasible-subtree}, find a subtree $H(\ell)$ of $T$
such that

\[
|V(H(\ell))|
=
|V(X(t))| + |V(X(\ell))| + k + 1 ,
\]

$X(\ell)$ is a subgraph of $H(\ell)$,
the maximum isthmus size of $H(\ell)$ is at most $k$,
and all vertices of $H(\ell)$ are black.
Move the $k+1$ empty spaces into $V(H(\ell))$
while recording the sequence of moves.
\label{mini-puzzle}

\item Apply Theorem~\ref{sbs3} to $H(\ell)$ and fill
$X(\ell)-\JoiV(X(\ell))$
with arbitrary pebbles that are not contained in
$g^{\prime}(V(X(t)))$,
leaving no empty space there.
Record the sequence of moves.
\label{mini-puzzle-finish}

\item Set $x := p(\ell)$ and let $y$ be the child of $x$ different from $\ell$.

\item Change the color of every vertex in
$V(X(\ell)-\JoiV(X(\ell)))$
to gray.

\item
If $T_{\mathrm{black}}$ contains an isthmus of size at least
$k+1$, then choose an edge $uv$ with
\[
    u\in V(X(\ell)-\JoiV(X(\ell)))
    \qquad\text{and}\qquad
    v\in V(X(y)).
\]
Change the color of $u$ from gray back to black,
find a leaf node $w$ of $\BT2$ such that $v\in V(X(w))$, and update
\[
    X(w):=X(w)\cup uv.
\]
\label{hige2}

\item
Delete the leaf node $\ell$ from $\BT2$ and suppress its parent $x$
by replacing $x$ with its remaining child $y$.
Thus $y$ inherits the position of $x$ in $\BT2$ and becomes the root
if $x$ was the root.
Relabel this surviving node as $x$.
\label{del-ell}

\item
If $x$ has become the root of $\BT2$, re-root the surviving binary tree at $x$ as described in the paragraph \textit{Re-rooting} in Section~\ref{splitting-tree}.

\item Go back to Step~\ref{initial-check2}.

\end{alg-enumerate}
\end{alg}

The following Lemma\tume\ref{subtrees-containing-v} is implicitly invoked in the complexity estimation 
within the proof of Theorem\tume\ref{main}, insofar as its computational contribution is asymptotically 
negligible compared to the other components.
\begin{la}\label{subtrees-containing-v}
Step\tume\maru{\ref{hige}} and Step\tume\maru{\ref{del-t}} of $\Alg.\CutOff(G,T,\BT,f,g^{\prime})$, as well as Step\tume\maru{\ref{hige2}} and Step\tume\maru{\ref{del-ell}} of 
$\Alg.\Extraction(T, \BT,g^{\prime},t)$, can each be achieved in $\Ord(k)$. 
\end{la}

The re-rooting operations do not alter the asymptotic running time. Let $m$ denote the number of nodes of the current binary decomposition tree; by Lemma\tume\ref{bunkatsu}, $m=\Ord((n+k)/k)$. Whenever the root is replaced, the new root is a child in the centroid decomposition, so the number of nodes below successive roots decreases geometrically. Hence all re-rooting operations during the whole execution of $\Alg.\CutOff$ take $\Ord(m)$ time, and those within one call of $\Alg.\Extraction$ also take $\Ord(m)$ time. Since $\Alg.\Extraction$ is called $\Ord(m)$ times, the total re-rooting time is $\Ord(m^2)=\Ord(n^2+1)$, which is absorbed by the bound in Theorem\tume\ref{main}.

\begin{la}\label{mini-jyunbi}
The total cost of Step\tume\maru{\ref{mini-puzzle}} of $\Alg.\Extraction(T,\BT,g^{\prime},t)$ is in $\Ord(k(k+n))$.  
%上記のサブルーチン $\Alg.\Extraction(T,\BT,h,t)$のステップ\maru{\ref{mini-puzzle}}にかかるトータルの手間は、$\Ord(k(k+n))$である。
\end{la}

%%%%%%%%%%%%%%%%%%%%%%%%%%%%%%%%%%%%%%%%%%%%%%%%
%%%%%%%%%%%%%%%%%%%%%%%%%%%%%%%%%%%%%%%%%%%%%%%%
%%%%%%%%%%%%%%%%%%%%%%%%%%%%%%%%%%%%%%%%%%%%%%%%
%%%%%%%%%%%%%%%%%%%%%%%%%%%%%%%%%%%%%%%%%%%%%%%%
%%%%%%%%%%%%%%%%%%%%%%%%%%%%%%%%%%%%%%%%%%%%%%%%

\subsection{Main Algorithm}

We now present the main algorithm, which combines the procedures described above.
%以上に述べた一連のアルゴリズム達の機能を集積することにより、下記のmain algorithmが構成される。
\begin{alg}[{\bf Main Algorithm:} $\Alg.\Main((G,P), f, g)$]\hfill
\begin{namelist}{OUTPUT}
\item[\sl Input]{A tree $G$ with $N$ vertices, with $k:=k(G)$ denoting
its maximum isthmus size; a set of pebbles $P:=\{1,\ldots,n\}$ such
that $n+k<N$; and an initial configuration $f$ and a goal
configuration $g$ on the pair $(G,P)$.}
\item[\sl Output]{An explicit sequence of moves to rearrange from the initial configuration $f$ to the goal configuration $g$.}
\end{namelist}
\begin{alg-enumerate}
\item{By using Lemma\tume\ref{feasible-subtree}, find a subtree $T$ of $G$ with $n+k+1$ vertices such that the maximum size of isthmuses in $T$ is at most $k$.}
%\item{Lemma\tume\ref{feasible-subtree}を用いて、位数が$n+k+1$で最大isthmus長が$k$以下となるような$G$の部分木$T$を見つける}
\item{By using the algorithm from Lemma\tume\ref{bunkatsu}, decompose the tree $T$ into subtrees of order $k/8$ or less, and store their data in the binary tree $\BT$.}
%\item{（最適化された）ボードグラフである位数$k+n+1$の木$T$を、centroidで位数$k/8$以下のsubtreeに分解し、それを二分木データベース$\BT$に納める}
\item{By using $\Alg.\PacCon(G,T,\BT, g)$, obtain the configuration $g^{\prime}:=\phi[r]\cdot\psi(g)$ on $T$.}
%\item{$\Alg.\PacCon(G,T,\BT, g)$を用いて、$T$上の配置$g^{\prime}:=\phi[r]\cdot\psi(g)$を得る}
\item{Call $\Alg.\CutOff(G,T,\BT,f,g^{\prime})$.}
%\item{$\Alg.\CutOff(G,T,\BT,f,g^{\prime})$を呼ぶ}
\item{Return the sequence of moves $\psi^{-1}\cdot\phi[r]^{-1}\cdot\sigma$. (\#\# Note that $g=\psi^{-1}\cdot\phi[r]^{-1}\cdot\sigma(f)$.)}
%\item{$\psi^{-1}\cdot\phi[r]^{-1}\cdot\sigma$というmoveの系列を出力して停止する（$g=\psi^{-1}\cdot\phi[r]^{-1}\cdot\sigma(f)$に注意）}
\end{alg-enumerate}
\end{alg}

\section{Board Trees: Proof of Theorem\tume\ref{main}}\label{Main1Proof}

In this section, we present a proof of Theorem\tume\ref{main}.

\noindent
\begin{Proofof}{Theorem\tume\ref{main}}
Since the correctness of $\Alg.\Main((G, P), f, g)$ is straightforward, 
we focus here on analyzing its time complexity, which is shown to be 
$\Ord(N + n\D(G) + n^2\log(1+\min\{n,k\}))$, 
and on proving that the length of the output solution sequence is bounded by 
$\Ord(n\D(G) + n^2\log(1+\min\{n,k\}))$. 

It is worth noting that the initial $\Ord(N)$ term in the time complexity accounts 
for the cost of scanning the input tree $G$ and the configurations $f$ and $g$, 
in particular the positions of the pebbles. 
This cost does not contribute to the length of the solution sequence, and thus 
does not appear in the asymptotic bound for the solution length.

From a given input tree $G$, the effort of finding its subtree $T$ of order $n+k+1$ using Lemma\tume\ref{feasible-subtree} is $\Ord(n+k)$.
%与えられた入力tree $G$から、Lemma\tume\ref{feasible-subtree}を使ってその位数$n+k+1$の部分木$T$を見つける手間は$\Ord(n+k)$である。
The task of creating a binary tree database $\BT$ from the obtained $T$ using Lemma\tume\ref{BT} has a cost of $\Ord((n+k)\log(n+k))$.
%得られた$T$から、Lemma\tume\ref{BT}を用いて二分木データベース$\BT$を作る作業には、$\Ord((n+k)\log(n+k))$の手間が生じる。
Using Lemma\tume\ref{sbs0} to convert the arrangements $f$ and $g$ on $(G,P)$ into strongly related arrangements on $(T,P)$ requires a time complexity and number of moves of at most $\Ord(n\D(G))$.
%Lemma\tume\ref{sbs0}を用いて$(G,P)$上の配置$f$や$g$をstrongly relatedな$(T,P)$上の配置に変換する作業には、a time complexity and number of moves of at most $\Ord(Nn)$を要する。

In $\Alg.\Main((G,P), f, g)$, the other parts are spent on repeatedly
solving sub-puzzles of size $\Ord(k)$
(Steps \maru{\ref{mini-puzzle}} and
\maru{\ref{mini-puzzle-finish}} of
$\Alg.\Extraction(T,\BT,g^{\prime},t)$) on $(T,P)$.
%Mainアルゴリズムにおけるそれ以外のpartは、$(T,P)$上で行われるサイズ$\Ord(k)$のサブパズル（$\Alg.\Extraction(T, \BT,h,t)$のSteps \maru{\ref{mini-puzzle}} and \maru{\ref{mini-puzzle-finish}}）を繰り返し解く操作により費やされる。
According to Lemma\tume\ref{mini-jyunbi}, the repeated process of
designating a suitable area of $\Ord(k)$ and gathering $(k+1)$
unoccupied spaces in that location
(Step\tume\maru{\ref{mini-puzzle}} of
$\Alg.\Extraction(T,\BT,g^{\prime},t)$)
for preparing those sub-puzzles has a time complexity and number of
moves of at most $\Ord(k(k+n))$ for the entire
$\Alg.\Extraction(T,\BT,g^{\prime},t)$.
%Lemma\tume\ref{mini-jyunbi}によれば、それらのサブパズルを準備するために、適切な$\Ord(k)$の領域を定めてその場所に$k+1$個のunoccupied spacesを寄せてくる作業（$\Alg.\Extraction(T, \BT,h,t)$のStep\tume\maru{\ref{mini-puzzle}}）の繰り返しは、$\Alg.\Extraction(T, \BT,h,t)$全体として、a time complexity and number of moves of at most $\Ord(k(k+n))$を要することが分かっている。
Therefore, in $\Alg.\Main((G,P), f, g)$, the most effort and number of moves required is from the repetition of sub-puzzles of size $\Ord(k)$ (Step\tume\maru{\ref{mini-puzzle-finish}} of $\Alg.\Extraction(T, \BT,g^{\prime},t)$).
%したがって、メインアルゴリズムにおいて、もっとも多くの手間とmoveの数を必要とするのは、サイズ$\Ord(k)$のサブパズル（$\Alg.\Extraction(T, \BT,h,t)$のStep\tume\maru{\ref{mini-puzzle-finish}}）の繰り返しである。
This task is repeated $\Ord((n+k)/k)$ times each time $\Alg.\Extraction(T, \BT,g^{\prime},t)$ is called. If each of these sub-puzzles of size $\Ord(k)$ meets the prerequisites of Theorem\tume\ref{sbs3} (the proof for this part will be postponed until the end), then the total effort and number of moves for one call to
$\Alg.\Extraction(T,\BT,g^{\prime},t)$ is at most $\Ord(((n+k)/k) \cdot (k^2\log k))=\Ord((n+k)k\log k)$. 
%この作業は、$\Alg.\Extraction(T, \BT,h,t)$を呼ぶ度に$\Ord((n+k)/k)$回繰り返されるが、これら一つ一つのサイズ$\Ord(k)$のサブパズルがTheorem\tume\ref{sbs3}の前提条件を満たすのであれば（この部分の証明は最後に回す）、その手間とmoveの総数は、一回の$\Alg.\Extraction(T, \BT,h,t)$について、高々$\Ord(((n+k)/k) \cdot (k^2\log k))=\Ord((n+k)k\log k)$である。
$\Alg.\Extraction(T, \BT,g^{\prime},t)$ is called $\Ord((n+k)/k)$ times within $\Alg.\CutOff(G,T,\BT,f,g^{\prime})$, so the total cost and estimate of moves here is $\Ord(((n+k)/k)(k(k+n)+(n+k)k\log k ))=\Ord((n+k)^2\log k )$.
%$\Alg.\Extraction(T, \BT,h,t)$は、$\Alg.\CutOff(G,T,\BT,f,h)$の中で$\Ord((n+k)/k)$回呼ばれるから、ここでの総コストとmoveの手数の見積は、
%$\Ord(((n+k)/k)(k(k+n)+(n+k)k\log k ))=\Ord((n+k)^2\log k )$となる。

Since the operation of scaling down the input puzzles $f,g$ on $(G,P)$ to their strongly related puzzles on $(T,P)$ requires $\Ord(n\D(G))$ effort and moves, assuming $n > k$, the estimate for the total effort and number of moves in the main algorithm is $\Ord(n\D(G) +(n+k)^2 \log k)=\Ord(n\D(G)+n^2 \log k)$.
%$(G,P)$上のパズルを$(T,P)$上のパズルにスケールダウンする操作に$\Ord(Nn)$の手間がかかるので、$n>k$を仮定すると、メインアルゴリズムの手間とmoveの総数の見積は、$\Ord(Nn +(n+k)^2 \log k)=\Ord(Nn+n^2 \log k)$となる。
On the other hand, if $n \leq k$, then Step\tume\maru{\ref{mini-puzzle-finish}} of $\Alg.\Extraction(T, \BT,g^{\prime},t)$ is repeated only a constant number of times for the entire main algorithm, so the total effort estimate is $\Ord(N+n\D(G) + n^2\log n)$, just as stated in Theorem\tume\ref{sbs3}.
%一方、$n\leq k$であれば、（$\Alg.\Extraction(T, \BT,h,t)$のStep\tume\maru{\ref{mini-puzzle-finish}}）の繰り返しは、メインアルゴリズム全体として定数回しか繰り返されないため、トータルの手間の見積は、Theorem\tume\ref{sbs3}そのままに、$\Ord(Nn + n^2\log n )$となる。

Lastly, we will confirm that repeatedly performing Step\tume\maru{\ref{hige}} of $\Alg.\CutOff(G,T,\BT,f,g^{\prime})$ or Step\tume\maru{\ref{hige2}} of $\Alg.\Extraction(T, \BT,g^{\prime},t)$ does not increase the order of $X(w)$ beyond twice its original value. In this process, the edge $uv$ added to $X(w)$ is a pendant edge, and since $u$ is a leaf, by definition it is not a vertex on any isthmus. Therefore, even if this operation is repeated, no pendant edge will be added to $u$. By adding the edge $uv$, the size of a maximum isthmus has strictly decreased, which means that the degree of $v$ in the current graph is $3$ or more. Consequently, the degree of $v$ in $X(w)$ is either $3$ or a unique vertex $p$ outside $X(w)$ is adjacent to $v$. By definition, a vertex with a degree of $3$ or more is not an interior vertex of any isthmus, so as long as the degree of $v$ in the current graph is $3$, no further pendant edges will be added to $v$ through the repetition of this operation. For the degree of $v$ to again become $2$, the only option is to delete the vertex $p$ outside $X(w)$ and the edge $pv$, but since $u$ is not a vertex on any isthmus, the size of a maximum isthmus in the remaining graph will not increase by this deletion of $pv$. Therefore, the number of pendant edges that can potentially be added to each vertex of $X(w)$ is at most $1$, and since no pendant edges will be added to the newly added leaves, the order of $X(w)$ can increase by at most a factor of $2$. This implies that $k + 1 \geq k=2(2(k/8)+2(k/8)) \geq 2(|V(X(t))| + |V(X(\ell))|)$ holds, and therefore the condition of Theorem\tume\ref{sbs3} (i.e. $|V(X(t))| + |V(X(\ell))| + k + 1 \geq 3(|V(X(t))| + |V(X(\ell))|)$) is always satisfied in Step\tume\maru{\ref{mini-puzzle-finish}} of $\Alg.\Extraction(T, \BT,g^{\prime},t)$. 
\end{Proofof}

%%%%%%%%%%%%%%%%%%%%%%%%%%%%%%%%%%%%%%%%%%%%%%%%%%%%%%%%%%%%%%%%%%%%%%%%%%%%%%%%
%%%%%%%%%%%%%%%%%%%%%%%%%%%%%%%%%%%%%%%%%%%%%%%%%%%%%%%%%%%%%%%%%%%%%%%%%%%%%%%%
%%%%%%%%%%%%%%%%%%%%%%%%%%%%%%%%%%%%%%%%%%%%%%%%%%%%%%%%%%%%%%%%%%%%%%%%%%%%%%%%
\section{General Board Graphs: Proof of Theorem\tume\ref{main2}}\label{Main2Proof}

In this section, we provide a proof of Theorem\tume\ref{main2}. 
Now that the efficient algorithm for the \textsc{Pebble Motion Problem on Trees} and 
its computational complexity evaluation has been clarified in the previous sections, 
the actual proof of this theorem can be obtained in a simple manner. 

\begin{Proofof}{Theorem\tume\ref{main2}}
First, computing the $2$-connected component decomposition of the input graph $G$, as well as 
constructing a breadth-first search spanning tree $T_G$ of $G$, requires 
$\Ord(|E(G)|)$ time. 
However, this computational cost is purely algorithmic overhead and does not affect the length of 
the output solution sequence in any way.

In Step\tume\maru{1} of the Main Algorithm ($\Alg.\Main((G, P), f, g)$), when finding a subtree $T$ 
of order $n + k + 1$ of $T_G$, we apply Lemma\tume\ref{feasible-subtree} separately to each part of $T_G$ 
corresponding to the $2$-connected component decomposition of $G$. 
That is, for each $2$-connected component $X$ of $G$, it suffices to ensure that if a part of the subtree $Y$ 
of $T_G$ that spans $X$ is included in $T$, then the size of that part does not exceed the maximum isthmus 
length of $Y$.
It is the same in that $T$ is decomposed into $\Ord(k)$ size subtrees and represented as binary tree data. 
The difference appears in the computational complexity evaluation. 
%我々がここで議論する2-連結ボードグラフGの上のアルゴリズムは、input board treeの代わりにGの適当な
%全域木Tを用い、パラメータkをmaxmum isthmus sizeではなくk:=N-nと定義しなおすこと以外は、木の上の
%アルゴリズムと変わるところはない。Tを細かいO(k)のsubtreeに分解し、それらを二分木データとして表現する
%ことも同じである。違いは計算量評価の方に現れる。
If the input board graph is $2$-connected then the puzzle is feasible if and only if the board graph is 
not a cycle and has at least two unoccupied spaces. Therefore, no matter how the tree $T$ is 
constructed, $T$ may inevitably have to possess an isthmus with at least $k+1$ vertices, 
and as a result, some of the subtrees in the decomposition of $T$ may have long isthmuses of size
$k+1$ or more. It should be noted that each subtree can contain at most one isthmus of size 
at least $k+1$. 
In such cases, to create a feasible sub-puzzle of size $\Ord(k)$
(c.f. Step\tume\ref{mini-puzzle} of
$\Alg.\Extraction(T,\BT,g^{\prime},t)$), the puzzle is constructed
around a vertex $v$ of degree at least $3$ lying on a cycle of $G$
that bypasses the obstructing isthmus.  We temporarily rotate the
pebbles along such a cycle only as far as is necessary to perform the
required local exchange near the isthmus, and after the exchange we
reverse the same rotation.  Hence every vertex outside the active
local part is restored to exactly the same pebble or blank as before
the repair.  In particular, a cycle used for this purpose may pass
through a part that has already been completed or is temporarily
inactive: the configuration on that part is unchanged after the
repair, so the recursive state maintained by the tree algorithm is
preserved.

The $\Ord(k)$ pebbles forming the sub-puzzle can therefore be brought
around $v$ with
\[
    \Ord\bigl(k\min\{n,\CL(G)\}\bigr)
\]
moves for each such sub-puzzle.  Since the main algorithm solves
$\Ord((n/k)^2)$ such sub-puzzles, the total additional cost is
\[
    \Ord\!\left(
        \left(\frac{n}{k}\right)^2
        k\min\{n,\CL(G)\}
    \right)
    =
    \Ord\!\left(
        \frac{n^2\min\{n,\CL(G)\}}{k}
    \right).
\]
%In such cases, to create a feasible sub-puzzle of size $\Ord(k)$
%(c.f. Step\tume\ref{mini-puzzle} of $\Alg.\Extraction(T, \BT,g^{\prime},t)$), 
%the puzzle must be constructed around a vertex $v$ with degree $3$ or higher located somewhere 
%on a cycle of $G$ that includes the isthmus. To achieve this, the $\Ord(k)$ pebbles that make up 
%the sub-puzzle must be brought around $v$. This operation requires $\Ord((k)\times \min\{n,\CL(G)\})$
%moves per each sub-puzzle.
%In our algorithm, since we have to solve these sub-puzzles $\Ord((\frac{n}{k})^2)$ times, 
%the cost of the above operations for the entire algorithm can be estimated as 
%$\Ord((\frac{n}{k})^2 \times k \times \min\{n,\CL(G)\})=\Ord(\frac{n^2\min\{n,\CL(G)\}}{k})$.
%cycle以外の$2$-connected board graphにおいては、$2$個以上の空白があればfeasibleであるから、
%$G$のspanning subtree $T$をどのように工夫して取ろうとしても、$T$が長いisthmusを含み、その結果
%$T$の分解におけるsubtreesの一部がそのisthmus上のO(N-n)のsub path達になってしまう可能性は
%避けられない。そうした場合、O(N-n)サイズのfeasibleなsub-puzzleを作るには、そのisthmusを含む
%Gの部分サイクル上のどこかにある次数３以上の頂点vを中心にしてpuzzleを構成するしかない。
%そのためには、当該sub-puzzleを構成するO(N-n)個のpebble達を、vの周囲に持ってくる必要がある。
%この操作には、一つのsub-puzzle当たり、O((N-n)n)のmoveの手間が必要となる。
%我々のアルゴリズムにおいては、これらのsub-puzzle達をO((n/(N-n))^2)回解かなければならないから、
%上記の手数のコストはアルゴリズム全体ではO((n/(N-n))^2 * (N-n)n)=O(n^3/(N-n))と見積もられる。

Note that, for a given instance $(G, P)$, the feasibility of the puzzle is guaranteed as long as 
$k$ satisfies $k \leq N - n - 1$. 
Therefore, in the above discussion, the constant $k$ in the proof may be replaced with an arbitrary 
fixed integer $\theta$ satisfying $k \leq \theta \leq N - n - 1$, without affecting the correctness of the proof.

The statement of our theorem can be directly derived from the above fact and Theorem\tume\ref{main}. 
\end{Proofof}

%%%%%%%%%%%%%%%%%%%%%%%%%%%%%%%%%%%%%%%%%%%%%%%%%%%%%%%%%%%%%%%%%%%%%%%%%%%%%%%%
%%%%%%%%%%%%%%%%%%%%%%%%%%%%%%%%%%%%%%%%%%%%%%%%%%%%%%%%%%%%%%%%%%%%%%%%%%%%%%%%
%%%%%%%%%%%%%%%%%%%%%%%%%%%%%%%%%%%%%%%%%%%%%%%%%%%%%%%%%%%%%%%%%%%%%%%%%%%%%%%%
\section{Diameters, Shortest Cycles, and Refinements}\label{refinement}

In this section, we present a refined version of the algorithm introduced in the preceding sections.
This refinement enables us to demonstrate that, for any instance, our method achieves provable asymptotic improvements—not only in the order
of the length of the resulting move sequence, but also in the overall computational complexity—when compared with existing approaches that
place pebbles sequentially, one at a time, onto predetermined vertices.

\subsection{Refined Analysis for Board Trees: Proof of Theorem\tume\ref{refinement1}}\label{refinement-tree}

We first note that the time complexity stated in Theorem\tume\ref{main} can be improved to 
$\Ord(N+n\D(G) + n^2 \log(1+\min\{n,k\}))$, and that the length of the solution sequence produced by 
the algorithm admits a refined upper bound of $\Ord(n\D(G) + n^2 \log(1+\min\{n,k\}))$.
In fact, this has already been explicitly demonstrated in the proof of Theorem\tume\ref{main}.

% 我々が本論文のprevious sectionまでに提案したアルゴリズムでは、その中のAlgorithm 4 $(\Alg.\Extraction(T,BT,g',t))$の\maru{4}において、
% 現在のboard treeから除去してもreminder graphが連結のままであるようなsub tree $X(\ell)$を指定し、
% この$X(\ell)$を含み頂点数が$|V(X(t))|+|V(X(\ell))|+k+1$である$H(\ell)$というsubtree (ここで$X(t)$は配置を確定させるターゲットのsubtree)を探して、
% $H(\ell)$というboard tree上のsub puzzleを解き、$X(\ell)$上にゴミ（$X(t)$に乗せるべきpebble以外の碁石達）を詰めて切り落とす手続き（\maru{5}）を繰り返します。
In the algorithm proposed up to the previous section of this paper, Step\tume\maru{\ref{mini-puzzle}} of Algorithm\tume\ref{alg:extraction} ($\Alg.\Extraction(T, \BT,g^{\prime},t)$) performs the following procedure iteratively:
we identify a subtree $X(\ell)$ of the current board tree such that removing $X(\ell)$ does not disconnect the remainder graph. We then search for a subtree $H(\ell)$ that
contains $X(\ell)$ and has exactly $|V(X(t))| + |V(X(\ell))| + k + 1$ vertices, where $X(t)$ denotes the target subtree on which the placement of pebbles is to be correctly finalized.
Next, we solve the sub-puzzle on $H(\ell)$, pack irrelevant pebbles (i.e., pebbles not intended to be placed on $X(t)$) into $X(\ell)$, and cut off $X(\ell)$ from the board tree (\maru{5}).

% このときまず、
% １） $X(\ell)$の付け根$\JoiV(X(\ell))$と、$X(t)$の付け根$\JoiV(X(t))$とを結ぶ（木ですから『唯一の』）パス $P(\ell,t)$を考えると、$X(\ell)$のどの頂点から$X(t)$のどの頂点に至る最短路も$P(\ell,t)$を必ず使います。
% なので、$P(\ell,t)$のJointVeretx$(X(\ell))$をその一方の端点とする、$P(\ell,t)$の十分長い（少なくとも$2|V(X(t)|+1$以上の頂点数を持つ）サブパス$P'(H(\ell))$が$H(\ell)$に含まれるように、$H(\ell)$を取ることにします。
% これは、$H(\ell)$の頂点数が$|V(X(t))|+|V(X(\ell))|+k+1$ということで、十分余裕がありますので、$X(\ell)$が$X(t)$に極めて近い特別なケース以外では、必ずそのようなサブパス$P'(H(\ell))$を含むように$H(\ell)$を取ることが出来ます。
% 逆に、$P(\ell,t)$の頂点数が$2|V(X(t)|$以下であれば、$H(\ell)$として、$X(t)$を含むように取ってくることが出来ますので、その場合にはそうします。
% ２）上記１）の条件を満たすように取った$H(\ell)$のパス$P'(H(\ell))$の上に、もし、$X(t)$に置かれるべきpebbleが何個か既に乗っていた場合には、それらのpebble達達が乗っている頂点は触らずに、今$X(\ell)$上にある「$X(t)$に置かれるべきpebble達」を、$P'(H(\ell))$の「それ以外の」「できる限り$X(\ell)$から離れた位置にある」頂点のどこかに乗っている（ゴミ）pebble達と入れ替え（位置の交換）をします。
% この処理が終了した段階で、必然的に、$X(\ell)$には「ゴミのみ」が詰められることになります。
% なお、$P(\ell,t)$の頂点数が$2|V(X(t)|$以下であれば、$H(\ell)$は$X(t)$を含みますから、$X(t)$上の正しい位置の頂点に$H(\ell)$上のすべてのpebble達を正しく配置してしまえばいいです。
In the refinement of the algorithm proposed in this section, we proceed as follows:
\begin{enumerate}
      \item Consider the unique path $P(\ell, t)$ in the board tree connecting the joint vertex of $X(\ell)$, denoted by $\JoiV(X(\ell))$, and the joint vertex of $X(t)$, denoted by $\JoiV(X(t))$.
Since the underlying structure is a tree, this path is uniquely determined. Any shortest path from a vertex in $X(\ell)$ to a vertex in $X(t)$ must pass through $P(\ell, t)$.
Based on this observation, we construct $H(\ell)$ so that it contains a sufficiently long subpath $P'(H(\ell)) \subseteq P(\ell, t)$, where $P'(H(\ell))$ starts at $\JoiV(X(\ell))$
and contains at least $2|V(X(t))| + 1$ vertices. Since $H(\ell)$ is allowed to have $|V(X(t))| + |V(X(\ell))| + k + 1$ vertices, it generally has enough capacity to include such a subpath
unless $X(\ell)$ happens to be located very close to $X(t)$. In most cases, therefore, we can ensure that $H(\ell)$ includes $P'(H(\ell))$ as required. Conversely, if the total number
of vertices in $P(\ell, t)$ is at most $2|V(X(t))|$, we instead construct $H(\ell)$ so that it includes $X(t)$ itself.
      \item Suppose that some of the pebbles that are supposed to be placed on $X(t)$ are already located on vertices in $P'(H(\ell))$. In that case, we avoid moving those pebbles.
Instead, for the pebbles currently located in $X(\ell)$ that are eventually intended to be placed on $X(t)$, we swap their positions with "garbage" pebbles (i.e., pebbles not
destined for $X(t)$) located on other vertices of $P'(H(\ell))$, choosing such vertices to be as far from $X(\ell)$ as possible. After this process, all pebbles remaining in $X(\ell)$
are guaranteed to be garbage. Note that if the number of vertices in $P(\ell, t)$ is at most $2|V(X(t))|$, then $H(\ell)$ includes $X(t)$ by construction. In this case, we can directly
solve the sub-puzzle on $H(\ell)$ by placing all pebbles correctly on their target vertices within $X(t)$.
\end{enumerate}

% ２）の手続きは、「本来$X(t)$に置かれるべき$X(\ell)$上のpebble達」の個数$k_{\ell}$が、$\log k$未満である場合には、$\Ord(k^2)$の手間を$k_{\ell}$回繰り返して、トータルで$\Ord(k_{\ell} k^2)$の計算量で達成可能です。
% （ここに、我々の論文のSection 3で述べている「十分な量（pebbleの数の２倍以上）の空白の個数がある場合」の$N$-vertex board treeのケースで、任意の２つのpebbleの位置『のみを』交換する手間が、$\Ord(Nn + n^2)$であること、を使います。
% 実際、交換すべき２つのpebbleを、髭一本の場合であれば髭の両脇あたりに、髭二本のケースでは２つの髭の上にそれぞれを一つずつ、配置した後、髭を使ってそれら２つのpebbleの位置の交換をし、その後、他のpebble達について行った前処理の逆操作で状態を戻せばいいです。
% Section 3の設定では、スペースが十分確保されているおかげで、目的の２つのpebbleの交換作業中は、他のすべてのpebble達は適当な「仮置き場」に避難させておく（前処理）ことができるので、この避難を元に戻す作業を含めて、この作業は全体として$\Ord(Nn+n^2)$で達成可能です。（←この事実の証明については、Section 3のどこかに、中上川先生の方で記載して頂くことは可能でしょうか？）
% もちろん、$k_{\ell}$が、$\log k$以上である場合には、我々がすでにSection 3のTheorem8で示しているアルゴリズムを使って$\Ord(k^2 \log k)$の手間で処理が完了します。

% この１）、２）の付随条件を課しますと、$X(t)$上に置かれるべきpebble達は、初期状態で$G$のどの頂点上に置かれていたとしても、$\Ord(k)$のサイズのboard tree上のサブパズルを1回解くごとに、スタート時の頂点と自身が置かれるべき$X(t)$上のゴールの頂点とを結ぶ$G$上の最短路に沿って、$\Theta(k)$ずつ着実に歩を進める、ことになります。
% なので、この方法であれば、moveの手数は（$G$の直径を$D(G)$とすると）、$\Ord((N-n)n+\min\{k n D(G), n^2 \log(1+\min\{ n, k \})\})$ で達成可能であることが分かります。
The procedure described in step 2 of the above refinement can be executed with a total computational cost of $\Ord(k_{\ell} k^2)$
when the number of pebbles on $X(\ell)$ that are supposed to be placed on $X(t)$, denoted $k_{\ell}$, is less than $\log k$.
In this case, we perform a swap operation of two pebbles at most $k_{\ell}$ times, each of which costs $\Ord(k^2)$.
To achieve this, we make use of the result from Lemma\tume\ref{sbs4} in Section\tume\ref{base-cases} of this paper:
in a board tree with $N$ vertices and a sufficient number of unoccupied vertices (at least twice the number of pebbles),
swapping the positions of only two pebbles can be done in $\Ord(Nn)$ time. 
Of course, when $k_{\ell} \geq \log k$, we instead apply the algorithm already presented in Theorem\tume\ref{sbs3} of 
Section\tume\ref{base-cases}, which achieves the operation in $\Ord(k^2 \log k)$ time.

By incorporating the refinements 1) and 2) described above, the pebbles destined for $X(t)$ will, regardless of
their initial positions in the board graph $G$, make consistent progress along the shortest path in $G$ from their initial positions
to their target vertices in $X(t)$, by a distance of $\Theta(k)$ each time we solve a sub-puzzle on $H(\ell)$ such that $X(\ell)$
contains such a pebble. 

%%%%%%%%%%%%%%%%%%%%%%%%%%%%%%%%%%%%%%%%%%%%%%%%%%%%%%%%%%%%%%%%%%%%%%%%%%%%%%%%
%%%%%%%%%%%%%%%%%%%%%%%%%%%%%%%%%%%%%%%%%%%%%%%%%%%%%%%%%%%%%%%%%%%%%%%%%%%%%%%%
%%%%%%%%%%%%%%%%%%%%%%%%%%%%%%%%%%%%%%%%%%%%%%%%%%%%%%%%%%%%%%%%%%%%%%%%%%%%%%%%
At the beginning of this section, we stated that the refinement leads to an algorithm whose output exhibits an asymptotically shorter solution
--- compared to existing methods that sequentially reposition pebbles, one by one, onto predetermined locations.
In what follows, we examine this fact in greater detail.
%%%%%%%%%%%%%%%%%%%%%%%%%%%%%%%%%%%%%%%%%%%%%%%%%%%%%%%%%%%%%%%%%%%%%%%%%%%%%%%% 
%%%%%%%%%%%%%%%%%%%%%%%%%%%%%%%%%%%%%%%%%%%%%%%%%%%%%%%%%%%%%%%%%%%%%%%%%%%%%%%%
%%%%%%%%%%%%%%%%%%%%%%%%%%%%%%%%%%%%%%%%%%%%%%%%%%%%%%%%%%%%%%%%%%%%%%%%%%%%%%%%
Let $P_{X(t)}$ denote the set of pebbles that are to be placed on the target subtree $X(t)$.
For each pebble $p$ in $P_{X(t)}$, let $\SP(p)$ denote the unique path in the board tree $G$ from the vertex at which $p$ initially resides
to the vertex in $X(t)$ where $p$ is to be placed.
We associate a formal variable $x_v$ with each vertex $v$ of $G$, and define the polynomial in ${\mathbf x}:=(x_v)_{v\in V(G)}$ 
corresponding to the path $\SP(p)$ as
$$
F_{\SP(p)}({\mathbf x}) := \sum_{v \in V(\SP(p))} x_v,
$$
and 
$$
F({\mathbf x}) := \sum_{p \in P_{X(t)}} F_{\SP(p)}({\mathbf x}).
$$
By expanding $F({\mathbf x})$, let $a_v$ denote the coefficient of the variable $x_v$ corresponding to each vertex $v$ of $G$; that is,
$$
F({\mathbf x}) = \sum_{v \in V(G)} a_v x_v.
$$
Note that each $a_v$ is a non-negative integer. 
Then, we consider the following quantity:
$$
\#\SP(X(t)) := \sum_{v \in V(G)} \min\{a_v, \log k\}.
$$
The total cost of Step\tume\maru{\ref{mini-puzzle-finish}} in the refinement of Algorithm\tume\ref{alg:extraction} ($\Alg.\Extraction(T, \BT, g', t)$) 
can be regarded as performing $\Ord(k^2)$ units of computational effort (and moves) up to $\Ord\big(\frac{\#\SP(X(t))}{k}\big)$ times.
Hence, we have:
\begin{prop}\label{refinement1-path-weight}
The above refinement of  Algorithm\tume\ref{alg:extraction} $(\Alg.\Extraction(T, \BT,g^{\prime},t))$  can be stopped with a time complexity in $\Ord(k \cdot \#\SP(X(t)))$.
\end{prop}
Note that in existing approaches in the literature, where pebbles are placed one by one on the board tree,
the above complexity can grow to the order of $\Omega(k\sum_{p\in P_{X(t)}} |\SP(p)|)=\Omega(k \sum_{v \in V(G)} a_v)$.

%%%%%%%%%%%%%%%%%%%%%%%%%%%%%%%%%%%%%%%%%%%%%%%%%%%%%%%%%%%%%%%%%%%%%%%%%%%%%%%% 
%%%%%%%%%%%%%%%%%%%%%%%%%%%%%%%%%%%%%%%%%%%%%%%%%%%%%%%%%%%%%%%%%%%%%%%%%%%%%%%%
%%%%%%%%%%%%%%%%%%%%%%%%%%%%%%%%%%%%%%%%%%%%%%%%%%%%%%%%%%%%%%%%%%%%%%%%%%%%%%%%
%%%%%%%%%%%%%%%%%%%%%%%%%%%%%%%%%%%%%%%%%%%%%%%%%%%%%%%%%%%%%%%%%%%%%%%%%%%%%%%%
%%%%%%%%%%%%%%%%%%%%%%%%%%%%%%%%%%%%%%%%%%%%%%%%%%%%%%%%%%%%%%%%%%%%%%%%%%%%%%%%
%%%%%%%%%%%%%%%%%%%%%%%%%%%%%%%%%%%%%%%%%%%%%%%%%%%%%%%%%%%%%%%%%%%%%%%%%%%%%%%%
For each pebble $p \in P_{X(t)}$, it is clear that $|\SP(p)| \leq \D(T)+1 \leq \D(G)+1$, where 
$\D(T)$ and $\D(G)$ denote the diameters of the board graphs $T$ and $G$, respectively.
Therefore, we have $\#\SP(X(t)) \leq \sum_{p\in P_{X(t)}} |\SP(p)| \in \Ord(k \D(T))$.
Within the main algorithm, Algorithm\tume\ref{alg:extraction} ($\Alg.\Extraction(T, \BT, g', t)$) is called at most $\Ord(n / k)$ times.
Therefore, the overall computational complexity contributed by this refinement is $\Ord(n k \D(T))$, which is 
at most $\Ord(n k \D(G))$.
Together with the result of Theorem\tume\ref{main}, this completes the proof of Theorem\tume\ref{refinement1}. \qed

We conclude this subsection with a brief remark on Theorem\tume\ref{refinement1}.
In our main algorithm $\Alg.\Main((G,P), f, g)$, we first reduce the pebble motion puzzle on 
the instance's board tree $G$ to a smaller puzzle on a subtree $T$ of $G$, where $T$ is 
a board graph with $k+n+1$ vertices. 
Then, we solve the puzzle on this reduced subtree $T$ and restore the resulting solution to 
the original board tree $G$ with $\Ord(N+n\D(G))$ effort and $\Ord(n\D(G))$ moves. 
Therefore, more precisely, the occurrence of $\D(G)$ inside the $\min\{\}$ expression 
in the second term of the asymptotic analysis in Theorem\tume\ref{refinement1} should 
instead be written as $\D(T)$.
However, since the instance graph is $G$ and the tree $T$ is merely an internal variable used 
within the algorithm, we have expressed the bound using $\D(G)$ rather than $\D(T)$.

As noted in the previous section, here as well, all occurrences of the parameter $k$ can be 
simultaneously replaced by any fixed integer $\theta$ satisfying $k \leq \theta \leq N - n - 1$, 
without invalidating the statement.

\subsection{Refined Analysis for General Board Graphs}\label{refinement-general}

As in Section~\ref{Main2Proof}, we first select an arbitrary vertex in the
connected board graph $G$ and construct a spanning tree $T_G$ by performing
a breadth-first search.
There is a natural relationship between the diameter $\D(G)$ of $G$
and the diameter $\D(T_G)$ of $T_G$, namely
$\D(G) \le \D(T_G) \le 2\D(G)$,
which implies $\Theta(\D(T_G)) = \Theta(\D(G))$.
We then apply the algorithm from the above 
subsection\tume\ref{refinement-tree} to this BFS spanning tree $T_G$. 
%% board graph $G$が２連結成分を含む場合にも、任意の頂点から１回幅優先探索を掛けて幅優先探索木を作れば、$G$の直径の高々２倍の直径を有する全域木$T_G$が得られます。
%% $G$と$T_G$の直径のオーダーは変わりませんので、これからは、$T_G$上での直径$D(T_G)$を$D(G)$とみなして、議論を進めます。
%% 先に示した「board graph $G$が木である場合のアルゴリズム」を$T_G$に適用することを考えます。

As already discussed in the previous section, the BFS spanning tree $T_G$ of $G$ may contain 
an isthmus of size at least $k+1$ in the subtrees that span the $2$-connected components of $G$.
Thus, continuing the same line of reasoning as in the previous section, we again encounter a potential issue:
the candidates for $H(\ell)$ may contain an isthmus of size at least $N-n$.
In this case, we fix such a candidate and denote it by $H(\ell)$.
This $H(\ell)$ contains exactly one isthmus $I(\ell)$ of size at least $N - n$, which lies along the subpath 
of the path $P(\ell, t)$ connecting the root $\JoiV(X(\ell))$ of $X(\ell)$ and the root $\JoiV(X(t))$ of $X(t)$.
If it is possible to augment $H(\ell)$ by adding appropriate edges from
$G \setminus T_G$ so as to reduce its maximum isthmus size and make
the puzzle on the resulting augmented subgraph feasible, then we simply
add such edges to $H(\ell)$ and solve the resulting sub-puzzle.
If even this is not possible, then the graph $G$ must contain a chordless cycle $C$ of length at most $\CL(G)$
that passes through the isthmus $I(\ell)$.

In this case, we temporarily add such a cycle $C$ to the current
tree structure.  Whenever the progress along $P(\ell,t)$ is blocked
by the long isthmus, we rotate the pebbles along $C$ only as far as
is necessary to perform the corresponding local exchange on the
active path.  Immediately after that exchange, we execute the cycle
rotation in reverse.  Thus the net effect of the repair is confined
to the active path: every vertex used only by the detour cycle is
restored to the same pebble or blank as before the repair.

Consequently, the detour may temporarily use vertices belonging to
completed or temporarily inactive parts of the tree decomposition
without changing their configurations after the repair.  Hence these
cycle repairs preserve the recursive state of
$\Alg.\Extraction(T,\BT,g^{\prime},t)$ and do not interfere with
previously completed target subtrees.

A repair performed at a vertex $v$ costs
\[
    \Ord\bigl(\min\{n,\CL(G,v)\}\bigr)
    \subseteq
    \Ord\bigl(\min\{n,\CL(G)\}\bigr).
\]
During the construction associated with one $X(\ell)$, this local
repair is repeated along at most $\Ord(k)$ successive positions of
the active path.  Therefore the total preparation cost for one
$X(\ell)$ is
\[
    \Ord\bigl(k\min\{n,\CL(G)\}\bigr),
\]
as in the proof of Theorem\tume\ref{main2}.
%In this case, we temporarily consider the structure obtained by adding $C$ to the tree $T_G$.
%By simultaneously shifting the pebbles along the cycle $C$, we can reposition the pebbles on $I(\ell)$ 
%around the vertices of degree at least three in $T \cup C$, thereby enabling the construction of a feasible 
%sub-puzzle $H'(\ell)$.
%The process of constructing $H'(\ell)$ by collectively shifting the pebbles along $C$ requires at most 
%$\Ord(k \cdot \min\{\CL(G), n\})$ moves per instance of $X(\ell)$.
%This reasoning is entirely analogous to the argument used in the proof of Theorem\tume\ref{main2} 
%in Section\tume\ref{Main2Proof}. 
% このとき、唯一の問題点は、$H(\ell)$を作る際に、$N-n$以上のサイズのisthmusを持つものしか得られない可能性があることです。
% このとき$T_G$は、 $X(\ell)$の付け根$\JoiV(X(\ell))$と$X(t)$の付け根$\JoiV(X(t))$とを結ぶパス$P(\ell,t)$の部分パスとして、サイズ$N-n$以上のisthmus $I(\ell)$を持ちます。
% $H(\ell)$に（$T_G$以外の）$G$の適切な辺を付け加えてisthmusを短くできるのであれば、そのような（$T_G$の外の）辺を$H(\ell)$に追加して、$H(\ell)$をfeasibleにしてから、そのサブパズルを解けばよいです。
% それすらできない場合には、$G$には、このisthmus $I(\ell)$を含む長さ$\CL(G)$以下のコードレスサイクル$C$が存在しますので、$C$上でpebble達を一斉に動かすことにより、$I(\ell)$上のpebble達を$G$の次数３以上の頂点の周りに持ってくると、feasibleなサブパズル$H^{\prime}(\ell)$を作れます。
% この$C$上でpebble達を一斉に動かして$H^{\prime}(\ell)$を構成する作業は、１つの$X(\ell)$について、高々$\Ord((N-n)\min\{\CL(G),n\})$回のmoveを必要とします。
% これは、Theorem 4の証明における議論と全く同じ理屈です。

The total number of moves involved in the sequence of operations can be calculated by focusing on 
one of the $\Ord(k)$ pebbles, denoted by $p$, out of the $|X(t)|$ pebbles. Specifically:
\begin{itemize}
\item When pebble $p$ is placed on $X(\ell)$, a series of preparation moves, as described earlier, are required to construct $H^{\prime}(\ell)$. 
Afterward, solving the subproblem on $H^{\prime}(\ell)$ involves traversing a shortest path on the graph $G$ from the initial position of $p$ 
(the vertex $f^{-1}(p)$) to its goal position on $X(t)$ (the vertex ${g^{\prime}}^{-1}(p)$). 
Each step in this path moves $\Ord(k)$ units towards the goal position ${g^{\prime}}^{-1}(p)$, and this process is repeated.
\item For each pebble, the total number of subproblems to be solved is at most $\Ord(\D(G)/k)$, leading to the goal vertex.
\end{itemize}
Therefore, when aggregating the total effort for all $\Ord(k)$ pebbles on $X(t)$, the overall number of moves required 
in this process can be expressed as:
$$
\Ord(k \min\{\CL(G),n\}) \cdot \Ord(\D(G)/k) \cdot \Ord(k) = \Ord(k\min\{\CL(G),n\}\D(G)).
$$
This quantity should be regarded as an overhead added to the cost
obtained under the assumption that $G$ is a tree.

Since there are $\Ord(n/k)$ candidates for $X(t)$, 
the total overhead of the algorithm is therefore
$$\Ord(k\min\{\CL(G),n\}\D(G)) \cdot \Ord(n/k) = \Ord(n\min\{\CL(G),n\}\D(G)).$$

% 一連の作業により掛かるmoveは、$|X(t)|$個（$\Ord(N-n)$個）のpebbleの一つ $p$ に着目して計算すると、
% ・$p$が$X(\ell)$上に置かれている場合に、$H^{\prime}(\ell)$を構成するために上記の（準備の）手間のmoveが発生し、その後$H^{\prime}(\ell)$上のサブパズルを解くことで、$p$のinitial positionとしての$G$の頂点$f^{-1}(p)$と、$p$のgoal positionとなる$X(t)$上の頂点${g^{\prime}}^{-1}(p)$を結ぶ$G$上の（長さ$\Ord(D(G))$の）最短パスの上をgoal position ${g^{\prime}}^{-1}(p)$に向かって$\Ord(N-n)$ずつ歩を進める、という一連の作業が繰り返される
% ので、
% ・pebble１個あたり、トータルとして高々$\Ord(D(G)/(N-n))$回のサブパズルを解いてゴールとなる頂点上に行き着く
% 訳です。
% よって$X(t)$上のすべてのpebble（は$\Ord(N-n)$個ある）の手間を合算すると、この一連の過程で、
% $\Ord((N-n)\min\{\CL(G),n\}(D(G)/(N-n))(N-n))=\Ord((N-n)\min\{\CL(G),n\}D(G))$のmoveの手間が（board graph Gが木の場合の見積に加えて）余分に発生する、と考えることが出来ます。
% さらに、ターゲットである$X(t)$の候補は、全部で$\Ord(n/(N-n))$個あるので、この一連の過程は高々$\Ord(n/(N-n))$回繰り返されます。
% したがってこの余分な作業は、アルゴリズム全体として、$\Ord((N-n)\min\{\CL(G),n\}D(G)(n/(N-n)))=\Ord(n\min\{\CL(G),n\}D(G))$のmoveであると見積もることが出来ます。

In the previous subsection (Subsection\tume\ref{refinement-tree}), we established 
that the length of the solution sequence is 
$$
\Ord(n\D(G) + \min\{kn\D(G), n^2 \log(1+\min\{n, k\})\})
$$
when the board graph $G$ is a tree.
By adding the order of the computational cost for the preparatory moves estimated above, 
we obtain the following order for the move operations:
$$
\Ord(n\D(G) + n \min\{\CL(G), n\} \D(G) + \min\{kn \D(G), n^2 \log(1+\min\{n, k\})\}).
$$
Although similar remarks have been repeatedly made before, all occurrences of the parameter $k$ 
above can also be simultaneously replaced by any fixed integer $\theta$ satisfying $k \leq \theta \leq N - n - 1$, 
without invalidating the statement.
By comparing this result with the one derived in Theorem\tume\ref{main2} and taking the minimum of 
both orders, we obtain the overall time complexity as follows:
\begin{thm}\label{refinement2}
Assume that a shortest cycle passing through each vertex $v$ of $G$ (or $\emptyset$ if no such cycle exists) 
is provided in advance.
Let $\theta$ be an arbitrary integer such that $k \leq \theta \leq N-n-1$.
Then, for the \textsc{Pebble Motion Problem on General Graphs}, a solution sequence of length
$\Ord(n\D(G) + n \min\{\CL(G), n\} \min\{\D(G), \frac{n}{\theta}\} + \min\{n\theta\D(G), n^2 \log(1+\min\{n, \theta\})\})$
can be computed in time 
$\Ord(|E(G)| + n\D(G)+n \min\{\CL(G), n\} \min\{\D(G), \frac{n}{\theta}\} + \min\{n\theta\D(G), n^2 \log(1+\min\{n, \theta\})\}).$
\end{thm}

Theorem\tume\ref{refinement2-general-constant-blanc}, stated earlier, follows 
immediately as a corollary of Theorem~\ref{refinement2}.
\begin{Proofof}{Theorem\tume\ref{refinement2-general-constant-blanc}}
Note that $k \leq \theta < N-n \in \Ord(1)$ and hence $\{\CL(G), \D(G)\} \subset \Ord(n)$ holds in this setting.  
As a consequence, the time complexity characterization in Theorem\tume\ref{refinement2} can be simplified as 
$\Ord(|E(G)| + n\D(G)+ n \CL(G) \D(G) + n\D(G))=\Ord(|E(G)| + n \CL(G) \D(G)). $
Similarly, the solution sequence length also simplifies to $\Ord(n \CL(G) \D(G))$.
\end{Proofof}

The preceding proof also gives a more local accounting of the cycle
repairs.  Consider the processing of one target subtree $X(t)$.
The relevant pebbles progress along the active paths toward $X(t)$.
Whenever the next portion of such a path is blocked by an isthmus
that is too long for the available blank space, choose a vertex $v$
of that obstructing portion and borrow a shortest cycle of $G$
through $v$.  Perform the necessary cycle rotation, make the local
exchange on the active path, and immediately reverse the rotation.
Thus this repair changes the required local state on the active path
while restoring every vertex used only by the detour cycle.

For the reader's convenience, recall that $V(G,q)$ denotes the set of
degree-two vertices of $G$ that belong to a contiguous sequence of at
least $q$ degree-two vertices.  Only vertices in
\[
    V(G,N-n-2)
\]
can occur in the portions of $T_G$ that require such a cycle repair.
At a path position corresponding to $v$, the repair costs
$\Ord(\CL(G,v))$.  The active paths are processed monotonically toward
the target subtree, so the repair cost for one processing of $X(t)$
is bounded, by charging the possible repair at each path position, by
\[
    \Ord\!\left(
        \sum_{v\in V(G,N-n-2)} \CL(G,v)
    \right).
\]
This is simply the path-length accounting used above, with the uniform
bound $\CL(G)$ replaced by the actual shortest-cycle length
$\CL(G,v)$ at each position.

There are $\Ord(n/k)$ target subtrees.  Taking
$k=N-n$, their number is
\[
    \Ord\!\left(\frac{n}{N-n}\right).
\]
Hence the total cycle-repair overhead is
\[
    \Ord\!\left(
        \frac{n}{N-n}
        \sum_{v\in V(G,N-n-2)} \CL(G,v)
    \right).
\]
Adding the tree-routing terms gives the following theorem.
%In the proof of Theorem~\ref{refinement2}, we estimated the overhead required to eliminate the infeasibility 
%caused by an excessively long isthmus in $T_G$.
%More precisely, when such an isthmus appears, we consider the shortest cycle in $G$ containing this isthmus 
%and restore feasibility by simultaneously shifting all pebbles along the cycle by a distance equal to the length 
%of the isthmus.
%
%This overhead also admits another upper bound.
%For the reader's convenience, we recall the definition of $V(G, q)$, which denotes the set of degree-two vertices in $G$ 
%that lie on paths in $G$ and are part of a contiguous sequence of at least $q$ degree-two vertices.
%Vertices in $V(G, N-n-2)$ may create overly long isthmuses in the spanning tree $T_G$, serving as sources of infeasibility for the PMT.  
%In the worst case, every vertex in this set induces such an isthmus in $T_G$. 
%Consequently, the overhead required to move pebbles to their designated positions on a single target subtree $X(t)$ is bounded by
%$\Ord\!\left(\sum_{v \in V(G, N-n-2)} \CL(G, v)\right)$.
%Since the number of target subtrees is generally $\Ord(n/k)$,
%setting $k = N-n$ gives that the number of target subtrees is $\Ord(n/(N-n))$.
%Multiplying this bound by the number of target subtrees then gives the following estimate for the total overhead of the algorithm:
%$\Ord\!\left(\frac{n}{N-n} \sum_{v\in V(G,N-n-2)} \CL(G,v)\right)$.
%From this observation, we obtain the following theorem as a corollary.
\begin{thm}\label{refinement3}
Let $(G,P)$ be a feasible puzzle. For every pair of configurations $f$ and
$g$ of $P$ on $G$, the length of a shortest solution sequence from $f$ to
$g$ is in
$$\Ord\Bigl(n\D(G) + \frac{n}{N-n}\sum_{v\in V(G,N-n-2)} \CL(G,v) + n^2 \log(1+\min\{n, N-n\})\Bigr).$$
\end{thm}

Theorem\tume\ref{refinement3} immediately yields the following
necessary conditions for cubic complexity.

\begin{thm}\label{Restriction}
Let $(G,P)$ be a feasible puzzle on a general board graph.
If the worst-case shortest solution length of the puzzle $(G,P)$ is
$\Theta(N^3)$, then the board graph $G$ satisfies the following conditions:
\begin{enumerate}
\item\label{res:constant-blanks}
$N-n=\Ord(1)$;

\item\label{res:CLL}
$\displaystyle
\sum_{v\in V(G,N-n-2)}\CL(G,v)=\Theta(N^2)$;

\item\label{res:CLL3}
$G$ contains a set
$S\subseteq V(G,N-n-2)$ of $\Theta(N)$ vertices;

\item\label{res:CLL4}
for every $v\in S$, the length of a shortest cycle in $G$ passing
through $v$ is $\Theta(N)$.
\end{enumerate}
\end{thm}

\begin{Proofof}{Theorem~\ref{Restriction}}
Conditions~\ref{res:constant-blanks} and \ref{res:CLL} follow
immediately from Theorem~\ref{refinement3}.
Since both $\CL(G)$ and $|V(G,N-n-2)|$ are in $\Ord(N)$,
Condition~\ref{res:CLL} implies
Conditions~\ref{res:CLL3} and \ref{res:CLL4}.
\end{Proofof}

\noindent
This observation shows that the $\Theta(N^3)$ behavior of the PMP
can arise only under highly restrictive structural conditions on the
instance $(G,P)$.

%%%%%%%%%%%%%%%%%%%%%%%%%%%%%%%%%%%%%%%%%%%%%%%%%%%%%%%%%%%%%%%%%%%%%%%%%%%%%%%%
%%%%%%%%%%%%%%%%%%%%%%%%%%%%%%%%%%%%%%%%%%%%%%%%%%%%%%%%%%%%%%%%%%%%%%%%%%%%%%%%
%%%%%%%%%%%%%%%%%%%%%%%%%%%%%%%%%%%%%%%%%%%%%%%%%%%%%%%%%%%%%%%%%%%%%%%%%%%%%%%%

\section{Pebble Presentation Algorithm}
\label{sec:pebble-presentation}

The next several sections investigate structural restrictions on feasible
puzzles whose worst-case shortest solution length is of order
$\Theta(N^3)$.  By Theorem\tume\ref{Restriction}, such puzzles have only
$\Ord(1)$ unoccupied vertices.  Accordingly, throughout this section
we write
\[
    b:=N-n
\]
and assume
\[
    b=\Ord(1).
\]

The purpose of this section is to establish the \textit{Pebble
Presentation Algorithm}.  Given a connected board graph $G$ and a fixed
vertex $u$, the algorithm makes every pebble visit $u$ exactly once.  It
will be used in Section~\ref{sec:structural-restrictions} to move
prescribed sets of pebbles through cut vertices and into ears.

%%%%%%%%%%%%%%%%%%%%%%%%%%%%%%%%%%%%%%%%%%%%%%%%%%%%%%%%%%%%%%%%%%%%%%%%%%%%%%%%

\subsection{Biconnected--bridge-tree decomposition}
\label{subsec:biconnected-bridge-tree-decomposition}

We first introduce the decomposition used by the algorithm.

A \textit{biconnected component} is a nontrivial maximal biconnected
subgraph of $G$.  A cut vertex that belongs to a biconnected component
and is also incident with a bridge is called an \textit{attachment
vertex}.

Let $B(G)$ be the forest consisting of all bridges of $G$ and their
endpoints, and let $A$ be the set of attachment vertices.  Delete $A$
from $B(G)$.  For each nonempty connected component $C$ of $B(G)-A$,
let $\bar C$ be the subtree obtained from $C$ by adding back every vertex
of $A$ adjacent to $C$ in $B(G)$ together with the corresponding bridge.
We call $\bar C$ a \textit{bridge-tree component}.  If a bridge has both
endpoints in $A$, that bridge together with its two endpoints is itself
regarded as a bridge-tree component.

The collection of all biconnected components and all bridge-tree
components is denoted by
\[
    \mathcal{B}(G)
\]
and is called the \textit{biconnected--bridge-tree decomposition} of
$G$.  Its members are edge-disjoint, although two or more members may
share a cut vertex.  The corresponding component--attachment incidence
graph is the bipartite graph whose vertices are the members of
$\mathcal{B}(G)$ together with the attachment vertices, with a component
adjacent to each attachment vertex it contains.  This incidence graph is
a tree.

The members of $\mathcal{B}(G)$ are structural units only; a member need
not be the board graph of a feasible puzzle when considered in isolation.
All moves in
what follows are moves of the original graph $G$.  Accordingly, when we
refer to a reconfiguration on a component $X\in\mathcal{B}(G)$, we mean
the moves required to complete the prescribed configuration on $X$
within the original graph $G$, rather than a PMP restricted to $X$ as
an isolated board graph.  When necessary, the local working region may
temporarily use, at each cut vertex of $X$, at most two incident edges
of $G$ lying outside $X$, together with their adjacent vertices and the
pebbles occupying them. Since there are at most $|V(X)|$ cut vertices
in $X$, this enlarges the working region by only $\Ord(|V(X)|)$ vertices
and therefore does not change its asymptotic order.
Any vertices outside $X$ that are temporarily included in the working
region, together with the pebbles occupying them, are restored to their
original configuration when the reconfiguration of $X$ is completed.
We call the corresponding worst-case number of moves the
\textit{reconfiguration cost of $X$}.
%The members of $\mathcal{B}(G)$ are structural units only; a member need
%not be a feasible PMP board when considered in isolation.  All moves in
%what follows are moves of the original graph $G$.

%%%%%%%%%%%%%%%%%%%%%%%%%%%%%%%%%%%%%%%%%%%%%%%%%%%%%%%%%%%%%%%%%%%%%%%%%%%%%%%%

\subsection{The Pebble Presentation Algorithm}

Fix a vertex $u$, called the \textit{interface vertex}.  The algorithm
will bring every pebble to $u$ exactly once.  To record its progress, we
color every pebble blue initially.  When a blue pebble is presented at
$u$, it is recolored red.  Thus blue means ``not yet presented'' and red
means ``already presented.''  After a pebble is presented, a storage
subroutine moves that red pebble away from the current interface so that
a blue pebble can subsequently reach it.  The storage subroutines never
change colors.

\begin{alg}[
$\Alg.\mathrm{PebblePresentation}
(G,\mathcal{B}(G),\sigma,X,u)$]\hfill
\label{alg:pebble-presentation}

\begin{namelist}{OUTPUT}

\item[\sl Input]
A connected graph $G$, its biconnected--bridge-tree decomposition
$\mathcal{B}(G)$, a configuration $\sigma$ of a pebble set $P$ such that
the puzzle $(G,P)$ is feasible, a component $X\in\mathcal{B}(G)$, and an
interface vertex $u\in V(X)$.

\item[\sl Output]
A sequence of moves starting from $\sigma$ in which every pebble on $G$
is presented at $u$ exactly once.

\end{namelist}

\begin{alg-enumerate}

\item
\label{alg:ppa-step1}
Color every pebble blue.

\item
\label{alg:ppa-step2}
If a blue pebble occupies $u$, choose it.  Otherwise move a blue pebble
to $u$.  Denote the chosen pebble by $p$.

\item
\label{alg:ppa-step3}
Present $p$ at $u$ and recolor it red.

\item
\label{alg:ppa-step4}
If every pebble is red, return the sequence of moves performed so far.
Otherwise, if $X$ is a bridge-tree component, call
\[
    \Alg.\mathrm{BridgeTreeStorage}
    (G,\mathcal{B}(G),X,u,p),
\]
and if $X$ is a biconnected component, call
\[
    \Alg.\mathrm{BiconnectedStorage}
    (G,\mathcal{B}(G),X,u,p).
\]

\item
\label{alg:ppa-step5}
Return to Step~\ref{alg:ppa-step2}.

\end{alg-enumerate}
\end{alg}

%%%%%%%%%%%%%%%%%%%%%%%%%%%%%%%%%%%%%%%%%%%%%%%%%%%%%%%%%%%%%%%%%%%%%%%%%%%%%%%%

\subsection{Storage in a bridge-tree component}

Suppose that the current component $X$ is a bridge-tree component, that
$u\in V(X)$ is the current interface, and that a red pebble $p$ occupies
$u$.  The purpose of the storage operation is simply to move $p$ away
from $u$ through the part of $X$ that still contains blue pebbles or
blanks.  The tree operations needed here are the same ones used in the
proof of Theorem\tume\ref{refinement1}.

\begin{alg}[
$\Alg.\mathrm{BridgeTreeStorage}
(G,\mathcal{B}(G),X,u,p)$]\hfill
\label{alg:bridge-tree-storage}

\begin{namelist}{OUTPUT}

\item[\sl Input]
A bridge-tree component $X\in\mathcal{B}(G)$, an interface vertex
$u\in V(X)$, and a red pebble $p$ occupying $u$.

\item[\sl Output]
A sequence of moves that moves $p$ away from $u$, without changing the
color of any pebble.

\end{namelist}

\begin{alg-enumerate}

\item
\label{alg:bridge-storage-step1}
Let $Y$ be the subtree of $X$ induced by $u$ together with all vertices
that are unoccupied or occupied by blue pebbles.  Choose a leaf $r$ of
$Y$ at maximum distance from $u$.

\item
\label{alg:bridge-storage-step2}
Move the red pebble $p$ from $u$ to $r$, using the tree reconfiguration
operations established earlier.

\item
\label{alg:bridge-storage-step3}
If $r$ is a boundary attachment vertex and there exists another
component in $\mathcal{B}(G)$ containing $r$ that still contains a blue
pebble, choose one such component $Z$.  Call
$\Alg.\mathrm{BridgeTreeStorage}$ or
$\Alg.\mathrm{BiconnectedStorage}$ on $Z$ with interface $r$, according
as $Z$ is a bridge-tree or biconnected component.
%\item
%\label{alg:bridge-storage-step3}
%If $r$ is a boundary attachment vertex and another component
%$Z\in\mathcal{B}(G)$ containing $r$ still contains a blue pebble, call
%$\Alg.\mathrm{BridgeTreeStorage}$ or
%$\Alg.\mathrm{BiconnectedStorage}$ on $Z$ with interface $r$, according
%as $Z$ is a bridge-tree or biconnected component.

\item
\label{alg:bridge-storage-step4}
Return the sequence of moves performed by the subroutine.

\end{alg-enumerate}
\end{alg}

%%%%%%%%%%%%%%%%%%%%%%%%%%%%%%%%%%%%%%%%%%%%%%%%%%%%%%%%%%%%%%%%%%%%%%%%%%%%%%%%

\subsection{Storage in a biconnected component}

We next describe the only new local motion used by the algorithm.  If a
blank is moved once around a cycle, then every pebble on that cycle is
shifted by one position in the direction opposite to the motion of the
blank.  We call this operation a \textit{cycle rotation}.  The blank
that travels around the cycle will be called the \textit{active blank}.
The storage procedure in a biconnected component uses one active blank
to move an already presented pebble away from the interface and make
room for the remaining blue pebbles.

Let $X$ be a biconnected component, let $u\in V(X)$ be the current
interface, and let $p$ be a red pebble occupying $u$.

\begin{alg}[
$\Alg.\mathrm{BiconnectedStorage}
(G,\mathcal{B}(G),X,u,p)$]\hfill
\label{alg:biconnected-storage}

\begin{namelist}{OUTPUT}

\item[\sl Input]
A biconnected component $X\in\mathcal{B}(G)$, an interface vertex
$u\in V(X)$, and a red pebble $p$ occupying $u$.

\item[\sl Output]
A sequence of moves that moves a red pebble away from the current
interface, without changing the color of any pebble.

\end{namelist}

\begin{alg-enumerate}

\item
\label{alg:biconnected-storage-step1}
Choose an edge $e=uv$ of $X$ incident with $u$.  Orient $X-e$ acyclically
with $v$ as its unique source and $u$ as its unique sink, so that every
vertex lies on a directed path from $v$ to $u$.  Orient $e$ from $u$ to
$v$.  Denote the resulting orientation of $X$ by $D(X)$.

\item
\label{alg:biconnected-storage-step2}
Choose one blank as the active blank and move it to $v$.  Place every
other blank, if any, in a connected region adjacent to $u$, obtained by
exploring $D(X-e)$ backward from $u$ in breadth-first order.  Do this so
that, for each such blank $z$, a shortest directed path from $z$ to $u$
contains no blue pebble.  Equivalently, along every directed path toward
$u$ containing only blue pebbles and blanks, the blanks form a terminal
segment adjacent to $u$.  Return $p$ to $u$ if it was displaced while
the active blank was moved to $v$.

\item
\label{alg:biconnected-storage-step3}
Let $H$ be the subgraph of $D(X-e)$ induced by $u$ together with all
vertices $x$ from which $u$ is reachable by a directed path whose
vertices other than $u$ are all unoccupied or occupied by blue pebbles.

Choose a directed path $P$ in $H$ ending at $u$ that cannot be extended
backward in $H$ at its initial vertex, and let $w$ be the initial vertex
of $P$.  Choose a shortest directed path $Q$ from $u$ to $w$ in $D(X)$.
Move the active blank once around the directed cycle $P\cup Q$ in the
direction opposite to its orientation.

If $P$ contains a blue pebble and one or more inactive blanks lie
consecutively between the first blue pebble on $P$ and $u$, move that
blue pebble forward along $P$ through these blank vertices until it
reaches $u$.  No other blue pebble is moved in this final sliding
operation.  Since there are only $b=\Ord(1)$ blanks altogether, this
requires only $\Ord(b)=\Ord(1)$ additional moves.
%If one or more inactive blanks lie consecutively between the first blue
%pebble on $P$ and $u$, move that blue pebble forward along $P$ through
%these blank vertices until it reaches $u$.  No other blue pebble is
%moved in this final sliding operation.  Since there are only
%$b=\Ord(1)$ blanks altogether, this requires only $\Ord(b)=\Ord(1)$
%additional moves.
%\item
%\label{alg:biconnected-storage-step3}
%Let $H$ be the subgraph of $D(X-e)$ induced by $u$ together with all
%vertices that are unoccupied or occupied by blue pebbles.  Choose a
%directed path $P$ in $H$ ending at $u$ that cannot be extended backward
%at its initial vertex, and let $w$ be the initial vertex of $P$.
%Choose a shortest directed path $Q$ from $u$ to $w$ in $D(X)$.
%Move the active blank once around the directed cycle $P\cup Q$ in the
%direction opposite to its orientation.

\item
\label{alg:biconnected-storage-step4}
If this rotation places a red pebble at a cut vertex $r\ne u$ and there
exists another component in $\mathcal{B}(G)$ containing $r$ that still
contains a blue pebble, choose one such component $Z$.  Call
$\Alg.\mathrm{BridgeTreeStorage}$ or
$\Alg.\mathrm{BiconnectedStorage}$ on $Z$ with interface $r$, according
as $Z$ is a bridge-tree or biconnected component.
%\item
%\label{alg:biconnected-storage-step4}
%If this rotation places a red pebble at a cut vertex $r\ne u$ and
%another component $Z\in\mathcal{B}(G)$ containing $r$ still contains
%a blue pebble, call $\Alg.\mathrm{BridgeTreeStorage}$ or
%$\Alg.\mathrm{BiconnectedStorage}$ on $Z$ with interface $r$, according
%as $Z$ is a bridge-tree or biconnected component.

\item
\label{alg:biconnected-storage-step5}
Return the sequence of moves performed by the subroutine.

\end{alg-enumerate}
\end{alg}

%%%%%%%%%%%%%%%%%%%%%%%%%%%%%%%%%%%%%%%%%%%%%%%%%%%%%%%%%%%%%%%%%%%%%%%%%%%%%%%%

\subsection{Correctness and complexity}

We first justify the orientation used in
Step~\ref{alg:biconnected-storage-step1}.
Since $X$ is biconnected and $uv\in E(X)$, an $st$-numbering of $X$
with $v$ first and $u$ last exists.
Orient every edge of $X-e$ from the smaller to the larger vertex in
this numbering.
The resulting orientation of $X-e$ is acyclic, has $v$ as its unique
source and $u$ as its unique sink, and every vertex lies on a directed
path from $v$ to $u$.
Finally orient $e$ from $u$ to $v$, as prescribed in
Step~\ref{alg:biconnected-storage-step1}.

We next record why the cycle in Step~\ref{alg:biconnected-storage-step3} 
has the required form. Since $u$ is the unique sink of $D(X-e)$, every 
directed path leaving $u$ in $D(X)$ starts with the directed edge $uv$.  
Hence $Q$ contains $uv$. Moreover, $P$ and $Q$ have no common vertex 
other than $u$ and $w$; otherwise their subpaths between two common 
vertices would form a directed cycle in the acyclic digraph $D(X-e)$. 
Consequently $P\cup Q$ is a directed cycle containing $uv$.

The definition of $H$ records precisely the part of the current
biconnected component that is still accessible from the interface
through blue pebbles and blanks.  Throughout the procedure we maintain
the following invariant: if a red pebble occupies a vertex $r\ne u$,
then every ancestor of $r$ in $D(X-e)$ other than the source vertex
$v$ is also occupied by a red pebble.  
The vertex $v$ is excluded
because the active blank is placed at $v$ in
Step~\ref{alg:biconnected-storage-step2}.

This invariant is preserved by
Step~\ref{alg:biconnected-storage-step3}.  By the backward maximality
of $P$, every in-neighbor of its initial vertex $w$ in $D(X-e)$,
except possibly $v$, is occupied by a red pebble; otherwise that
vertex would belong to $H$ and $P$ could be extended backward.
Consequently, along the directed path $Q$ from $v$ to $w$, the
invariant implies successively that every vertex other than $v$ and
$w$ is occupied by a red pebble.  On the other hand,
$P-\{u\}\subseteq H$, so its vertices are unoccupied or occupied by
blue pebbles.  Hence the cycle rotation preserves the invariant.
The subsequent sliding of the first blue pebble through the terminal
block of blanks changes no red pebble and therefore preserves it as
well.

If the current presentation frontier reaches a cut vertex $r$ and 
an adjacent component $Z$ beyond $r$ still contains a blue pebble, 
Step~\ref{alg:biconnected-storage-step4} immediately invokes 
the corresponding storage subroutine on $Z$ with interface $r$. 
Thus the presentation frontier does not pass beyond $r$ 
until all blue pebbles in $Z$ have been processed.

Consequently, every blue pebble that remains in the current
biconnected component belongs to $H$.  Indeed, since $u$ is the unique
sink of $D(X-e)$, every vertex has a directed path to $u$.  If a
directed path from a blue vertex $x$ to $u$ met a red vertex
$r\ne u$, then $x$ would be an ancestor of $r$.  Since $x$ is occupied
by a blue pebble, we have $x\ne v$, and the invariant would therefore
imply that $x$ itself is red, a contradiction.  Thus $x$ reaches $u$
through vertices that are all blue or unoccupied, except possibly for
$u$ itself.
%The definition of $H$ records precisely the part of the current
%biconnected component that is still accessible from the interface
%through blue pebbles and blanks.  Throughout the procedure we maintain
%the following invariant: if a red pebble occupies a vertex $r\ne u$,
%then every ancestor of $r$ in $D(X-e)$ is also occupied by a red pebble.
%Equivalently, no blue pebble lies upstream of a red pebble distinct
%from the current interface pebble.
%
%Consequently, every blue pebble that remains in the current
%biconnected component belongs to $H$.  Indeed, since $u$ is the unique
%sink of $D(X-e)$, every vertex has a directed path to $u$.  If a
%directed path from a blue vertex $x$ to $u$ met a red vertex
%$r\ne u$, then $x$ would be an ancestor of $r$ and hence would itself
%be red by the invariant, a contradiction.  Thus $x$ reaches $u$
%through vertices that are all blue or unoccupied, except possibly for
%$u$ itself.

The backward maximality of $P$ gives the required separation at the
presentation frontier.  There is no directed path
\[
    x_0\to x_1\to\cdots\to x_k=w,
    \qquad k\ge 1,
\]
in $D(X-e)$ for which $x_0$ contains a blue pebble and the intervening
vertices $x_1,\ldots,x_{k-1}$, if any, are all blank.  Such a path
would lie in $H$ and could be prepended to $P$, contradicting the
choice of $P$.  Hence a red pebble that has already been moved beyond
the current presentation frontier cannot subsequently be reintroduced
among the blue pebbles transported toward $u$.
%The backward maximality of $P$ also gives the separation needed by the
%presentation procedure.  There is no directed path
%\[
%    x_0\to x_1\to\cdots\to x_k=w,
%    \qquad k\ge 1,
%\]
%in $D(X-e)$ for which $x_0$ contains a blue pebble and the intervening
%vertices $x_1,\ldots,x_{k-1}$, if any, are all blank.  Thus a blue
%pebble cannot reach $w$ along a directed path with zero or more blanks
%between it and $w$.  Such a path would lie in $H$ and could be prepended
%to $P$, contradicting the choice of $P$.  Hence a red pebble that has
%been moved past $w$ cannot subsequently be reintroduced among the blue
%pebbles transported toward $u$.

\begin{thm}[Pebble Presentation Theorem]
\label{thm:pebble-presentation}
Let $(G,P)$ be a feasible puzzle on a connected $N$-vertex board graph
$G$, and let $\sigma$ be a configuration of the $n=|P|$ pebbles on $G$.
Suppose that
\[
    2\leq N-n=\Ord(1).
\]
Then, for every vertex $u\in V(G)$, there exists a move sequence of
length
\[
    \Ord(N^2)
\]
in which every pebble is presented at $u$ exactly once.

Moreover, the Pebble Presentation Algorithm constructs such a sequence
in $\Ord(N^2)$ computational time, whether or not
$\mathcal{B}(G)$ is supplied as part of the input.
\end{thm}

\begin{Proofof}{Theorem\tume\ref{thm:pebble-presentation}}
Each iteration of Algorithm~\ref{alg:pebble-presentation} presents
exactly one blue pebble at $u$ and recolors it red.  We first verify
that, whenever a blue pebble remains, another blue pebble can indeed
be brought to the interface.

In a bridge-tree component this follows from the construction of the
subtree induced by the interface vertex together with vertices that
are unoccupied or occupied by blue pebbles.  In a biconnected
component, the invariant established above implies that every remaining
blue pebble belongs to the subgraph $H$ of
Step~\ref{alg:biconnected-storage-step3}.  Hence the storage operation
can advance the presentation frontier until a blue pebble lies
immediately behind the terminal block of blank vertices adjacent to
$u$.  After the single cycle rotation in
Step~\ref{alg:biconnected-storage-step3}, that blue pebble is moved
forward through this block of blanks to $u$.  This final operation moves
only that blue pebble and requires at most $\Ord(b)=\Ord(1)$ moves.

No storage subroutine changes the color of any pebble.  In a
bridge-tree component the presented red pebble is moved away from the
interface through the subtree induced by blue pebbles and blanks; if
the selected leaf is an attachment vertex, storage continues in the
adjacent component.  In a biconnected component, the backward-maximal
choice of $P$ prevents an already stored red pebble from re-entering
the blue presentation frontier.  Hence no pebble is presented twice.

Therefore, whenever a blue pebble remains, one further blue pebble can
be presented at $u$, and exactly that pebble is recolored red.  After
$n$ iterations every pebble has consequently been presented exactly
once.
%Each iteration of Algorithm~\ref{alg:pebble-presentation} selects exactly
%one blue pebble at $u$ and recolors it red.  No storage subroutine changes
%a color.  In a bridge-tree component the red pebble is moved away from
%the interface through the subtree induced by blue pebbles and blanks;
%if the selected leaf is an attachment vertex, storage continues in the
%adjacent component.  In a biconnected component, the maximality property
%proved above prevents a red pebble from re-entering the set of blue
%pebbles that can subsequently be transported to the interface.  Hence a
%pebble is never presented twice.  Since one new pebble is presented in
%each iteration, after $n$ iterations every pebble has been presented
%exactly once.

We next bound the number of moves.  In a biconnected component $X$,
constructing the orientation, finding the paths $P$ and $Q$, positioning
the $\Ord(1)$ blanks, and performing the single cycle rotation require
$\Ord(|V(X)|)$ moves in total.  After that rotation, the next blue
pebble may have to be moved forward through a terminal block of at most
$b-1=\Ord(1)$ blank vertices before reaching $u$.  This requires only
$\Ord(b)=\Ord(1)$ additional moves.  Hence the total contribution of
the biconnected storage operation remains
\[
    \Ord(|V(X)|).
\]
%We next bound the number of moves.  In a biconnected component $X$, the
%orientation, the paths $P$ and $Q$, the positioning of the $\Ord(1)$
%blanks, and one cycle rotation together require $\Ord(|V(X)|)$ moves.
In a bridge-tree component $X$, the tree operations used in
Theorem\tume\ref{refinement1} move the red pebble along the required tree
path with $\Ord(1)$ work per bounded isthmus segment.  Thus the
contribution of $X$ is $\Ord(|V(X)|)$.

During the storage of one presented pebble, recursive calls follow one
path in the component--attachment incidence tree and never revisit a
component.  If the visited components are $X_1,\ldots,X_k$, then
\[
    |V(X_1)|+
    \sum_{i=2}^{k}(|V(X_i)|-1)\le N.
\]
Hence one complete storage operation costs $\Ord(N)$ moves.  There are
at most $N$ presentations, and therefore the total number of moves is
$\Ord(N^2)$.

The same accounting gives the running time.  Bridges, cut vertices,
biconnected components, and the incidence tree can be found in linear
time.  For a biconnected component one may preprocess a spanning
biconnected subgraph with $\Ord(|V(X)|)$ edges, after which the
source--sink orientation and the paths $P,Q$ are found in linear time.
Tree operations are likewise linear in the portion of the tree visited.
Thus one presentation costs $\Ord(N)$ computational time and all
presentations cost $\Ord(N^2)$.
\end{Proofof}

%%%%%%%%%%%%%%%%%%%%%%%%%%%%%%%%%%%%%%%%%%%%%%%%%%%%%%%%%%%%%%%%%%%%%%%%%%%%%%%%
%%%%%%%%%%%%%%%%%%%%%%%%%%%%%%%%%%%%%%%%%%%%%%%%%%%%%%%%%%%%%%%%%%%%%%%%%%%%%%%%
%%%%%%%%%%%%%%%%%%%%%%%%%%%%%%%%%%%%%%%%%%%%%%%%%%%%%%%%%%%%%%%%%%%%%%%%%%%%%%%%

\section{Structural Restrictions for the Cubic Regime}
\label{sec:structural-restrictions}

Throughout this section we retain the standing assumption
\[
    N-n=\Ord(1).
\]
Throughout this section, reconfiguration costs of components are
understood in the sense specified in
Subsection~\ref{subsec:biconnected-bridge-tree-decomposition}.

\paragraph{Asymptotic convention.}
All asymptotic statements in Sections~\ref{sec:structural-restrictions}
and~\ref{lower bound} are understood for families of feasible puzzles
\[
    \{(G_N,P_N)\}_{N\to\infty},
    \qquad |V(G_N)|=N,
\]
and, when relevant, for families of subgraphs, components, or cycles
chosen from them.  We suppress the subscript $N$ whenever no confusion
can arise.  Thus symbols such as $G$, $X$, $Y$, and $C$ in an
asymptotic statement denote $G_N$, $X_N$, $Y_N$, and $C_N$,
respectively; no identification of these objects for different values
of $N$ is intended.  The constants implicit in $\Ord(\cdot)$,
$\Omega(\cdot)$, and $\Theta(\cdot)$ are uniform along the family.
In particular, a universally quantified assertion that components or
attachments have order $\Ord(1)$ means that one common constant bounds
all of them for all sufficiently large $N$.  The assertion that the
worst-case shortest solution length is $\Theta(N^3)$ means that it is
bounded above and below by positive constant multiples of $N^3$ for all
sufficiently large $N$, whereas $\ord(N^3)$ means that its ratio to
$N^3$ tends to zero.  In contradiction arguments we may pass to an
infinite subsequence and relabel it by $N$.

We use the Pebble Presentation Algorithm as a tool for moving prescribed
pebbles through a cut vertex or into an ear.
We first treat a cut vertex that separates a subgraph of order 
$\omega(1)\cap\ord(N)$.  We then analyze the effect of adjoining an ear 
and finally apply these results to the components of $\mathcal{B}(G)$.
Whenever one of the constructions below requires the blank vertices in
a particular unfinished side or local working region, we first move the
required blanks there along paths of the unfinished graph.  Since the
total number of blanks is $\Ord(1)$, one such relocation costs at most
$\Ord(N)$ moves.  These relocation costs will be included implicitly
below; in every application their total is dominated by the stated
subcubic bound.

%%%%%%%%%%%%%%%%%%%%%%%%%%%%%%%%%%%%%%%%%%%%%%%%%%%%%%%%%%%%%%%%%%%%%%%%%%%%%%%%

\subsection{Transcription across a cut vertex}
\label{subsec:transcription-cut}

In this subsection we develop an efficient reconfiguration mechanism
based on the Pebble Presentation Algorithm of the preceding section.
The mechanism is based on a key new idea, which we call
\textit{transcription}.

The significance of the Pebble Presentation Algorithm is not merely
that every pebble on $G$ is presented exactly once at the interface
vertex $u$.  A more important feature is that it induces a total order
on the pairs consisting of an initially occupied vertex and the pebble
placed on that vertex: the order is the order in which those pebbles are
presented at $u$.  Equivalently, every presentation position is associated
with a unique vertex of the initial configuration.

This observation yields the following fundamental replacement principle.
Suppose that a pebble $p$ is presented at $u$ at some specified position
of a recorded presentation sequence.  If $p$ is replaced there by another
pebble $q$, and the recorded moves leading to that presentation position
are then executed in reverse, without undoing the replacement of $p$ by
$q$, the reverse sequence carries $q$ to the vertex initially occupied by
$p$.  All other pebbles affected by that reversed part of the presentation
return to their previous positions.  Thus a substitution made locally at
the interface can be transferred back to the vertex corresponding to the
chosen presentation position.

For later use we record one elementary technical observation.

\begin{la}[Interface clearing]
\label{lem:interface-clearing}
Let a Pebble Presentation Algorithm be running on one side of a cut
vertex $u$, with $u$ as its interface vertex.  The procedure can be
paused with $u$ unoccupied, without changing the order in which the
remaining blue pebbles will subsequently be presented and without
changing the asymptotic cost.
\end{la}

\begin{Proofof}{Lemma\tume\ref{lem:interface-clearing}}
If the current component is biconnected, pause a cycle rotation at the
moment when the pebble at $u$ has moved away and $u$ is blank.  If the
current component is a bridge-tree component, move the red pebble at $u$
away by the same bounded local tree operation used by the storage
subroutine.  Resuming the interrupted operation restores exactly the
state from which the presentation continues.  The extra work is
absorbed in the cost of the corresponding storage operation.
\end{Proofof}

We now describe transcription across a cut vertex.  Let $X$ and $Y$ be
connected subgraphs such that
\[
    V(X)\cap V(Y)=\{u\},
    \qquad
    V(X)\cup V(Y)=V(G),
\]
and suppose that an efficient reconfiguration algorithm $A(X)$ is known
for $X$.  Run Pebble Presentation Algorithms on $X$ and $Y$, both with
interface $u$, pausing either one whenever the other side is operated on.

At every switching point, the inactive blanks are kept in the
constant-size local workspace adjacent to the common interface vertex
$u$, and the active presentation operation can be paused and restored
by the interface-clearing convention.  Since the total number of
blanks is $\Ord(1)$, transferring the required blank workspace from
one side of $u$ to the other requires only $\Ord(1)$ additional moves.
Hence alternating the two presentation procedures introduces no
additional asymptotic cost.

Consider a consecutive interval $S$ of presentation positions in the
recorded presentation sequence of $Y$, whose size does not exceed the
number of pebbles that can be stored on the $X$ side.  Suppose that a
prescribed set $P_S$ of pebbles is to be assigned, in a prescribed
order, to the positions of $S$.  By advancing the presentation
procedures and exchanging pebbles through $u$, the pebbles of $P_S$ can
first be collected on $X$.  The algorithm $A(X)$ can then be used to
reconfigure the $X$ side so that these pebbles are supplied at $u$ in
the prescribed order.  Meanwhile, the recorded presentation on $Y$ can
be returned to the first position of $S$.

Starting from this state, advance the recorded presentation on $Y$
through the positions of $S$.  Whenever the next position of $S$ is
presented at $u$, exchange its current pebble with the next prescribed
pebble supplied from $X$.  Thus the prescribed ordered set $P_S$ is
written into the presentation positions of $S$.  We call this operation
the \textit{transcription} of $S$ from $X$ to $Y$.

The replacement principle above determines what happens when the
recorded presentation on $Y$ is subsequently executed backward: each
pebble written into a presentation position is carried to the vertex
corresponding to that position.  Importantly, the reversal need not be
performed immediately after a transcription.  Several consecutive
intervals may first be transcribed while the presentation state is
carried forward, and the corresponding substitutions may then be
realized by a later backward execution.  This batched use of
transcription is the basis of the next theorem.

%%%%%%%%%%%%%%%%%%%%%%%%%%%%%%%%%%%%%%%%%%%%%%%%%%%%%%%%%%%%%%%%%%%%%%%%%%%%%%%%

\subsection{Subcubic reconfiguration across a cut vertex}

\begin{thm}
\label{thm:intermediate-cut-component}
Let $\{(G_N,P_N)\}_{N\to\infty}$ be a family of feasible
puzzles.  For each $N$, suppose that
$G_N=X_N\cup Y_N$, where $X_N$ and $Y_N$ are connected,
\[
    V(X_N)\cap V(Y_N)=\{u_N\},
\]
and $u_N$ is a cut vertex.  If
\[
    |V(X_N)|=\omega(1)
    \qquad\text{and}\qquad
    |V(X_N)|=\ord(N),
\]
then the family of worst-case shortest solution lengths of the PMP on
$G_N$ is
\[
    \ord(N^3).
\]
\end{thm}

\begin{Proofof}{Theorem\tume\ref{thm:intermediate-cut-component}}
Suppress the subscript $N$ according to the asymptotic convention above.
Put
\[
    f=f(N):=|V(X)|.
\]
Then $f=\omega(1)$ and $f=\ord(N)$, while $|V(Y)|=\Theta(N)$.
Fix arbitrary initial and target configurations.

Since $b=\Ord(1)$ and $|V(X)|=\omega(1)$, for all sufficiently
large $N$ choose a set
\[
    B^{\ast}\subseteq V(X),
    \qquad
    |B^{\ast}|=b,
\]
and put
\[
    S^{\ast}:=V(G)\setminus B^{\ast}.
\]
Let $T_0$ be an arbitrary spanning tree of $G$.
Apply Lemma\tume\ref{sbs0} on $T_0$ to the supports of the initial
and target configurations, in both cases using $S^{\ast}$ as the
target support.
We thereby obtain normalized initial and target configurations whose
blank vertices are exactly $B^{\ast}$.
Each normalization requires
\[
    \Ord(n\D(T_0))=\Ord(N^2)
\]
moves.

It is therefore sufficient to construct an $\ord(N^3)$-move
reconfiguration between these two normalized configurations.
Prepending the normalization sequence for the initial configuration
and appending the reverse of the normalization sequence for the
target configuration changes the total length by only
$\Ord(N^2)$.

Run a Pebble Presentation Algorithm on $Y$ with interface $u$ and record
both its move sequence and the correspondence between presentation
positions and the vertices of $Y$.  Divide the presentation positions,
from the beginning of the presentation order, into consecutive blocks
\[
    S_1,S_2,\ldots,S_m,
\]
each of size $\Theta(f)$ except possibly the last.  Thus
\[
    m=\Ord\!\left(\frac{N}{f}\right).
\]
For each $i$, let $P_i$ be the set of pebbles that must finally occupy
the vertices corresponding to the presentation positions in $S_i$.

We use the transcription mechanism of the preceding subsection in a
batched form.  The essential point is that, after a block has been
written into the presentation sequence of $Y$, we do not execute the
presentation backward to realize that block immediately.  Instead, we
leave the presentation at the end of the block and continue with the
next block.  A single backward execution will be performed only after
all blocks have been written.

We process the blocks in the order
\[
    S_1,S_2,\ldots,S_m.
\]
At the beginning of the processing of $S_i$, maintain the following
invariants.

\begin{enumerate}
\item[(I1)]
Every presentation position belonging to
$S_1,\ldots,S_{i-1}$ has already been assigned its correct target
pebble.

\item[(I2)]
The recorded presentation on $Y$ is in the state immediately after
$S_{i-1}$, equivalently at the boundary immediately preceding $S_i$.
For $i=1$, this is the initial presentation state on $Y$.
\end{enumerate}

The block $S_i$ is processed by the following four operations.

\medskip
\noindent
\textbf{(B1) Collect $P_i$ on $X$.}
Starting from the state in (I2), run the recorded Pebble Presentation
Algorithm on $Y$ forward through its unprocessed suffix.  Whenever a
pebble of $P_i$ that is not already on the $X$ side is presented at
$u$, pause the presentation on $Y$.  Using the interface exchange
described in the preceding subsection, transfer that pebble through
$u$ to $X$, exchanging it with an arbitrary pebble supplied from the
$X$ side.  Then resume the presentation on $Y$.  Continue until every
pebble of $P_i$ has been collected on $X$.  Record the terminal
presentation state reached by this forward scan.  In the worst case,
the scan reaches the end of the recorded presentation sequence.

\medskip
\noindent
\textbf{(B2) Return $Y$ to the beginning of $S_i$.}
Keep the $X$ side fixed and execute backward exactly the portion of the
recorded presentation on $Y$ traversed in (B1), stopping as soon as the
state in (I2) is recovered.  Thus the reversal stops at the boundary
immediately preceding $S_i$ and never passes through a presentation
position belonging to $S_1,\ldots,S_{i-1}$.  By
Theorem~\ref{thm:pebble-presentation}, the forward scan in (B1) together
with this partial reversal costs $\Ord(N^2)$ moves.

\medskip
\noindent
\textbf{(B3) Order $P_i$ on $X$.}
Pause the presentation on $Y$.  All pebbles of $P_i$ are now stored on
the $X$ side.  Reconfigure the $f+\Ord(1)$ vertices available there so
that the pebbles of $P_i$ can be supplied at the interface $u$ in
exactly the order prescribed by the presentation positions of $S_i$.
By the general cubic bound on a graph of order $f+\Ord(1)$, this
requires
\[
    \Ord(f^3)
\]
moves.

\medskip
\noindent
\textbf{(B4) Write $S_i$.}
Starting again from the boundary preceding $S_i$, apply the
transcription mechanism of the preceding subsection while advancing
the recorded presentation on $Y$ through exactly the positions of
$S_i$.  At each presentation position, exchange the currently
presented pebble at $u$ with the next prescribed pebble supplied from
$X$.  After the last position of $S_i$ has been processed, stop the
presentation on $Y$ immediately after $S_i$.

No backward execution carrying the newly substituted pebbles to their
corresponding vertices of $Y$ is performed at this stage.  Consequently
(I1) holds with $i$ in place of $i-1$, while the state at which (B4)
stops is exactly the state required by (I2) for $S_{i+1}$ when $i<m$.
Thus the two invariants are restored for the next iteration.

After (B1)--(B4) have been performed for every block, the recorded
presentation on $Y$ has reached its final state, and every presentation
position of $Y$ has been assigned the pebble that must occupy its
corresponding target vertex.  Only now execute the entire recorded
presentation sequence on $Y$ backward once, from its final presentation
state to its initial presentation state.  During this final reversal no
new exchange with $X$ is made.  By the replacement principle from the
preceding subsection, this single backward execution carries every
pebble substituted at a presentation position to the corresponding
target vertex of $Y$.  Hence the target configuration on $Y$ is
obtained.  A final reconfiguration on the $X$ side completes the target
configuration on $G$.

There are $m=\Ord(N/f)$ blocks.
For each block, operations (B1), (B2), and (B4) together contribute at
most $\Ord(N^2)$ moves.  Hence, over all $m=\Ord(N/f)$ blocks, the
presentation and transcription operations contribute
\[
    \Ord\!\left(N^2\frac{N}{f}\right)
    =\Ord\!\left(\frac{N^3}{f}\right)
    =\ord(N^3),
\]
because $f=\omega(1)$.  The single final backward execution of the
complete recorded presentation sequence on $Y$ costs only
$\Ord(N^2)$ moves and is absorbed in this bound.

Operation (B3) costs $\Ord(f^3)$ for each block.  Hence all
$X$-side reconfigurations used in (B3) contribute
\[
    \Ord\!\left(f^3\frac{N}{f}\right)
    =\Ord(Nf^2)
    =\ord(N^3),
\]
because $f=\ord(N)$.
The final reconfiguration on $X$ costs
$\Ord(f^3)=\ord(N^3)$.  Hence the complete solution has length
\[
    \Ord\!\left(\frac{N^3}{f}+Nf^2\right)
    =\ord(N^3).
\]
\end{Proofof}

The proof of Theorem~\ref{thm:intermediate-cut-component} also shows
that the assumption $|V(X)|=\ord(N)$ is used only to make the general
cubic bound for the reconfigurations on the $X$ side subcubic.  If a
subcubic reconfiguration bound on $X$ is available independently, the
same batched-transcription argument applies even when $X$ has linear
order.  We record this consequence for later use.

\begin{cor}
\label{cor:linear-cut-subcubic-side}
Let $G$ be an $N$-vertex board graph of a feasible puzzle.  Suppose that
$G=X\cup Y$, where $X$ and $Y$ are connected,
\[
    V(X)\cap V(Y)=\{u\},
\]
and $u$ is a cut vertex.  Suppose further that
\[
    |V(X)|=\Theta(N),
\]
and that the reconfiguration cost of $X$ is $\ord(N^3)$.
%and that the worst-case shortest solution length of the PMP on $X$ is
%$\ord(N^3)$.  
Then the worst-case shortest solution length of the PMP
on $G$ is also
\[
    \ord(N^3).
\]
\end{cor}

\begin{Proofof}{Corollary\tume\ref{cor:linear-cut-subcubic-side}}
Apply the same batched-transcription construction as in the proof of
Theorem~\ref{thm:intermediate-cut-component}.  Since
$|V(X)|=\Theta(N)$, the presentation positions of $Y$ can be partitioned
into only $\Ord(1)$ blocks, each of size at most $\Theta(|V(X)|)$.
Consequently, all presentation and transcription operations on the
$Y$ side together require only $\Ord(N^2)$ moves.

For each block, the required reconfiguration on the $X$ side can, by
hypothesis, be performed in $\ord(N^3)$ moves.  Since only $\Ord(1)$
blocks are processed, the total cost of these reconfigurations remains
$\ord(N^3)$.  The final reconfiguration on $X$ is likewise
$\ord(N^3)$.  Therefore the complete solution sequence on $G$ has
length $\ord(N^3)$.
\end{Proofof}

%%%%%%%%%%%%%%%%%%%%%%%%%%%%%%%%%%%%%%%%%%%%%%%%%%%%%%%%%%%%%%%%%%%%%%%%%%%%%%%%

\subsection{Adding ears}

Let
\[
    E=u\,v\cdots r\,w
\]
be an ear whose endpoints $u,w$ belong to a connected graph $H$ and
whose internal vertices are new.  During the constructions below, a
vertex of the ear sometimes only has to be occupied, while the identity
of the occupying pebble is irrelevant.  A pebble used only for this
purpose will be called a \textit{filler pebble}.

%%%%%%%%%%%%%%%%%%%%%%%%%%%%%%%%%%%%%%%%%%%%%%%%%%%%%%%%%%%%%%%%%%%%%%%%%%%%%%%%

\subsubsection{A linear-size core}

\begin{thm}
\label{thm:linear-ear}
Let $H$ be an $N$-vertex connected board graph for which the worst-case
shortest solution length of the PMP is $\ord(N^3)$.  Let $E$ be an ear
of order $\Ord(N)$ whose endpoints belong to $H$, and let
\[
    G=H\cup E.
\]
Arbitrary new pebbles may be placed on the internal vertices of $E$.
Then the worst-case shortest solution length of the PMP on $G$ is also
$\ord(N^3)$.
\end{thm}

\begin{Proofof}{Theorem\tume\ref{thm:linear-ear}}
Put
\[
    \beta:=\left\lceil \frac{2|V(E)|}{N}\right\rceil .
\]
Since $|V(E)|=\Ord(N)$, we have $\beta=\Ord(1)$.  Read the target
pebble sequence on $E$ from $u$ toward $w$, and partition it into
consecutive blocks
\[
    S_1,S_2,\ldots,S_\beta,
\]
each containing at most $N/2$ pebbles.  Let
\[
    m:=\left\lceil\frac{\beta}{2}\right\rceil .
\]
We construct these blocks from the middle outward.

We begin with $S_m$.  The first task is to collect in $H$ all
pebbles that belong to $S_m$.  Reconfigure $H$ so that a batch of
pebbles not belonging to $S_m$, of size at most $N/2$, occurs first
in the presentation order at $u$.  Run the corresponding initial part
of the Pebble Presentation Algorithm and send these pebbles through
$uv$ into $E$.  As these filler pebbles enter $E$, the same number
of pebbles leave $E$ through its opposite end $w$ and enter $H$.
Any pebble belonging to $S_m$ that enters $H$ is retained there; the
remaining pebbles may subsequently be used as fillers.

The choice of the block size is important here.  Since
\[
    |S_m|\leq N/2
\]
and $H$ has $N$ vertices, $H$ can retain all pebbles of $S_m$
while still using its remaining capacity to supply the filler pebbles
sent into the ear.  Repeating this operation moves successive batches
of order $N$ through $E$.  Since
\[
    |V(E)|\leq \frac{\beta N}{2}
\]
and $\beta=\Ord(1)$, after only $\Ord(\beta)=\Ord(1)$ repetitions
all pebbles of $S_m$ have entered $H$.

Now reconfigure $H$ so that the pebbles of $S_m$ occur in their
target order in the presentation order at $u$, and insert them through
$uv$ into $E$.  Thus $S_m$ is formed as one consecutive ordered
block at the $u$-side of the ear.

We next prepare $S_{m+1}$, when it exists.  Continue to send filler
pebbles from $H$ into $E$ through $u$.  Their purpose is to expel
through $w$ the pebbles lying farther along the ear.  Whenever a pebble
needed for $S_{m+1}$ leaves $E$ through $w$, retain it in $H$.
Thus, during this sweep, $H$ discards filler pebbles into the ear and
replaces them by pebbles required for $S_{m+1}$.  Because
\[
    |S_{m+1}|\leq N/2,
\]
all pebbles of $S_{m+1}$ can be retained simultaneously in $H$.

The same sweep has a second effect.  As filler pebbles are inserted from
$u$, the already completed block $S_m$ is carried along the ear
toward $w$.  By the time $S_m$ reaches the $w$-end, every pebble
lying beyond it has been expelled through $w$, and hence all pebbles
required for $S_{m+1}$ have been collected in $H$.  Reconfigure
$H$ so that these pebbles occur in the target order of $S_{m+1}$,
and insert $S_{m+1}$ from the $w$-side.  The two completed blocks
now occur consecutively in the correct relative order
\[
    S_m,S_{m+1}.
\]

We then perform the same operation in the opposite direction.  Send
filler pebbles from the $w$-side into $E$, retain in $H$ the
pebbles expelled through $u$ that are required for $S_{m-1}$, and
continue until the already completed sequence $S_m,S_{m+1}$ has been
carried to the $u$-end of the ear.  At that moment all pebbles required
for $S_{m-1}$ have been collected in $H$.  Reconfigure $H$ to
place them in their target order and insert $S_{m-1}$ from the
$u$-side.  The completed sequence is then
\[
    S_{m-1},S_m,S_{m+1}.
\]

Continue in the same way, alternately from the two ends of the ear.  Thus
the blocks are added in the order
\[
    S_m,\ S_{m+1},\ S_{m-1},\ S_{m+2},\ S_{m-2},\ldots ,
\]
whenever the indicated block exists.  At each stage, filler pebbles are
sent into the ear only in order to expel, through the opposite end, the
pebbles needed for the next block.  The same operation simultaneously
transports the already completed sequence to that opposite end, which
is precisely where the next block must be attached.  Consequently the
completed blocks always form one consecutive segment of $E$ in their
correct relative order.  Eventually the whole target sequence on $E$
is obtained.

It remains to estimate the number of moves.  All reconfigurations
performed inside $H$ use only edges of $H$, and therefore remain
valid after the ear is adjoined.  Since $\beta=\Ord(1)$, only a
constant number of sweeps and block insertions are required.  Hence the
assumed reconfiguration procedure on $H$ is invoked only $\Ord(1)$
times, and its total contribution is
\[
    \Ord(1)\,\ord(N^3)=\ord(N^3).
\]
Likewise, only $\Ord(1)$ executions, or initial portions of executions,
of the Pebble Presentation Algorithm are required.  By
Theorem\tume\ref{thm:pebble-presentation}, their total contribution is
\[
    \Ord(N^2)=\ord(N^3).
\]
After the target sequence on $E$ has been completed, one final
reconfiguration of $H$ places the remaining pebbles in their target
positions, at an additional cost of $\ord(N^3)$.  Therefore the
complete solution on $G$ has length
\[
    \ord(N^3).
\]
\end{Proofof}

%%%%%%%%%%%%%%%%%%%%%%%%%%%%%%%%%%%%%%%%%%%%%%%%%%%%%%%%%%%%%%%%%%%%%%%%%%%%%%%%

\subsubsection{A sublinear core}

\begin{thm}
\label{thm:sublinear-core-linear-ear}
Let $G$ be an $N$-vertex board graph obtained by adjoining
an ear $E$ to a connected graph $H$, where $H$ is the board graph of a
feasible puzzle. Suppose that
\[
|V(H)|=f(N), \qquad f(N)=\omega(1)\cap o(N),
\]
and that $|V(E)|=\Ord(N)$.
Arbitrary new pebbles may be placed on the internal vertices of $E$.
Then the worst-case shortest solution length of the PMP on $G$ is
$\ord(N^3)$.
\end{thm}

\begin{Proofof}{Theorem\tume\ref{thm:sublinear-core-linear-ear}}
Write
\[
    f:=f(N)=|V(H)|.
\]
Use exactly the same middle-outward construction as in
Theorem\tume\ref{thm:linear-ear}, but now partition the target sequence
on $E$ into consecutive blocks
\[
    S_1,S_2,\ldots,S_\beta
\]
of size at most $f/2$, except possibly for a harmless adjustment at
one end.  Since $|V(E)|=\Ord(N)$,
\[
    \beta=\Ord\!\left(\frac{N}{f}\right).
\]

The mechanism is the same as in the preceding proof.  While one block
is being prepared, $H$ retains the pebbles already collected for that
block and uses its remaining capacity to supply filler pebbles.  The
reconfiguration algorithm on $H$ selects and arranges the filler
pebbles to be presented at the current endpoint of the ear.  Sending
them into $E$ expels pebbles from the opposite endpoint; those needed
for the block currently being prepared are retained in $H$, while the
others may later be used as fillers.  Once all pebbles of that block
have been collected, $H$ is reconfigured so that they are presented
in their target order and the block is inserted.  The same sweep also
carries the already completed block sequence to the opposite end of the
ear, so the next adjacent block can be attached there.  The direction
is then reversed, exactly as in the proof of
Theorem\tume\ref{thm:linear-ear}.

The difference is quantitative.  One reconfiguration of $H$ can
prepare only $\Theta(f)$ filler pebbles at a time.  Hence one such
batch advances the relevant portion of the ear by $\Theta(f)$
positions.  In the worst case, collecting all pebbles required for one
prescribed block therefore needs
\[
    \Ord\!\left(\frac{N}{f}\right)
\]
such batches, and hence the same number of reconfigurations of $H$.
Since there are
\[
    \Ord\!\left(\frac{N}{f}\right)
\]
blocks, the total number of invocations of the reconfiguration algorithm
on $H$ is
\[
    \Ord\!\left(\left(\frac{N}{f}\right)^2\right).
\]
Each invocation concerns a graph of order $f$, and the general cubic
bound therefore gives $\Ord(f^3)$ moves.  Thus all reconfigurations of
$H$ together cost
\[
    \Ord\!\left(
        \left(\frac{N}{f}\right)^2 f^3
    \right)
    =
    \Ord(N^2f)
    =
    \ord(N^3),
\]
because $f=\ord(N)$.

For the presentation cost, group all presentation segments used while
constructing one prescribed block into one block-presentation
procedure.  As in the preceding theorem, this procedure only sends
successive filler batches through one endpoint and finally inserts the
ordered block; by Theorem\tume\ref{thm:pebble-presentation}, its total
cost is at most
\[
    \Ord(N^2).
\]
There are $\Ord(N/f)$ blocks, so the total presentation cost is
\[
    \Ord\!\left(\frac{N}{f}N^2\right)
    =
    \Ord\!\left(\frac{N^3}{f}\right)
    =
    \ord(N^3),
\]
because $f=\omega(1)$.

After the target sequence on $E$ has been completed, one final
reconfiguration of $H$ costs
\[
    \Ord(f^3)=\ord(N^3).
\]
Combining the estimates gives
\[
    \Ord\!\left(\frac{N^3}{f}+N^2f\right)
    =
    \ord(N^3).
\]
\end{Proofof}

%%%%%%%%%%%%%%%%%%%%%%%%%%%%%%%%%%%%%%%%%%%%%%%%%%%%%%%%%%%%%%%%%%%%%%%%%%%%%%%%

\subsection{The scale of the components in the decomposition}

\begin{thm}
\label{thm:no-intermediate-decomposition-component}
Let $\{(G_N,P_N)\}_{N\to\infty}$ be a family of feasible
puzzles whose worst-case shortest solution length is $\Theta(N^3)$.
Then there exist constants $K<\infty$, $\delta>0$, and $N_0$ such
that, for every $N\geq N_0$ and every component
$X\in\mathcal{B}(G_N)$,
\[
    |V(X)|\leq K
    \qquad\text{or}\qquad
    |V(X)|\geq \delta N.
\]
Equivalently, there is no increasing sequence
$N_j\to\infty$ and no choice of components
$X_{N_j}\in\mathcal{B}(G_{N_j})$ such that
\[
    |V(X_{N_j})|\to\infty
    \qquad\text{and}\qquad
    \frac{|V(X_{N_j})|}{N_j}\to 0.
\]
%Equivalently, there is no choice of components
%$X_N\in\mathcal{B}(G_N)$ satisfying
%\[
%    |V(X_N)|=\omega(1)\cap\ord(N).
%\]
\end{thm}

\begin{Proofof}{Theorem\tume\ref{thm:no-intermediate-decomposition-component}}
Suppose that the stated uniform gap does not exist.  Then, recursively,
we can choose an increasing sequence $N_j$ and components
$X_{N_j}\in\mathcal{B}(G_{N_j})$ such that
\[
    j<|V(X_{N_j})|<\frac{N_j}{j}.
\]
Passing to this subsequence and relabeling it by $N$, we may therefore
assume that a component $X\in\mathcal{B}(G)$ satisfies
\[
    f(N):=|V(X)|=\omega(1)\cap\ord(N).
\]
Let $u\in V(X)$ be a cut vertex through which other components are
attached.  Let $H_1(u)$ be the connected subgraph consisting of $u$ and
all components of $G-u$ lying on the side away from the interior of
$X$, and let $H_2(u)$ be the complementary connected subgraph containing
$X$.  Thus
\[
    V(H_1(u))\cap V(H_2(u))=\{u\},
    \qquad
    V(H_1(u))\cup V(H_2(u))=V(G).
\]

Since $X\subseteq H_2(u)$ and $|X|=\omega(1)$, $H_2(u)$ is not bounded.
If $|V(H_2(u))|=\ord(N)$, then
Theorem~\ref{thm:intermediate-cut-component} gives a subcubic solution,
contradicting the hypothesis.  Hence $|V(H_2(u))|=\Theta(N)$.  The same
theorem shows that $H_1(u)$ has either order $\Ord(1)$ or order
$\Theta(N)$.

Consider first those vertices $u$ for which
$|V(H_1(u))|=\Theta(N)$.  The interiors of distinct such subgraphs are
disjoint, so only $\Ord(1)$ of them exist.

Choose a set $B^\ast$ of blank vertices in the part that will remain
unfinished, with
\[
    |B^\ast|=N-n,
\]
and put
\[
    S^\ast:=V(G)\setminus B^\ast.
\]
We choose $B^\ast$ so that the interiors
$H_1(u)-\{u\}$ of all the above linear-order attached subgraphs are
contained in $S^\ast$.  Such a choice is possible because
$|V(X)|=f(N)=\omega(1)$ whereas the number of blanks is $\Ord(1)$.

Let $T$ be an arbitrary spanning tree of $G$.  Applying
Lemma\tume\ref{sbs0} to $T$, with the support of $g$ as $S_1$ and
$S^\ast$ as $S_2$, gives a move sequence of length
\[
    \Ord(n\D(T))=\Ord(N^2)
\]
that transforms $g$ into a configuration $g^\ast$ whose support is
exactly $S^\ast$.  We record this sequence and work below toward
$g^\ast$.
%We first replace the prescribed target configuration $g$ by a packed
%target configuration $g^\ast$.  In $g^\ast$, the interiors
%$H_1(u)-\{u\}$ of all the above linear-order attached subgraphs are
%filled completely with pebbles, while all blank vertices are placed in
%the part that will remain unfinished.  This is possible because
%$|V(X)|=f(N)=\omega(1)$ whereas the total number of blanks is
%$\Ord(1)$.
%
%The configuration $g^\ast$ may be obtained from $g$ by the same
%elementary packing operation used earlier in the paper: using a
%spanning tree, move the pebbles one by one into the prescribed packed
%positions.  This requires at most
%\[
%    \Ord(Nn)=\Ord(N^2)
%\]
%moves.  We record this move sequence and work below toward $g^\ast$.

Now process the linear-order subgraphs $H_1(u)$ one at a time.
After the interior of one such subgraph has been completed, delete that
interior from the current unfinished graph; its attachment vertex in
$X$ is retained.  Thus the unfinished graph remains connected, and all
subsequent presentation procedures and exchanges are performed entirely
within this unfinished graph.

Before the current $H_1(u)$ is completed, use presentation-based
exchanges through $u$, restricted to the current unfinished graph, to
bring to the $H_1(u)$-side every pebble assigned there by $g^\ast$.
%Now process the linear-order subgraphs $H_1(u)$ one at a time.  Before
%the current $H_1(u)$ is completed, use presentation-based exchanges
%through $u$ to bring to the $H_1(u)$-side every pebble assigned there
%by $g^\ast$.  
At the same time, any blank occurring in
$H_1(u)-\{u\}$ is simply pushed to the unfinished side.  No attempt is
made to place the blanks at their final positions individually.  This
preliminary separation requires $\Ord(N^2)$ moves.

After the separation, every pebble required for the packed target on
$H_1(u)$ is available on that side, and no such pebble remains hidden
in another component attached to $X$.  
We may therefore apply the batched-transcription construction from the
proof of Theorem~\ref{thm:intermediate-cut-component}, using $X$ as the
working region through $u$.
Since
\[
    |V(X)|=f(N),
\]
the packed target configuration on the current $H_1(u)$ is completed in
\[
    \Ord\!\left(\frac{N^3}{f(N)}+Nf(N)^2\right)
    =\ord(N^3)
\]
moves.

Once $H_1(u)-\{u\}$ has been completed, its interior is never modified
again.  Since there are only $\Ord(1)$ linear-order subgraphs
$H_1(u)$, all of them can be completed in $\ord(N^3)$ moves.
%Consider first those vertices $u$ for which
%$|V(H_1(u))|=\Theta(N)$.  The interiors of distinct such subgraphs are
%disjoint, so only $\Ord(1)$ of them exist.  For each of them, apply the
%block construction from the proof of
%Theorem~\ref{thm:intermediate-cut-component}, using the intermediate-size
%component $X$ to hold the blocks that are exchanged through $u$.  Since
%$|V(X)|=f(N)$, the target configuration on each such $H_1(u)$ is
%completed in
%\[
%    \Ord\!\left(\frac{N^3}{f(N)}+Nf(N)^2\right)
%    =\ord(N^3)
%\]
%moves.  Because there are only $\Ord(1)$ such subgraphs, all of them can
%be completed in $\ord(N^3)$ moves.

For every remaining cut vertex $u$ of $X$ we now have
$|V(H_1(u))|=\Ord(1)$.  Let $C(X)$ be the set of these vertices and put
\[
    X':=X\cup\bigcup_{u\in C(X)}H_1(u).
\]
Since $|C(X)|\le |V(X)|$ and every added $H_1(u)$ has bounded order,
\[
    |V(X')|=\Ord(|V(X)|)=\ord(N).
\]
All parts of $G$ outside $X'$ have already been completed according to
the packed target configuration $g^\ast$.  The remaining configuration
on $X'$ can therefore be solved at once by the general cubic
reconfiguration bound in
\[
    \Ord(|V(X')|^3)=\ord(N^3)
\]
moves.  We have thus obtained the packed target configuration
$g^\ast$ on the whole graph.

Finally, reverse the move sequence supplied by
Lemma\tume\ref{sbs0} that transformed $g$ into $g^\ast$.
This recovers the prescribed target configuration $g$ in
\[
    \Ord(N^2)=\ord(N^3)
\]
moves.
%Finally, reverse the previously recorded packing sequence from $g$ to
%$g^\ast$.  This transforms $g^\ast$ into the prescribed target
%configuration $g$ in
%\[
%    \Ord(Nn)=\Ord(N^2)=\ord(N^3)
%\]
%moves.

Together with the work already spent on the $\Theta(N)$ subgraphs
$H_1(u)$ and on $X'$, this gives an $\ord(N^3)$ solution of the
original instance, a contradiction.
%All parts of $G$ outside $X'$ have already been completed.  The remaining
%configuration on $X'$ can therefore be solved at once by the general
%cubic reconfiguration bound in
%\[
%    \Ord(|V(X')|^3)=\ord(N^3)
%\]
%moves.  Together with the work already spent on the $\Theta(N)$
%subgraphs $H_1(u)$, this gives an $\ord(N^3)$ solution of the original
%instance, a contradiction.
\end{Proofof}

%%%%%%%%%%%%%%%%%%%%%%%%%%%%%%%%%%%%%%%%%%%%%%%%%%%%%%%%%%%%%%%%%%%%%%%%%%%%%%%%

\subsection{Linear-size components}

Theorem~\ref{thm:no-intermediate-decomposition-component} shows that every
component of $\mathcal{B}(G)$ has order either $\Ord(1)$ or $\Theta(N)$.
We now show that a component of the latter type must itself carry cubic
reconfiguration complexity.

Linear-order subgraphs attached through a single cut vertex can be handled
by Corollary~\ref{cor:linear-cut-subcubic-side}.  The remaining difficulty
is that a linear-size component may have many bounded subgraphs attached at
different cut vertices.  To handle all such attachments without
accumulating a cubic overhead, we use the core itself as a temporary
exchange buffer.

\begin{la}[Core-buffer redistribution]
\label{lem:core-buffer-redistribution}
Let $H$ be a core containing $b=\Ord(1)$ blanks, and let
$C_1,\ldots,C_r$ be pairwise disjoint full containers attached to $H$
at pairwise distinct interface vertices $u_1,\ldots,u_r\in V(H)$, with
\[
    r\le |V(H)|,
    \qquad
    |C_i|\le M
\]
for every $i$, where $M$ is fixed.  Suppose that a local operation can
exchange a prescribed pebble of $C_i$ with a prescribed pebble supplied
from $H$ through $u_i$, while restoring every vertex of $H$ outside
$u_i$ and every other container, and that the pebbles on $H$ can be
reconfigured arbitrarily at a cost charged separately.

Then the target set of pebbles can be placed in every $C_i$ using
$\Ord_M(1)$ reconfigurations of $H$ and $\Ord_M(|V(H)|)$ local exchanges.
\end{la}

\begin{Proofof}{Lemma\tume\ref{lem:core-buffer-redistribution}}
Put
\[
    h:=|V(H)|,
    \qquad
    q:=h-b,
\]
so that $q$ pebbles occupy the core.  For each $i$, let $P_i$ be the set
of pebbles that must finally occupy $C_i$, and let $P_H$ be the set of
pebbles that must finally occupy $H$.  We ignore the exact positions
inside the containers and consider only the container to which each pebble
belongs.

We use the local exchanges in rounds.  In one exchange cycle, at most one
exchange is made at each interface.  Since the interfaces are distinct,
all pebbles supplied from $H$ in such a cycle can be positioned
simultaneously by a single reconfiguration of $H$.  Hence, because
$|C_i|\le M$, any prescribed collection of pebbles from the containers, or
any prescribed distribution of pebbles from $H$ to their target
containers, can be carried out in at most $M$ exchange cycles.

A pebble once inserted into its correct container is never removed.  A
pebble of $P_H$ may be used temporarily as a filler when a misplaced
attachment-target pebble is collected from a container.  Every $C_i$ is
full and $|P_i|=|C_i|$, so when $C_i$ has received its complete target set,
no filler can remain there; it has necessarily been returned to $H$ by a
later exchange.  Thus no separate filler-recovery step is required.

We now bound the number of rounds.  Let $m$ be the number of pebbles
belonging to some $P_i$ that are not currently in their correct container,
and let $k$ of them currently lie on $H$.  Suppose first that $m\ge q$.
Return the $k$ attachment-target pebbles currently on $H$ to their correct
containers.  Since each container receives at most $M$ pebbles, this takes
at most $M$ exchange cycles and gives $k$ updates.  Some of these exchanges
may simultaneously return to $H$ a misplaced pebble destined for another
container; let $t$ be the number of these additional updates.

The $q-k$ pebbles of $P_H$ that were on $H$ at the beginning of the round
have not been used in this distribution step and remain available as
fillers.  If $t<q-k$, then after the distribution step at least
\[
    m-k-t\ge q-k-t
\]
misplaced attachment-target pebbles still lie in wrong containers.  Use
$q-k-t$ of the available fillers to collect that many of them onto $H$.
Again, because each container has order at most $M$, this requires at most
$M$ exchange cycles.  Thus one round with $m\ge q$ makes at least
\[
    k+t+(q-k-t)=q
\]
updates, using at most $2M$ reconfigurations of $H$.

Here an update means either moving a misplaced attachment-target pebble
from a wrong container to $H$, or moving such a pebble from $H$ into its
correct container.  Each attachment-target pebble can undergo each of
these two updates at most once.  Initially
\[
    m\le \sum_{i=1}^r |C_i|\le Mr\le Mh,
\]
so at most $2Mh$ updates occur in total.  Hence the number of rounds while
$m\ge q$ is at most
\[
    \frac{2Mh}{q}
    =\frac{2Mh}{h-b}
    =\Ord_M(1).
\]
For all sufficiently large $h$ we have $h\ge 2b$, so this is at most
$4M$.

Once $m<q$, the $m-k$ misplaced attachment-target pebbles still lying in
wrong containers satisfy $m-k<q-k$, so the core has enough fillers to
collect all of them.  At most $M$ further exchange cycles bring them all to
$H$, and at most another $M$ exchange cycles place all remaining
attachment-target pebbles into their correct containers.  Thus the total
number of reconfigurations of $H$ is $\Ord_M(1)$, and the total number of
local exchanges is $\Ord_M(h)$.  When the process stops, every $C_i$
contains exactly $P_i$, and consequently the pebbles remaining on $H$ are
exactly $P_H$.
\end{Proofof}

\begin{thm}
\label{thm:linear-component-cubic}
Let $\{(G_N,P_N)\}_{N\to\infty}$ be a family of feasible
puzzles whose worst-case shortest solution length is $\Theta(N^3)$.
Let $X_N\in\mathcal{B}(G_N)$ be a choice of components satisfying
\[
    |V(X_N)|=\Theta(N).
\]
Then the reconfiguration cost of $X_N$ is
\[
    \Theta(|V(X_N)|^3)=\Theta(N^3).
\]
\end{thm}

\begin{Proofof}{Theorem\tume\ref{thm:linear-component-cubic}}
Suppress the subscript $N$ according to the asymptotic convention above.
Suppose, to the contrary, that the reconfiguration cost of $X$ is
$\ord(N^3)$.  We show that then every pair of configurations of the
puzzle $(G,P)$ can be connected in $\ord(N^3)$ moves, contradicting the
hypothesis.  Fix arbitrary initial and target configurations $f$ and $g$
on $G$.  Let
\[
    b:=N-n.
\]
Since $b=\Ord(1)$ and $|V(X)|=\Theta(N)$, choose a set
\[
    B^\ast\subseteq V(X),
    \qquad
    |B^\ast|=b,
\]
and put
\[
    S^\ast:=V(G)\setminus B^\ast.
\]
Thus every vertex outside $X$ belongs to $S^\ast$.

Let $T$ be an arbitrary spanning tree of $G$.  Apply
Lemma\tume\ref{sbs0} on $T$ first to the support of $f$ and then to the
support of $g$, in both cases using $S^\ast$ as the target support.
We thereby obtain configurations $f^\ast$ and $g^\ast$ such that
\[
    \operatorname{sup}(f^\ast)
    =
    \operatorname{sup}(g^\ast)
    =
    S^\ast.
\]
Each transformation requires
\[
    \Ord(n\D(T))=\Ord(N^2)
\]
moves.  We record both move sequences.  It is therefore sufficient to
transform $f^\ast$ into $g^\ast$ in $\ord(N^3)$ moves: prepending the
recorded sequence from $f$ to $f^\ast$ and appending the reverse of the
recorded sequence from $g$ to $g^\ast$ changes the total length by only
$\Ord(N^2)$.
%Fix arbitrary feasible initial and target configurations $f$ and $g$ on $G$.
%As in the proof of Theorem~\ref{thm:no-intermediate-decomposition-component},
%we first replace them by packed configurations $f^*$ and $g^*$ in which all
%blank vertices lie in $X$.  This is possible because the number of blanks is
%$\Ord(1)$ whereas $|V(X)|=\Theta(N)$.  Using a spanning tree, each packing
%can be performed in $\Ord(N^2)$ moves.  We shall transform $f^*$ into $g^*$;
%at the end we reverse the packing sequence from $g$ to $g^*$.

For every cut vertex $u$ of $X$ through which other components are attached,
define $H_1(u)$ and $H_2(u)$ as in the proof of
Theorem~\ref{thm:no-intermediate-decomposition-component}: $H_1(u)$ consists
of $u$ together with the whole side attached through $u$ away from the
interior of $X$, while $H_2(u)$ is the complementary connected subgraph
containing $X$.  
Since $H_2(u)$ contains $X$, we have
$|V(H_2(u))|=\Theta(N)$.
By Theorem~\ref{thm:intermediate-cut-component},
$H_1(u)$ cannot have order $\omega(1)\cap\ord(N)$.
Hence, after passing to a subsequence if necessary,
each $H_1(u)$ has order either $\Ord(1)$ or $\Theta(N)$.
%By Theorem~\ref{thm:no-intermediate-decomposition-component},
%each $H_1(u)$ has order either $\Ord(1)$ or $\Theta(N)$.

We first deal with the subgraphs $H_1(u)$ of order $\Theta(N)$.  Their
interiors are pairwise disjoint, so there are only $\Ord(1)$ of them.
Process them one at a time.
For the current $H_1(u)$, first use the presentation-based exchange
through $u$ described in Subsection\tume\ref{subsec:transcription-cut} 
to bring to the $H_1(u)$-side exactly the pebbles assigned there by 
the packed target $g^\ast$.
Because substitutions made at different presentation positions can be
realized together by reversing the recorded presentation sequence, this
separation requires only $\Ord(N^2)$ moves.
After this separation, apply the batched-transcription argument of
Corollary~\ref{cor:linear-cut-subcubic-side}, with $X$ as the
linear-size working region through $u$.  Indeed,
\[
    |V(X)|=\Theta(N),
\]
and, under the present contradictory assumption, the worst-case shortest
solution length of the PMP on $X$ is $\ord(N^3)$.  Since also
$|V(H_1(u))|=\Theta(N)$, the same argument completes the prescribed
configuration on $H_1(u)$ in $\ord(N^3)$ moves while leaving the
previously completed parts untouched.
Once $H_1(u)-\{u\}$ has been completed, fix its interior permanently and
remove it from the unfinished graph.  Since there are only $\Ord(1)$ such
linear-order subgraphs, all of them can be completed in
$\ord(N^3)$ moves altogether.

We are left with $X$ together with the attached subgraphs $H_1(u)$ of bounded
order.  Put
\[
    C_u:=V(H_1(u))\setminus\{u\}.
\]
The sets $C_u$ are pairwise disjoint, and there is a constant $M$ such that
\[
    |C_u|\le M
\]
for every remaining $u$.  Because both $f^*$ and $g^*$ are packed, every
vertex of every $C_u$ is occupied by a pebble in both configurations; all
blanks lie in $X$.

We now ignore the exact positions inside the bounded sets $C_u$ and record
only which container each pebble belongs to.  Apply
Lemma~\ref{lem:core-buffer-redistribution} with core $H=X$ and with the
containers $C_u$.  Its local exchange operation is provided by the bounded
presentation operation through the interface $u$: a prescribed pebble of
$C_u$ can be presented at $u$ and exchanged with a prescribed pebble supplied
from $X$, while $X$ outside $u$ and all other bounded attachments are restored
to their previous states.

The hypotheses of the lemma hold because there is at most one container
$C_u$ for each vertex $u\in V(X)$ and $|C_u|\le M$.  Since
$|V(X)|=\Theta(N)$ and the number of blanks is $b=\Ord(1)$, the lemma uses
only $\Ord(1)$ reconfigurations of $X$ and $\Ord(N)$ local exchanges.
The latter cost only $\Ord(N^2)$ moves altogether, including the motion of
the constant blank workspace between interfaces.  Under the contradictory
assumption, each reconfiguration of $X$ costs $\ord(N^3)$.  Hence, after
$\ord(N^3)$ moves, every bounded container $C_u$ contains exactly the set of
pebbles assigned to it by $g^*$.

Only the internal order inside the bounded attachments remains.  Visit their
interfaces once more with the constant blank workspace and complete the exact
target arrangement inside each $H_1(u)$.  Each local reconfiguration costs
$\Ord(1)$ moves and the traversal of $X$ costs $\Ord(N)$, so all bounded
attachments are completed in $\Ord(N)$ additional moves.  Finally reconfigure
$X$ itself to its packed target configuration in $\ord(N^3)$ moves.

We have therefore transformed $f^\ast$ into $g^\ast$ in
$\ord(N^3)$ moves.  Together with the transformation from $f$ to
$f^\ast$ and the reverse of the transformation from $g$ to $g^\ast$,
both supplied by Lemma\tume\ref{sbs0} and both of length $\Ord(N^2)$,
this yields a solution from $f$ to $g$ of length $\ord(N^3)$.
Thus every pair of configurations of the puzzle $(G,P)$ could be
connected by a subcubic solution sequence, contradicting the assumption
that the worst-case shortest solution length on $G$ is $\Theta(N^3)$.  

Therefore the reconfiguration cost of $X$ cannot be $\ord(N^3)$.
%Therefore the PMP on $X$ cannot have worst-case shortest solution
%length $\ord(N^3)$.  
Moreover, the construction above uses only a
constant number of reconfigurations of $X$, together with
$\Ord(N^2)$ additional moves.  
Hence, if $L_X(N)$ denotes the reconfiguration cost of $X$,
%Hence, if $L_X(N)$ denotes the
%worst-case shortest solution length of the PMP on $X$, 
then the worst-case shortest solution length on $G$ is
\[
    \Ord(L_X(N)+N^2).
\]
Since the latter is $\Theta(N^3)$, it follows that
\[
    L_X(N)=\Omega(N^3).
\]
Together with the general upper bound
\[
    L_X(N)=\Ord(|V(X)|^3)=\Ord(N^3),
\]
where $|V(X)|=\Theta(N)$, we obtain
\[
    L_X(N)=\Theta(N^3).
\]
%Therefore the PMP on $X$ cannot have
%worst-case shortest solution length $\ord(N^3)$.  Since the general
%upper bound is $\Ord(|V(X)|^3)=\Ord(N^3)$ and $|V(X)|=\Theta(N)$, its
%worst-case shortest solution length is $\Theta(N^3)$.
\end{Proofof}

\begin{cor}
\label{cor:bridge-tree-components-bounded}
Let $G$ be an $N$-vertex board graph for which the worst-case shortest
solution length of the PMP is $\Theta(N^3)$.  Then every bridge-tree component
of $\mathcal{B}(G)$ has order $\Ord(1)$.
\end{cor}

\begin{Proofof}{Corollary\tume\ref{cor:bridge-tree-components-bounded}}
By Theorem~\ref{thm:no-intermediate-decomposition-component}, every component
of $\mathcal{B}(G)$ has order either $\Ord(1)$ or $\Theta(N)$.  
Suppose that a
bridge-tree component $X$ has order $\Theta(N)$.  Since $X$ is a tree,
Theorem~\ref{main}, together with the convention of
Subsection~\ref{subsec:biconnected-bridge-tree-decomposition}, gives
\[
    \Ord\bigl(|V(X)|^2\log(1+|V(X)|)\bigr)
    =\ord(N^3)
\]
for the reconfiguration cost of $X$.
This contradicts Theorem~\ref{thm:linear-component-cubic}.
%Suppose that a
%bridge-tree component $X$ has order $\Theta(N)$.  Since $X$ is a tree,
%Theorem~\ref{main} gives, uniformly over all feasible instances on $X$, a
%solution sequence of length
%\[
%    \Ord\bigl(|V(X)|^2\log(1+|V(X)|)\bigr)
%    =\ord(N^3).
%\]
%This contradicts Theorem~\ref{thm:linear-component-cubic}, which requires
%every $\Theta(N)$-size component of $\mathcal{B}(G)$ to have cubic worst-case
%complexity. 
Hence no bridge-tree component has order $\Theta(N)$, and every
bridge-tree component has order $\Ord(1)$.
\end{Proofof}

%%%%%%%%%%%%%%%%%%%%%%%%%%%%%%%%%%%%%%%%%%%%%%%%%%%%%%%%%%%%%%%%%%%%%%%%%%%%%%%%

\subsection{Closure under arbitrary ear extensions}
\label{subsec:arbitrary-ear-extensions}

We now combine the preceding routing and decomposition mechanisms to obtain
the closure result for arbitrary successive ear extensions.

For the rest of this subsection, put
\[
    b:=N-n.
\]
Since throughout Section~\ref{sec:structural-restrictions} we have
$b=\Ord(1)$, all constants below may be chosen uniformly over the
finitely many possible values of $b$.  It is therefore enough, in the
proofs, to regard $b\geq 2$ as fixed.

For a connected board graph $F$ with $|V(F)|>b$, let $L(F)$ denote the
maximum, over all mutually reachable pairs of configurations of
$|V(F)|-b$ labeled pebbles on $F$, of the length of a shortest solution
sequence between the pair.  We call $F$ \textit{admissible} if $F$ is
connected, contains a cycle but is not itself a cycle, and has no
isthmus with at least $b$ vertices.  By Theorem~A, since $b\geq 2$,
this is equivalent to saying that $F$ contains a cycle and that the
puzzle with $|V(F)|-b$ labeled pebbles on $F$ is feasible.  In
particular, every two configurations of these pebbles on $F$ are
mutually reachable.

We shall use two elementary monotonicity observations repeatedly.  If
$F$ is admissible and edges are added between vertices of $F$, then the
resulting graph $F'$ is admissible and
\[
    L(F')\leq L(F).
\]
Indeed, adding edges creates neither a new bridge nor a new cut vertex,
and every move sequence on $F$ is also a move sequence on $F'$.  Also,
if $P$ is an $F$-path, that is, a path of length at least two whose
endpoints are distinct vertices of $F$ and whose internal vertices lie
outside $F$, then $F\cup P$ is admissible.  Every edge of $P$ lies on
the cycle obtained from $P$ and an $F$-path joining its endpoints, no
internal vertex of $P$ is a cut vertex of $F\cup P$, and any isthmus of
$F\cup P$ is therefore already an isthmus of $F$.

\begin{thm}
\label{thm:linear-ear-closure}
Let $G$ be an $N$-vertex biconnected board graph that is not a cycle,
and let $G_0\subseteq G$ be a connected subgraph with
\[
    f(N):=|V(G_0)|=\omega(1).
\]
Suppose that one of the following holds.
\begin{enumerate}
\item $f(N)=\Theta(N)$, the graph $G_0$ is admissible, and
\[
    L(G_0)=\ord(N^3).
\]
\item $f(N)=\omega(1)\cap\ord(N)$, and $G_0$ is either admissible or a
cycle.
\end{enumerate}
Then
\[
    L(G)=\ord(N^3)=\ord(|V(G)|^3).
\]
In particular, under the same hypotheses on $G_0$, if $G$ is a biconnected 
graph obtained from $G_0$ by adjoining an arbitrary sequence of ears, then 
the same conclusion holds; no restriction is imposed on the number or 
the individual lengths of the ears.
\end{thm}

The following elementary lemma is the only additional structural fact
needed for the closure proof.  It converts depth outside the current
core into a long ear that can be absorbed into the core.

\begin{la}[Deep-island ear]
\label{lem:core-remainder-reduction}
Let $G$ be biconnected, let $H\subseteq G$ be connected with
$|V(H)|\geq2$, and let $Z$ be a connected component of $G-V(H)$.
Let $u\in V(H)$ be adjacent to $Z$, let $v\in V(Z)$, and put
\[
    d:=d_{G[Z\cup\{u\}]}(u,v).
\]
Then there exists an $H$-path $P$ whose internal vertices lie in $Z$
and whose length is at least $d$.
\end{la}

\begin{Proofof}{Lemma\tume\ref{lem:core-remainder-reduction}}
By the fan lemma, there are two paths $F_1,F_2$ from $v$ to two
distinct vertices $a,b\in V(H)$ such that $F_1$ and $F_2$ meet only at
$v$ and all their internal vertices lie in $Z$.

If $a=u$ or $b=u$, then $F_1\cup F_2$ is an $H$-path and the branch
joining $u$ to $v$ has length at least
$d_{G[Z\cup\{u\}]}(u,v)=d$.  Hence the conclusion follows.

Assume now that $a,b\neq u$.  Let $Q$ be a shortest $u$--$v$ path in
$G[Z\cup\{u\}]$, and traverse $Q$ from $u$ until it first meets
$F_1\cup F_2$, say at $y$.  Suppose without loss of generality that
$y\in V(F_1)$.  Then
\[
    P:=Q[u,y]\cup F_1[y,v]\cup F_2[v,b]
\]
is an $H$-path from $u$ to $b$: by the choice of $y$, the first piece
meets the other two only at $y$, and the two fan paths meet only at
$v$.  Moreover,
\[
\begin{split}
    |E(P)|
    &\geq |E(Q[u,y])|+|E(F_1[y,v])|\\
    &\geq d_{G[Z\cup\{u\}]}(u,y)
           +d_{G[Z\cup\{u\}]}(y,v)\\
    &\geq d_{G[Z\cup\{u\}]}(u,v)=d.
\end{split}
\]
The case $y\in V(F_2)$ is symmetric.
\end{Proofof}

\begin{Proofof}{Theorem\tume\ref{thm:linear-ear-closure}}
We give one argument for both alternatives in the statement.

If $G_0$ is admissible, put $H_0:=G_0$.  If $G_0$ is a cycle, then
$G_0\neq G$ and there is an edge $cx$ with $c\in V(G_0)$ and
$x\notin V(G_0)$.  In this case put
\[
    H_0:=G_0+cx.
\]
Then $H_0$ is admissible: it is connected, contains a cycle but is not
a cycle, and its only possible new isthmus has one vertex.  Since
$b\geq2$, this is allowed.  In either case, writing
\[
    f_0:=|V(H_0)|,
\]
we have $f_0=\omega(1)$.  Moreover
\[
    L(H_0)=\ord(N^3).
\]
Indeed, this is assumed when $f_0=\Theta(N)$; when
$f_0=\ord(N)$ it follows from the classical cubic bound
$L(H_0)=\Ord(f_0^3)$.

We shall use the following quantitative form of
Theorem~\ref{thm:sublinear-core-linear-ear}.  Retaining the constants in
its proof, there are constants $a\geq1$ and $C$, depending only on $b$,
such that for every admissible $H$ and every $H$-path $P$, if
\[
    \rho:=\frac{|V(H\cup P)|}{|V(H)|},
\]
then
\begin{equation}
    L(H\cup P)
    \leq
    a\rho^2L(H)+C\rho |V(H\cup P)|^2.
    \label{eq:depth-ear-quantitative}
\end{equation}
This estimate also contains the linear-core situation of
Theorem~\ref{thm:linear-ear}; when $|V(H)|=\Theta(N)$, the ratio $\rho$
is bounded by a constant.

Fix $\varepsilon>0$.  Starting from $H_0$, construct an increasing
sequence of admissible subgraphs
\[
    H_0\subseteq H_1\subseteq\cdots\subseteq H_t
\]
as follows.  For the current $H_j$, call each connected component of
$G-V(H_j)$ an island.  For every island $Z$, choose one neighbor
$u_Z\in V(H_j)$ and call it the interface assigned to $Z$.  If for some
island $Z$ there is a vertex $v\in Z$ with
\[
    d_{G[Z\cup\{u_Z\}]}(u_Z,v)\geq\varepsilon N,
\]
Lemma~\ref{lem:core-remainder-reduction} gives an $H_j$-path $P_j$ of
length at least $\varepsilon N$ whose internal vertices lie in $Z$.
Put
\[
    H_{j+1}:=H_j\cup P_j.
\]
By the elementary observation preceding the theorem, $H_{j+1}$ is
admissible.  If no such island exists, stop.

Every absorption adds at least $\varepsilon N-1$ new vertices.  Hence,
for all sufficiently large $N$,
\[
    t=\Ord(1/\varepsilon).
\]
We next verify that the final core still has subcubic reconfiguration
cost.  If no absorption occurs, this is already true for $H_0$.
Suppose that an absorption occurs.  For the first one, if
$f_0=\Theta(N)$ then $\rho=\Ord(1)$ in
\eqref{eq:depth-ear-quantitative}, and hence
$L(H_1)=\ord(N^3)$.  If instead $f_0=\ord(N)$, then the classical bound
$L(H_0)=\Ord(f_0^3)$ and $\rho\leq N/f_0$ give
\[
\begin{split}
    L(H_1)
    &\leq
    \Ord\!\left(\frac{N^2}{f_0^2}f_0^3\right)
    +\Ord\!\left(\frac{N}{f_0}N^2\right)\\
    &=\Ord(N^2f_0)+\Ord\!\left(\frac{N^3}{f_0}\right)
    =\ord(N^3),
\end{split}
\]
because $f_0=\omega(1)\cap\ord(N)$.  After the first absorption we have
$|V(H_1)|\geq\varepsilon N-1$, so every later ratio satisfies
\[
    \rho\leq \frac{N}{\varepsilon N-1}=\Ord_\varepsilon(1).
\]
Since $t=\Ord(1/\varepsilon)$, repeated use of
\eqref{eq:depth-ear-quantitative} therefore gives, for a constant
$K(\varepsilon)$ independent of $N$,
\[
    L(H_t)\leq K(\varepsilon)\bigl(L(H_1)+N^2\bigr)=\ord(N^3).
\]

Put
\[
    H^*:=G[V(H_t)].
\]
Adding the missing edges on the same vertex set preserves
admissibility and cannot increase the reconfiguration cost, so
\[
    L(H^*)\leq L(H_t)=\ord(N^3).
\]
The islands of $G-V(H^*)$ are the same as those present at the last
step of the absorption construction.  Keep for each island $Z$ the
interface $u_Z$ assigned to it at that last step.  By the stopping rule,
for every such island $Z$ and every $v\in Z$,
\begin{equation}
    d_{G[Z\cup\{u\}]}(u,v)<\varepsilon N.
    \label{eq:shallow-islands}
\end{equation}

For an interface $u$, the union of the islands assigned to $u$ will be
called the bundle at $u$.  For each such island $Z$, choose a breadth-first
spanning tree of $G[Z\cup\{u\}]$ rooted at $u$, and let $T_B$ be the
union of these trees over all islands in the bundle $B$.  Thus $T_B$ is
a tree rooted at $u$, and by~\eqref{eq:shallow-islands} every vertex of
$B$ has depth less than $\varepsilon N$.  Consequently every path in
$T_B$ between two vertices of $B\cup\{u\}$ has length less than
$2\varepsilon N$.

We use two local facts.  First, if $B$ is any union of islands and
\[
    Y:=G[V(G)\setminus B],
\]
then $Y$ is admissible.  Indeed, $Y$ is connected and contains the
admissible graph $H^*$.  Every edge of $Y$ outside $H^*$ lies on a
cycle of $Y$: subdivide such an edge if necessary and apply the fan
lemma from the subdividing vertex to $H^*$, then join the two distinct
endpoints in $H^*$.  Similarly, if $H^*-v$ is connected, then every
component of every remaining island minus $v$ has a path to $H^*-v$ in
$G-v$, since $G$ is biconnected.  Hence every cut vertex, and therefore
every isthmus, of $Y$ already comes from $H^*$.  Thus $Y$ has no
isthmus with at least $b$ vertices.

Second, the one-path cycle-borrowing argument of
Section~\ref{Main2Proof}, refined in
Subsection~\ref{refinement-general}, gives the following.  Let
\[
    Q=(q_0,q_1,\ldots,q_d)
\]
be an occupied path in $T_B$.  The cyclic shift
\[
    q_0\to q_d,
    \qquad
    q_i\to q_{i-1}\quad(1\leq i\leq d)
\]
or its reverse can be performed while restoring every vertex outside
$Q$ in
\[
    \Ord(bdN+N)
\]
moves.  Indeed, at every obstructing step we may borrow a cycle of
$G$ through the relevant vertex, and because $G$ is biconnected its
length is at most $N$.  Since $d<2\varepsilon N$, every such push costs
\begin{equation}
    \Ord(b\varepsilon N^2+N).
    \label{eq:shallow-push}
\end{equation}
If the interface itself is temporarily blank, occupy it first by moving
the nearest pebble of $H^*$ to the interface and reverse these moves
afterward; since $b$ is fixed, this changes only the constant in the
bound.

We now complete arbitrary initial and target configurations.  Choose a
set $W\subseteq V(H^*)$ of $b$ vertices.  For later use, choose a
spanning tree $T_H$ of $H^*$.  For each island $Z$ of $G-V(H^*)$, choose
a spanning tree of $Z$ and add one edge joining it to its assigned
interface in $H^*$.  Together with $T_H$ these trees form a spanning
tree $T$ of $G$.  Apply Lemma~\ref{sbs0} on $T$ to the supports of the
initial and target configurations, in both cases with target support
$V(G)\setminus W$.  At an additive cost $\Ord(N^2)$ for each
configuration, it is therefore enough to transform two configurations
$f^*,g^*$ whose blank vertices are exactly $W$.  In particular, every
vertex outside $H^*$ is occupied in both configurations.

Let
\[
    s:=\max_B |B|
\]
where the maximum is over the nonempty bundles.

We finish by a subsequence dichotomy according to the size of the
largest bundle.  The estimates below will be used on arbitrary
subsequences of the given family.

\paragraph{An unbounded bundle.}
Suppose first that, along the subsequence under consideration,
\[
    s\longrightarrow\infty.
\]
Choose a bundle $B$ of size $s$ with interface $u$.  We may choose $W$
above so that $u\notin W$.  Put
\[
    Y:=G[V(G)\setminus B].
\]
By the preceding local fact, $Y$ is admissible.  Run the Pebble
Presentation Algorithm on $Y$ with interface $u$ and record the
presentation sequence and the correspondence between presentation
positions and vertices of $Y$.

Divide the presentation positions into consecutive blocks of size at
most $s$.  Since $s\to\infty$, there are
\[
    \Ord(N/s+1)
\]
blocks.  Process them exactly as in the batched-transcription proof of
Theorem~\ref{thm:intermediate-cut-component}, except that the bundle
$B$, rather than a feasible subpuzzle, is used as the storage device.
For one block, scan the unprocessed suffix of the recorded presentation
on $Y$ and, whenever a pebble prescribed for the block is presented at
$u$, exchange it with a pebble of $B$ not prescribed for that block by
one push along $T_B$.  Such a pebble always exists until all prescribed
pebbles of the block have been collected, because the block contains at
most $s=|B|$ pebbles.  Reverse the scan to the beginning of the block,
and then advance through the block once more, supplying its prescribed
pebbles from $B$ one by one, again by pushes along $T_B$.
As in the proof of Theorem~\ref{thm:intermediate-cut-component}, all
blocks are written before a single final backward execution of the
recorded presentation sequence realizes all substitutions on $Y$.

The presentation scans cost
\[
    \Ord\!\left(N^2\left(\frac{N}{s}+1\right)\right)
    =\Ord\!\left(\frac{N^3}{s}+N^2\right).
\]
Over all blocks, only $\Ord(N)$ exchanges with $B$ are made.  By
\eqref{eq:shallow-push}, their total cost is
\[
    \Ord(b\varepsilon N^3+N^2).
\]
After the target configuration on $Y$ has been realized, the pebbles
on $B\cup\{u\}$ are exactly the target set for $B\cup\{u\}$.  Root
$T_B$ at $u$ and place these pebbles in their exact target positions in
postorder, including the root $u$ and placing it last.  Thus temporary
changes at the interface made while sorting different islands are
removed at the final step.  This uses at most $|B|+1\leq N+1$ further
pushes and has the same
\[
    \Ord(b\varepsilon N^3+N^2)
\]
bound.  Consequently, along every subsequence for which $s\to\infty$,
\begin{equation}
    L(G)
    \leq
    \Ord\!\left(\frac{N^3}{s}
        +b\varepsilon N^3+N^2\right).
    \label{eq:large-shallow-bundle}
\end{equation}

\paragraph{Uniformly bounded bundles.}
Suppose instead that, along the subsequence under consideration, there
is a constant $M$ such that
\[
    s\leq M.
\]
There is one bundle for each assigned interface, so the number of nonempty
bundles is at most $|V(H^*)|$.  Since every bundle has order at most $M$ and
together the core and the bundles contain all $N$ vertices, it follows that
\[
    N\le (M+1)|V(H^*)|,
\]
and hence $|V(H^*)|=\Theta(N)$.

Ignore the exact positions inside the bundles and record only their
containers.  Apply Lemma~\ref{lem:core-buffer-redistribution} with core
$H^*$ and with the nonempty bundles as the attached containers.  The local
exchange required by the lemma is a push along the fixed tree of the bundle,
using the temporary cycle repair described above whenever necessary.  By the
restoration property of that repair, the core outside the current interface
and every other bundle are returned to their previous states.  Since
every bundle has order at most $M$, every such push has length $\Ord(1)$.
Even if the borrowed repair cycle has length $\Theta(N)$, one local exchange
therefore costs only $\Ord(N)$ moves.

The lemma uses only $\Ord_M(1)$ reconfigurations of $H^*$ and
$\Ord_M(N)$ local exchanges.  Consequently all local exchanges, including
all temporary cycle repairs, cost
\[
    \Ord(N^2)
\]
in total, while all reconfigurations of the core cost
\[
    \Ord_M(1)L(H^*)=\ord(N^3).
\]
After this redistribution every bundle contains its exact target set of
pebbles.  Finally sort the pebbles inside all bundles in postorder on their
trees $T_B$.  The bundles are disjoint and contain at most $N$ vertices in
total, so this uses only $\Ord(N)$ further bounded-length pushes and hence
another $\Ord(N^2)$ moves including cycle repairs.  One final
reconfiguration of $H^*$ completes $g^*$.  Therefore, along every
subsequence on which the bundle sizes are uniformly bounded,
\begin{equation}
    L(G)
    \leq
    \Ord_M(1)L(H^*)+\Ord(N^2)
    =\ord(N^3).
    \label{eq:bounded-shallow-bundles}
\end{equation}

We now conclude.  Suppose, to the contrary, that $L(G)\neq\ord(N^3)$.
Then there are $\eta>0$ and a subsequence, still indexed by $N$, such
that
\[
    L(G)\geq \eta N^3
\]
throughout that subsequence.  Fix $\varepsilon>0$ and perform the
absorption construction above.  For the resulting sequence of largest
bundle sizes $s$, either there is a further subsequence on which $s$ is
uniformly bounded, or there is a further subsequence on which
$s\to\infty$.

In the first case,~\eqref{eq:bounded-shallow-bundles} already gives
$L(G)=\ord(N^3)$, contradicting the choice of the subsequence.  In the
second case,~\eqref{eq:large-shallow-bundle} gives
\[
    \limsup_{N\to\infty}\frac{L(G)}{N^3}
    \leq C_b\varepsilon,
\]
where $C_b$ depends only on the fixed number $b$ of blanks.  Since
$\varepsilon>0$ was arbitrary, choose it so that
$C_b\varepsilon<\eta$, again contradicting the choice of the subsequence.
Thus
\[
    L(G)=\ord(N^3).
\]
This proves the theorem.
\end{Proofof}

%%%%%%%%%%%%%%%%%%%%%%%%%%%%%%%%%%%%%%%%%%%%%%%%%%%%%%%%%%%%%%%%%%%%%%%%%%%%%%%%

\subsection{Consequences}

We have obtained a sequence of increasingly restrictive structural conditions
for a board graph whose worst-case shortest solution length is
$\Theta(N^3)$.  First, no cut vertex can separate a subgraph of order
$\omega(1)\cap\ord(N)$.  Second, every component of
$\mathcal{B}(G)$ has order $\Ord(1)$ or $\Theta(N)$.
Moreover, Theorem~\ref{thm:linear-component-cubic} shows that every
$\Theta(N)$-size component must itself have cubic worst-case reconfiguration
complexity, while Corollary~\ref{cor:bridge-tree-components-bounded} excludes
bridge-tree components of linear order. 
Third, starting either from a linear-size subcubic core or from a core of
unbounded sublinear order, an arbitrary sequence of ear additions cannot
create cubic complexity provided that the resulting graph is biconnected;
there is no restriction on the number or the individual lengths of the ears.
Consequently every unbounded component that remains in the
cubic regime is a $\Theta(N)$-size biconnected component and itself lies in
the cubic regime.

%%%%%%%%%%%%%%%%%%%%%%%%%%%%%%%%%%%%%%%%%%%%%%%%%%%%%%%%%%%%%%%%%%%%%%%%%%%%%%%%
%%%%%%%%%%%%%%%%%%%%%%%%%%%%%%%%%%%%%%%%%%%%%%%%%%%%%%%%%%%%%%%%%%%%%%%%%%%%%%%%
%%%%%%%%%%%%%%%%%%%%%%%%%%%%%%%%%%%%%%%%%%%%%%%%%%%%%%%%%%%%%%%%%%%%%%%%%%%%%%%%

\section{Uniqueness of the Giant Component in the Cubic Regime}\label{lower bound}
%\section{Sufficient condition for $\Theta(N^3)$ solution length}\label{lower bound}

Let $(G, P)$ be a puzzle with the board graph $G$ and the pebble set $P$.
Let $N$ be the number of the vertices of $G$, and let $n$ be the cardinality of $P$. 
The number of blanks of the puzzle is $N-n$. 
Throughout this section, we assume $N-n = \Ord (1)$.

\subsection{$k$-isthmus-avoiding covering}

Let $G$ be a graph with the vertex set $V(G)$ and the edge set $E(G)$.
A connected induced subgraph $H$ of $G$ is called {\it $k$-isthmus-avoiding} if $H$ contains no pair of vertices $x$, $y$ of $G$ such that $x$ and $y$ is contained in a pair of different components of $G \setminus V(I)$ with some $k$-isthmus $I$.
Furthermore, a $k$-isthmus-avoiding subgraph $H$ of $G$ is called {\it maximally $k$-isthmus-avoiding}, if for any additional vertex $v \not\in V(H)$, $G[V(H) \cup \{ v \}]$ is disconnected or it is not $k$-isthmus-avoiding.
By definition, any maximally $k$-isthmus-avoiding subgraph $H$ has at least
$\min\{k+1,|V(G)|\}$ vertices, since any connected induced subgraph
of $G$ with at most $k+1$ vertices has no $k$-isthmus of $G$.

In the following, we assume $k \ge 2$. 
We build a graph $G_0$ from $G$ by deleting all inner edges and inner vertices of all isthmuses with at least $k$ vertices.
Let ${\cal X}_k$ be the set of all connected components of $G_0$.
For $X \in {\cal X}_k$, define $X_k^{+}$ as a set of vertices $V(X) \cup \{v : v \in V(I) {\rm ~for~some~} k$-$ {\rm isthmus~} I {\rm ~with~} V(I) \cap V(X) \ne \emptyset  \}$.  

\begin{la}\label{la-k-isthmus-avoiding}
Let $G$ be a connected graph.
Let $H$ be an induced subgraph of $G$.
Then $H$ is a maximally $k$-isthmus-avoiding subgraph of $G$ if and only if $H$ satisfies one of the following conditions.
\begin{enumerate}
\item{
$H$ is a $(k+1)$-isthmus of $G$.
}
\item{
$H$ is $G[X_k^{+}]$ for some $X \in {\cal X}_k$.
}
\end{enumerate}
\end{la}

\begin{Proofof}{Lemma\tume\ref{la-k-isthmus-avoiding}}
First, we show that if $H$ satisfies {\rm (1)} or {\rm (2)}, then $H$ is maximally $k$-isthmus-avoiding.

Assume that $H$ is a $(k+1)$-isthmus $I$.
$H$ is $k$-isthmus-avoiding, since $|V(H)| \le k+1$.

Consider $G[V(H) \cup \{ u \}]$ with an additional vertex $u$. 
By preserving the connectivity, we have $u v \in E(G)$ with an endvertex $v$ of $I$.
Let $w$ be the other endvertex of $I$ other than $v$.
Then $u$ and $w$ are separated by a $k$-isthmus $I \setminus \{ w \}$.
Hence, $H$ is maximally $k$-isthmus-avoiding.

Assume that $H = G[X_k^{+}]$ for some $X \in {\cal X}_k$.
Every $k$-isthmus $I$ contained in $H$ satisfies $|V(X) \cap V(I)| = 1$, since $X$ has no inner vertex of a $k$-isthmus.
Suppose that $I$ partitions $V(G) \setminus V(I)$ into two parts $L$ and $R$.
Since $V(X) \setminus V(I)$ is contained in one of $L$ and $R$, we have $H$ is not partitioned by $I$.
Hence, $H$ is $k$-isthmus-avoiding.

Consider $G[V(H) \cup \{ u \}]$ with an additional vertex $u$. 
By preserving the connectivity, we have $u v \in E(G)$ with an endvertex $v$ of some $k$-isthmus $I$ with $V(I) \cap V(X) \ne \emptyset$.
Let $x \in V(X)$ be the other endvertex of $I$ other than $v$.
Since $I$ is an isthmus of $G$, we have $|V(X)| \ge 2$ or $x$ has an endvertex of a $k$-isthmus $J$ other than $I$.
Hence, we have $y \in V(H) \setminus V(I)$ such that $xy \in E(H)$.
Then we have $u$ and $y$ are separated by $I$.
Hence, we have $H$ is maximally $k$-isthmus-avoiding.

Conversely, we will show that if $H$ is maximally $k$-isthmus-avoiding, then $H$ satisfies one of conditions {\rm (1)} or {\rm (2)}.

If $V(H) \cap V(X) = \emptyset$ for all $X \in {\cal X}_k$ then $H$ is contained in some isthmus of $G$.
By the maximality of $H$, $H$ is a $(k+1)$-isthmus of $G$.

In the following, we assume $V(H) \cap V(X) \ne \emptyset$ for some $X \in {\cal X}_k$.
\medskip\\
{\rm Claim }$1$. %
Let $X \in {\cal X}_k$.
If $V(H) \cap V(X) \ne \emptyset$,
then $V(X) \subset V(H)$ or $|V(X) \cap V(H)| = 1$.

Assume that $V(X) \not\subset V(H)$ and $|V(X)| \ge 2$.
Then we have $x, y \in V(X)$ such that $x \in V(H)$, $y \not\in V(H)$, and $xy \in E(G)$. 
By the maximality of $H$, $G[V(H) \cup \{ y \}]$ contains a $k$-isthmus $I$ separating $y$ and $w \in V(H)$.
Since $X$ contains no inner vertex of a $k$-isthmus, $V(I) \cap V(X) = \{ x \}$ and $x$ is an endvertex of $I$.
If there exists a vertex $y' \in V(X) \setminus \{ x \}$ with $y' \in V(H)$, then $I$ separates $y'$ and $w$, a contradiction.
Hence, we have $|V(X) \cap V(H)| = 1$.
\medskip\\
{\rm Claim }$2$. %
Let $X \in {\cal X}_k$ such that $V(X) \subset V(H)$ and $|V(X)| \ge 2$.
Then $H = G[X_k^{+}]$.

In this case, for $u \not\in X_k^{+}$, there exists a $k$-isthmus $I \in G[X_k^{+}]$ with an endvertex $x \in V(I) \cap V(X)$ and a vertex $y \in V(X) \setminus V(I)$ such that $I$ separates $u$ and $y$.
Hence, we have $V(H) \subset V(X_k^{+})$.
Since $G[X_k^{+}]$ is $k$-isthmus-avoiding, by the maximality of $H$, we have $H = G[X_k^{+}]$. 

By Claim $1$ and Claim $2$, we assume $|V(X) \cap V(H)| \le 1$ for all $X \in {\cal X}_k$ in the following.

We define an incidence graph ${\cal T}_H$ such that $V({\cal T}_H)$ is $\{ X \in {\cal X}_k : V(X) \cap V(H) \ne \emptyset \}$ and for $X, Y \in V({\cal T}_H)$, $(X, Y) \in E({\cal T}_H)$ if $X$ and $Y$ are linked by an isthmus with at least $k$ vertices.
\medskip\\
{\rm Claim }$3$. %
${\cal T}_H$ is a single vertex or a star with at least two vertices.

If ${\cal T}_H$ contains a cycle, then we have a cycle of $G$ containing a $k$-isthmus, a contradiction.
If ${\cal T}_H$ contains a path with $4$ vertices $X_1 X_2 X_3 X_4$, then we have a pair of vertices $x_1 \in V(X_1)\cap V(H)$ and $x_4 \in V(X_4)\cap V(H)$ and an isthmus $I$ of length at least $k$ linking $X_2$ and $X_3$.
Then $I$ separates $x_1$ and $x_4$, a contradiction.
Hence, we have ${\cal T}_H$ is a star.  
\medskip\\
{\rm Claim }$4$. %
Let $X$ be the center of the star ${\cal T}_H$ with at least three vertices.
Then $H = G[X_k^{+}]$.

In this case, $H$ contains at least two $k$-isthmuses $I_i$ such that $V(I_i) \cap V(X) = \emptyset$ for $i = 1,2$.
Hence, for $u \not\in X_k^{+}$, there exists a $k$-isthmus $I \in G[X_k^{+}]$ and a vertex $y \in V(I_1) \cup V(I_2)$ such that $I$ separates $u$ and $y$.
Hence, we have $V(H) \subset V(X_k^{+})$.
Since $G[X_k^{+}]$ is $k$-isthmus-avoiding, by the maximality of $H$, we have $H = G[X_k^{+}]$. 
\medskip\\
{\rm Claim }$5$. %
If ${\cal T}_H$ is a single edge $(X, Y)$, then  $H$ is one of $G[X_k^{+}]$, $G[Y_k^{+}]$, or a $(k+1)$-isthmus.

Let $I$ be a maximal isthmus linking $X$ and $Y$.
We have $|V(I)| \le k+1$, since $H$ is $k$-isthmus-avoiding.
If $|V(I)|$ is $k+1$, then $H$ is a $(k+1)$-isthmus $I$.
Hence, we may assume $|V(I)|$ is $k$.
Let $x, y$ be the endvertices of $I$ with $x \in X$, $y \in Y$.
Then exactly one of $x$ and $y$ has a neighbor $z \in V(H) \setminus V(I)$.
By symmetry, we assume $xz \in E(H)$.
Then for $u \not\in X_k^{+}$, there exists a $k$-isthmus
$J \in G[X_k^{+}]$ and a vertex
$v \in V(I) \cup \{z\}$ such that $J$ separates $u$ and $v$.
Hence, we have $V(H) \subset V(X_k^{+})$.
Since $G[X_k^{+}]$ is $k$-isthmus-avoiding, by the maximality of $H$, we have $H = G[X_k^{+}]$. 
\medskip\\
{\rm Claim }$6$. %
If ${\cal T}_H$ is a single vertex $X$, then we have $H = G[X_k^{+}]$ or a $(k+1)$-isthmus.

Let $x \in V(X) \cap V(H)$.
If there exists a vertex $y \in V(X) \setminus V(H)$, then we have $z \in V(H)$ such that a $k$-isthmus $I$ contained in $H$ separating $z$ and $y$.
Since $z$ is contained no component of ${\cal X}_k$, $I \cup \{ z \}$ is a $(k+1)$-isthmus. 
     
If $V(X) \subset V(H)$, we have $|X|=1$.
For $u \not\in X_k^{+}$, there exists a $k$-isthmus $I \in G[X_k^{+}]$ and a vertex $y \in  V(H) \setminus V(I)$ such that $I$ separates $u$ and $y$.
Hence, we have $V(H) \subset V(X_k^{+})$.
Since $G[X_k^{+}]$ is $k$-isthmus-avoiding, by the maximality of $H$, we have $H = G[X_k^{+}]$. 
\end{Proofof} 
  
A family of subsets $\{ V_i : 1 \le i \le s \}$ of $V(G)$ is called {\it $k$-isthmus-avoiding covering} of $G$ if $\{ G[V_i] : 1 \le i \le s \}$ is the family of all maximally $k$-isthmus-avoiding subgraphs of $G$.
Note that a $1$-isthmus-avoiding covering is a biconnected component decomposition. 
%Every graph has a unique $k$-isthmus-avoiding covering for all $k \ge 1$.

%\subsection{Puzzle on $k$-isthmus-avoiding covering}
Let us consider a puzzle $(G, P)$ with $k=N-n=\Ord(1)$ blanks on a connected graph $G$.
In the following, we assume that $G$ is not a cycle.
Let $V(G)  = \cup_{1 \le i \le s} V_i $ be the $k$-isthmus-avoiding covering of $G$.
%Let $N_i = |V_i|$ for $1 \le i \le s$. 
Since $G$ is not a cycle and $k \ge 2$, $(G, P)$ is feasible if and only if $s=1$.

Since each $G[V_i]$ has no $k$-isthmus of $G[V_i]$ and $G[V_i]$ is not a cycle, the puzzle with $k$ blanks on the board graph $G[V_i]$ is feasible.

Let $f$ be a configuration of $(G, P)$.
For a member $V_i$ of the $k$-isthmus-avoiding covering of $G$, by a sequence of moves starting from $f$ to a configuration $g$, $k$ blanks can be moved to $V_i$.
Let $P^f_i$ be a set of $|V_i| - k$ pebbles which occupy $V_i$ in $g$.
$P^f_i$ is independent of $g$, since $P^f_i$ and the other pebbles cannot be exchanged by the existence of $k$-isthmuses around $V_i$.      

Hence, we have a decomposition of the pebble set $P = P^f_1 \cup P^f_2 \cup \dots \cup P^f_s$ with $P^f_i \cap P^f_j = \emptyset$ for $i \ne j$
such that $(G[V_i], P^f_i)$ is feasible with $k$ blanks for all $i$.
We call $P^f_i$ a {\it pebble block} with respect to $f$.
A pebble block is called {\it trivial} if it consists of a single pebble. 
%If two pebbles $p$ and $q$ are contained in distinct %pebble blocks $P^f_i$ and $P^f_j$, they cannot be %changed their positions by any sequence of moves %starting from $f$.
We have $|P^f_i| = |P^g_i|$ for any two configurations $f$ and $g$ for all $i$, since $|P^f_i| = |V_i|-k$.
We define a graph $\tilde{T}^f$ such that $V(\tilde{T}^f) = \{ P^f_i : 1 \le i \le s \}$ and two vertices $P^f_i$ and $P^f_j$ are adjacent in $\tilde{T}^f$ if $G[V_i \cap V_j]$ is a $k$-isthmus of $G$, where $P^f_i$ and $P^f_j$ are on the board graphs $G[V_i]$ and $G[V_j]$, respectively.

\begin{la}\label{la-k-isthmus-block-tree}
$\tilde{T}^f$ is a tree.
\end{la}

\begin{Proofof}{Lemma\tume\ref{la-k-isthmus-block-tree}}
Let $P^f_i$ and $P^f_j$ be adjacent in $\tilde{T}^f$.
Then $G[V_i \cap V_j]$ is a $k$-isthmus of $G$, where vertex sets $V_i$ and $V_j$ are corresponding to $P^f_i$ and $P^f_j$.
Hence, if $\tilde{T}^f$ has a cycle, then a set of $V_i$'s are cyclically linked by $k$-isthmuses of $G$.
This contradicts that the definition of isthmus.
Furthermore, since $G$ is connected, $\tilde{T}^f$ is also connected.
Therefore, $\tilde{T}^f$ is a tree.
\end{Proofof}

For two configurations $f$ and $g$, by identifying $P^f_i \in \tilde{T}^f$ and $P^g_i \in \tilde{T}^g$ for all $i$, we have $\tilde{T}^f$ and $\tilde{T}^g$ are isomorphic as graphs.
Hence, we can denote the tree by $\tilde{T}$, which is independent of the choice of a configuration.
We call $\tilde{T}$ a {\it $k$-isthmus block tree}.
For two configurations $f$ and $g$, $f$ can be transformed to $g$ by a sequence of finite moves if and only if $P^f_i = P^g_i$ for all $i$.  

\subsection{Edge-amalgamation onto a cycle} 
 
Let $C$ be a cycle graph with the edge set $E(C) = \{e_1, e_2, \dots, e_r\}$, where each edge $e_{i}$ connects the vertices $v_{i}$ and $v_{i+1}$ with $v_{r+1} = v_1$.
Let ${\cal H} = (H_1, H_2, \dots, H_r)$ be a sequence of connected graphs, where each $H_{i}$ has a designated pair of distinct vertices $(x_i, y_i) \in V(H_i) \times V(H_i)$ such that $x_i y_i \in E(H_i)$ or $\min\{ {\rm deg}_{H_i} x_i, {\rm deg}_{H_i} y_i \} \ge 2$.
A graph $G$ is constructed by identifying the vertex pair $(x_i, y_i)$ of each graph $H_{i}$ with the endpoints $(v_i, v_{i+1})$ of the edge $e_{i} \in E(C)$ for all $i \in \{1, 2, \dots, r\}$.
Hence, we have $V(G) = V(C) \cup \cup_{1 \le i \le r} (V(H_i) \setminus \{ x_i, y_i \})$ and $E(G) = \cup_{1 \le i \le r} E(H_i)$.  
The graph $G$ is called an {\it edge-amalgamation} of ${\cal H}$ onto the cycle $C$. 

\subsection{$\Theta(N^3)$-puzzle} 

Under the asymptotic convention of
Section~\ref{sec:structural-restrictions}, we say that a family of
feasible puzzles $\{(G_N,P_N)\}_{N\to\infty}$ is a
$\Theta(N^3)$-puzzle family, or an $\ord(N^3)$-puzzle family, if its
worst-case shortest solution length is $\Theta(N^3)$, or
$\ord(N^3)$, respectively.  When the underlying family is clear, we
continue to suppress the index $N$ and simply say that $(G,P)$ is a
$\Theta(N^3)$-puzzle or an $\ord(N^3)$-puzzle.

\begin{thm}\label{thm-Theta-n3}
Let $G$ be an edge-amalgamation of a sequence of connected graphs ${\cal H}$ onto a cycle $C_0$.
Let $(G, P)$ be feasible.
%Let $h_i = |V(H_i)| - 2$ for all $i$, and $h = \sum_i h_i$. 
Let $h = |V(G)| - |V(C_0)|$. 
If $h = \Ord(1)$, then $(G, P)$ is a $\Theta(N^3)$-puzzle.
\end{thm}

\begin{Proofof}{Theorem\tume\ref{thm-Theta-n3}}
Take a cycle $C$ of $G$ with $V(C_0) \subseteq V(C)$.
Let $Q$ be a longest corridor of $G$.
Since $h = \Ord(1)$, $Q$ is a part of $C$ and the length of $Q$ is $\Theta(N)$.
Let $Q_0 \subset Q$ be a path obtained by omitting the $k$ vertices on either side of $Q$. 
%The length of $Q_0$ is still $\Theta(N)$. 
Put $t = |V(Q_0)|$.
We have $t = \Theta(N)$.

For an edge $e \in E(Q_0)$
by considering the $k$-isthmus-avoiding covering of 
$G \setminus e $, we have a $k$-isthmus block tree $\tilde{T}(G \setminus e)$.
Since $(G,P)$ is feasible, all $k$-isthmuses of $G\setminus e$ are contained in $C$. 
Since these $k$-isthmuses are arranged in a line, $\tilde{T}(G \setminus e )$ is a path.
Moreover, since we choose $e$ from $Q_0$, both endvertices $u$ and $v$ of the path $\tilde{T}(G \setminus e )$ correspond to trivial pebble blocks $P^f_1(G \setminus e )$ and $P^f_2(G \setminus e )$.
Now, build a cycle $\tilde{C}(G \setminus e )$ by connecting endvertices $u$ and $v$ of $\tilde{T}(G \setminus e )$ by an additional edge.

%Let $e_1=x x_1 \in E(Q_0)$ and $e_2 = x x_2 \in E(Q_0)$ be adjacent in $G$.
%We may assume there is a configuration $f$ such that $x$, $x_1$, $x_2$ are occupied by pebbles $p$, $p_1$, $p_2$, respectively.
%Then endvertices of $\tilde{T}(G \setminus e_i)$ corresponds to trivial pebble blocks $\{ p \}$ and $\{ p_i \}$ for $i=1, 2$.
%Hence, $\tilde{T}(G \setminus e_2 )$ is built by removing $\{ p \}$ from one end of $\tilde{T}(G \setminus e_1 )$ and move it to the other end of $\tilde{T}(G \setminus e_1 )$.
%Hence $\tilde{C}(G \setminus e )$ is independent of the choice of $e \in E(Q_0)$.
%We denote this cycle by $\tilde{C}$.

For a configuration $f$, an edge $e \in E(Q_0)$, and two pebbles $p, q \in P$, we define a weight $w^f_e (p,q)$ as the distance between $P^f_i$ and $P^f_j$ in $\tilde{C}(G \setminus e )$, where $p \in P^f_i$ and $q \in P^f_j$.

For a configuration $f$ and an edge $e \in E(Q_0)$, let us define a weight $W_e (f) = \sum_{1 \le i \le m} w^f_e (p_i, q_i)$, where $p_1, p_2, \ldots, p_m$ and $q_1, q_2, \ldots, q_m$ are pebbles of $P$ with $m = \lfloor t/2 \rfloor$.

Let $f_1$ be a configuration such that $p_1, p_2, \ldots, p_m, q_1, q_2, \ldots, q_m$ are arranged on $V(Q_0)$ in this order.
Since all these pebbles are contained in trivial pebble blocks, we have
$W_e (f_1) = m^2$ for all $e \in E(Q_0)$.
%Since $t = \Theta(N)$, we have $W_e (f_1) = \Theta(N^2)$.
Let $f_2$ be a configuration such that $p_1, q_1, p_2, q_2, \ldots, p_m, q_m$  are arranged on $V(Q_0)$ in this order.
Then we have $W_e (f_2) = m$ for all $e \in E(Q_0)$.
Hence, we have $W_e(f_1) - W_e(f_2) = m^2 - m = \Theta(N^2)$ for all $e \in E(Q_0)$.

We claim that the number of any solution sequence of moves from $f_1$ to $f_2$ is $\Omega(N^3)$.

Suppose, for a contradiction that there exists a sequence of moves $M$ from $f_1$ to $f_2$ with $|M| = \ord(N^3)$.
Since $|E(Q_0)| = \Theta(N)$, there exists an edge $e \in E(Q_0)$ such that $e$ is traversed at most $|M|/|E(Q_0)| = \ord(N^2)$ times. 

For a move $m$ not using $e$,  all the pebble blocks of $\tilde{C}(G \setminus e)$ are unchanged, and the weight $W_e$ is not changed.
On the other hand, for a move $m$ using $e$, each pebble may cyclically go from a current pebble block to an adjacent pebble block.
\medskip\\
Claim. Let $e=uv \in E(Q_0)$. 
Let $m$ be a move which transfers a pebble $p_0$ from $u$ to $v$, and let $f$ and $g$ denote two configurations immediately before and after $m$.
Let $\{ P^f_i \}$ and $\{ P^g_i \}$ be families of pebble blocks according to a $k$-isthmus-avoiding covering $\{ V_i \}$ such that $P^f_1 P^f_2 \ldots P^f_\alpha$ and $P^g_1 P^g_2 \ldots P^g_\alpha$ are arranged along a path $\tilde{T}(G \setminus e)$ with $P^f_\alpha = P^g_1 = \{ p_0 \}$.
Then, we have $P^g_i = (P^f_i \setminus \{ p_i \} ) \cup \{ p_{i-1} \}$ with some $p_i \in P^f_i$ for $1 \le i \le \alpha$, where $P^f_0 = P^f_\alpha$.

Since $p_0$ is moved to $v$ by $m$, we have $P^g_1 = \{ p_0 \}$, and $p_1 \in P^f_1$ is pushed out to the next pebble block $P^f_2$.
Hence, we have $P^g_2 = (P^f_2 \setminus \{ p_2 \}) \cup \{ p_1 \}$, with some pebble $p_2 \in P^f_2$.
In the same way, a pebble $p_{i-1} \in P^f_{i-1}$ is pushed out to $P^g_i$ one by one, we have $P^g_i = (P^f_i \setminus \{ p_i \} ) \cup \{ p_{i-1} \}$ for $1 \le i \le \alpha$.

By the above claim, for two pebbles $p$ and $q$, $w^f_e (p,q)$ is changed at most $1$, and moreover, if both $p$ and $q$ are contained in some trivial blocks both before and after $m$, then $w^f_e (p,q)$ is not changed.

Since the number of pebbles belonging to nontrivial blocks is at most $\Ord(h) = \Ord(1)$, at most $\Ord(1)$ of the chosen pairs can change their weight.
Hence, the difference of $W_e$ before and after $m$ is at most $\Ord(1)$.
Therefore, $W_e$ decreases at most $\ord(N^2) \Ord(1) = \ord(N^2)$ throughout $M$, which contradicts $W_e(f_1) - W_e(f_2) = \Theta(N^2)$.

Now, we have $\Omega(N^3)$ moves are necessary for reconfiguration.
Together with the general $\Ord(N^3)$ upper bound for feasible PMP instances, this proves the theorem.
\end{Proofof} 

We conjecture the converse of Theorem\tume\ref{thm-Theta-n3} is true.
\begin{conj}\label{conj-Theta-n3}
If $(G, P)$ is $\Theta(N^3)$-puzzle, then $G$ is an edge-amalgamation of a sequence of connected graphs ${\cal H}$ onto a cycle $C_0$ such that $|V(G)| - |V(C_0)| = \Ord(1)$.
\end{conj}

In the following, we will show some results supporting Conjecture\tume\ref{conj-Theta-n3}.

\begin{thm}\label{thm-cycle-o-n3}
Let $G$ be an edge-amalgamation of a sequence of connected graphs ${\cal H}$ onto a cycle $C_0$ with $\Theta(N)$ vertices.
Let $(G, P)$ be feasible.
Let $h = |V(G)| - |V(C_0)|$ and $h_i = |V(H_i) \setminus V(C_0)|$ for $H_i \in {\cal H}$. 
If $h = \omega(1)$ and, for every $i$, $h_i = \Ord(1)$ and $H_i$ has no $k$-isthmus, then $(G, P)$ is an $\ord(N^3)$-puzzle.
\end{thm}

\begin{Proofof}{Theorem\tume\ref{thm-cycle-o-n3}}
Put $V_{\cal H} = V(G) \setminus V(C_0)$.
Then $|V_{\cal H} | = h = \omega(1)$, while $|V(C_0)|=\Theta(N)$.
Fix arbitrary initial and target configurations.

First we assume $h = \ord(N)$.

Divide $V(C_0)$ into $S_1,S_2,\ldots,S_m$ such that $|S_i| \le h$ and $|S_i| = \Theta(h)$ for all $i$.
Thus $m=\Ord(N/h)$.
For each $i$, let $P_i$ be the set of pebbles that must finally occupy $S_i$.

Throughout the proof, we use a rotation process on $C_0$ that shifts all pebbles on $C_0$ one by one without changing a position of pebbles on $V_{\cal H} \setminus V(C_0)$.
Although $C_0$ is divided into segments by $H_i \in {\cal H}$, each $H_i$ becomes a feasible subpuzzle having at most $\Ord(h_i) = \Ord(1)$ vertices with $k$ blanks when nearby vertices are added.
Hence, each reconfiguration process on $H_i$ can be also performed with $\Ord(1)$ moves.
Since $|{\cal H}| = \Ord(N)$, one shift of a rotation process on $C_0$ can be achieved with $\Ord(N)$ moves, which contains moves necessary for one shift of a rotation process on a 
genuine cycle with $\Ord(N)$ vertices.
Therefore, a rotation of all the pebbles on $C_0$ going through all the vertices of $C_0$ can be achieved with $\Ord(N)\Ord(N) = \Ord(N^2)$ moves.

We process the sets of pebbles in the order $P_1,P_2,\ldots,P_m$. 
For the processing of $P_i$, 
first, transfer all pebbles of $P_i$ to $V_{\cal H}$ by a rotation of all the pebbles on $C_0$ going through all the vertices of $C_0$.
When a pebble $p \in P_i$ on $C_0$ passes through $H \in {\cal H}$ with $q \in P \setminus P_i$ on $V(H) \setminus V(C_0)$, exchange $p$ and $q$ by using a subpuzzle according to $H$.
Since $|P_i| \le h$, all pebbles of $P_i$ can be moved to $V_{\cal H}$. 
  
Next transfer $P_i$ to their target positions, again by a rotation of all the pebbles on $C_0$ going through all the vertices of $C_0$.
For a pebble $p \in P_i$, let $q \in P \setminus P_i$ be a pebble such that $q$ occupies a target position of $p$.
When $q$ passes through $H \in {\cal H}$ with $p \in V(H) \setminus V(C_0)$, exchange $q$ and $p$ by using a subpuzzle according to $H$.

Iterating this process for $P_1$, $P_2$, $\ldots$, $P_m$, we have completed reconfiguration for $V(C_0)$.

We apply the following claim for the final procedure for reconfiguration of $V_{\cal H}$.
\medskip\\
Claim. Two pebbles $p_1$ and $p_2$ on $V_{\cal H}$ can be swapped while keeping the other pebbles in their original positions by $\Ord(N^2)$ moves. 

Assume that $p_i$ occupies $V(H_i) \setminus V(C_0)$ for $i=1,2$.
First, by using a subpuzzle according to $H_1$, we move $p_1$ to $C_0$ with $\Ord(1)$ moves.
Secondly, move $p_1$ closer to $H_2$ by a rotation process on $C_0$ with $\Ord(N^2)$ moves.
Thirdly, exchange $p_1$ and $p_2$ by using a subpuzzle according to $H_2$ with $\Ord(1)$ moves.
Then by inverse of the first two sequence of moves, we can complete to exchange $p_1$ and $p_2$ by $\Ord(N^2)$ moves. 

By applying the above claim $|V_{\cal H}|-1$ times, we have completed reconfiguration on $V_{\cal H}$, which completes the target configuration on $G$.

Let us estimate the total number of moves for reconfiguration.  
There are $m=\Ord(N/h)$ blocks.
For each set $P_i$, the process can be achieved  by a rotation of all the pebbles on $C_0$ going through all the vertices of $C_0$ two times.
Hence, the procedure needs $\Ord(N^2)$ steps for each block.
Therefore, reconfiguration for these transfer operations contribute $\Ord(N^2) \Ord(N/h) = \ord(N^3)$.

Since we can reconfigure $V_{\cal H}$ by swapping two pebbles of $V_{\cal H}$ at most $|V_{\cal H}| - 1$ times, Reconfiguration of $V_{\cal H}$ needs $\Ord(h) \Ord(N^2) = \ord(N) \Ord(N^2) = \ord(N^3)$ moves. 

Hence the complete solution has length $\ord(N^3)$.

Next we assume $h = \Theta(N)$.

Let $\ell = \lceil h/|V(C_0)| \rceil$.
Divide $V_{\cal H}$ into $S'_0,S'_1,\ldots,S'_\ell$ such that $|S'_i| \le |V(C_0)|$ and $|S'_i| = \Theta(N)$ for $1 \le i \le \ell$ and $|S'_0| = \omega(1)$ and $|S'_0| = \ord(N)$.
For each $i$, let $P'_i$ be the set of pebbles that must finally occupy $S'_i$.

We process the sets of pebbles in the order $P'_1,P'_2,\ldots,P'_\ell$. 
first, transfer all pebbles of $P'_i$ to $V(C_0)$ by a rotation of all the pebbles on $C_0$ going through all the vertices of $C_0$.
When a pebble $q \in P \setminus P'_i$ on $C_0$ passes through
$H \in {\cal H}$ with $p \in P'_i$ on $V(H) \setminus V(C_0)$, 
exchange $p$ and $q$ by using a subpuzzle according to $H$.
Since $|P'_i| \le |V(C_0)|$, all pebbles of $P'_i$ can be moved to $V(C_0)$. 
  
Next transfer $P'_i$ to their target positions, again by a rotation of all the pebbles on $C_0$ going through all the vertices of $C_0$.
For a pebble $p \in P'_i$, let $q \in P \setminus P'_i$ be a pebble such that $q$ occupies a target position of $p$.
When $p$ passes through $H \in {\cal H}$ with $q \in V(H) \setminus V(C_0)$, exchange $q$ and $p$ by using a subpuzzle according to $H$.

Iterating this process for $P'_1$, $P'_2$, $\ldots$, $P'_\ell$, we have completed reconfiguration for $V_{\cal H} \setminus S'_0$ by a sequence of $\Ord(N^2)$ moves.

Finally, for a reconfiguration on $S'_0 \cup V(C_0)$, we can apply the same argument in the case where $h=\ord(N)$, since in each $H \in {\cal H}$, a reconfiguration on $V(H)$ can be done with $\Ord(1)$ moves.
Hence, the complete solution has length $\ord(N^3)$.
\end{Proofof} 

A board graph $G$ in the assumption of Theorem\tume\ref{thm-Theta-n3} contains a cycle of length $\Theta(N)$ such that any other cycle $C'$ of length $\Theta(N)$ satisfies $|V(C') \setminus V(C)| = \Ord(1)$.
The following theorem states that if $G$ contains two cycles meeting in a single vertex, then $(G, P)$ is an $\ord(N^3)$-puzzle.

\begin{thm}\label{thm-two-cycles-o-n3}
Let $H_1$, $H_2$ be cycles or board graphs of feasible puzzles with $k$ blanks.
Let $N_i = |V(H_i)|$ with $N_1 \le N_2$. 
Let $G$ be a graph constructed by identifying $u_1 \in V(H_1)$ and $u_2 \in V(H_2)$.
Then the worst-case shortest solution length of the PMP on $G$ is 
$\displaystyle \Ord (\frac{N^3}{N_1} + N^2 \log{N_1})$.
\end{thm}

\begin{Proofof}{Theorem\tume\ref{thm-two-cycles-o-n3}}
Let us denote the identified vertex $u_1 = u_2$ as $u$.
Fix an arbitrary feasible instance on $G$.

By applying a Pebble Presentation Algorithm on $H_i$ with interface $u$, a linear order of positions of $V(H_i)$ is determined for $i=1, 2$.
Let us denote the set of the first $t$ vertices, or the last $t$ vertices, of $H_i$ on this order by $V^{\rm first}_i(t)$, or $V^{\rm last}_i(t)$, respectively.  
Put $\alpha = \lceil (N-1)/(N_1-1) \rceil - 1$.
We have $\alpha \ge 1$.
Divide $V(G)$ into $S_1$, $S_2$, $\ldots$, $S_\alpha$ such that $|S_i| = N_1-1$ for $1 \le i \le \alpha - 1$ such that $S_1, S_2, \ldots, S_\alpha$ is consecutively arranged according to the presentation order of Pebble Presentation Algorithm.
%Note that $|S_\alpha| \le 2(N_1-1)$.
  
For each $i$, let $P_i$ be the set of pebbles that must finally occupy $S_i$.
We process the pebbles in the order $P_1,P_2, \ldots, P_\alpha$.  At the beginning of the processing of $P_i$, we maintain $P_1,\ldots, P_{i-1}$ has already been arranged in its correct order.

We provide two algorithms, which are applications of Pebble Presentation Algorithm.
In Algorithm \ref{alg:two-parts-separation} and \ref{alg:two-parts-sort}, two presentation procedures have a common interface vertex.
The algorithms are excuted according to transcription mechanism and interface-clearing convention described in Subsection \ref{subsec:transcription-cut}.

%%%%%%%%%%%%%%%%%%%%%%%%%%%%%
\begin{alg}[
$\Alg.\mathrm{TwoPartsSeparation}
(G,Q_0,\prec,f)$]\hfill
\label{alg:two-parts-separation}

\begin{namelist}{OUTPUT}

\item[\sl Input]
A board graph $G$ satisfying the assumption of Theorem \ref{thm-two-cycles-o-n3}.
Pebble sets $Q_0$.
A linear order $\prec$ on $P$. 
A configuration $f$ such that $f(V^{\rm first}_1(t_1) \cup V^{\rm first}_2(t_2)) = Q_0$, where $|Q_0| = t_1 + t_2$.
Put $Q_0 = Q_1 \cup Q_2$ such that $|Q_i|=t_i$ and $p_1 \prec p_2$ for any pair $p_1 \in Q_1$ and $p_2 \in Q_2$.

\item[\sl Output]
A sequence of moves which transfer $f$ to a configuration $g$ such that $g(V^{\rm first}_i(t_i)) = Q_i$ for $i=1, 2$ 
and
$g(v)=f(v)$ for all
$v\notin
V^{\rm first}_1(t_1)\cup V^{\rm first}_2(t_2)$.

\end{namelist}

\begin{alg-enumerate}

\item
\label{alg:two-parts-separation-step0}
Set $R_i$ as the set of pebbles on $V^{\rm first}_i(t_i)$ for $i=1,2$.

\item
\label{alg:two-parts-separation-step1}
If $R_i = \emptyset$ with some $i$, then run Pebble Presentation Algorithm for $H_{3-i}$ forward $|R_{3-i}|$ steps.
Reset $R_{3-i} = \emptyset$.
Goto Step \maru{\ref{alg:two-parts-separation-last-step}}.

\item
\label{alg:two-parts-separation-step2}
Let $r_i \in R_i$ be the current pebble on Pebble Presentation Algorithm for $H_i$ for $i=1,2$.
If $r_1 \in Q_1$ or $r_2 \in Q_2$, let $j$ be the index with $r_j \in Q_j$, and run Pebble Presentation Algorithm for $H_j$ forward $1$ step, and set $R_j = R_j \setminus \{ r_j \}$.
If $r_1 \in Q_2$ and $r_2 \in Q_1$, then swap their positions with a subpuzzle around $u$, run Pebble Presentation Algorithm for $H_1$ forward $1$ step, and set $R_1 = R_1 \setminus \{ r_1 \}$, $R_2 = (R_2 \setminus \{    r_2 \}) \cup \{ r_1 \}$.
Goto Step \maru{\ref{alg:two-parts-separation-step1}}.

\item
\label{alg:two-parts-separation-last-step}
Now, $Q_i$ occupies $V^{\rm last}_i(t_i)$ for $i=1,2$.
Run Pebble Presentation Algorithm for $H_1$ backward $t_1$ steps and 
Run Pebble Presentation Algorithm for $H_2$ backward $t_2$ steps.
Return the sequence of moves performed by the subroutine.

\end{alg-enumerate}
\end{alg}
%%%%%%%%%%%%%%%%%%%%%%%%%%%%%

Let $t = |Q_0|$.
The number of moves necessary for executing $\Alg.\mathrm{TwoPartsSeparation}$ is $\Ord(tN)$, since the algorithm applies Pebble Presentation Algorithm $\Ord(|Q_1| + |Q_2|)$ steps.

%%%%%%%%%%%%%%%%%%%%%%%%%%%%%
\begin{alg}[
$\Alg.\mathrm{TwoPartsSort}
(G,Q_0,\prec,f)$]\hfill
\label{alg:two-parts-sort}

\begin{namelist}{OUTPUT}

\item[\sl Input]
A board graph $G$ satisfying the assumption of Theorem \ref{thm-two-cycles-o-n3}.
A pebble set $Q_0$. 
A linear order $\prec$ on $P$.
A configuration $f$ such that $f(V^{\rm first}_1(t_1) \cup V^{\rm first}_2(t_2)) = Q_0$, where $|Q_0| = t_1 + t_2$,  $|t_2 - t_1| \le 1$.

\item[\sl Output]
A sequence of moves which transfer $f$ to a configuration $g$ such that $g(V^{\rm first}_1(t_1) \cup V^{\rm first}_2(t_2)) = Q_0$, $Q_0$ is sorted according to $\prec$, and $g(v) = f(v)$ for all $v \not\in V^{\rm first}_1(t_1) \cup V^{\rm first}_2(t_2)$.

\end{namelist}

\begin{alg-enumerate}

\item
\label{alg:two-parts-sort-step0}
Call Alg. \ref{alg:two-parts-separation} $(G,Q_0,\prec,f)$.
Set $Q_i = g(V^{\rm first}_i(t_i))$ for $i=1,2$, where $g$ is a current configuration.

\item
\label{alg:two-parts-sort-step1}
Combining Pebble Presentation Algorithm for $H_1$ forward and for $H_2$ backward, move $s_2$ pebbles of $Q_1$ to $V^{\rm first}_2(s_2)$ and move $s_1$ pebbles of $Q_1$ to $V^{\rm first}_1(s_1)$, where $s_1 + s_2 = t_1$, $|s_2-s_1| \le 1$.
Let $g$ be a current configuration.
Call Alg. \ref{alg:two-parts-sort} $(G,Q_1,\prec,g)$.

\item
\label{alg:two-parts-sort-step2}
Combining Pebble Presentation Algorithm for $H_1$ backward and for $H_2$ forward, move $s_1$ pebbles of $Q_2$ to $V^{\rm first}_1(s_1)$ and move $s_2$ pebbles of $Q_2$ to $V^{\rm first}_2(s_2)$, where $s_1 + s_2 = t_2$, $|s_2-s_1| \le 1$.
Let $g$ be a current configuration.
Call Alg. \ref{alg:two-parts-sort} $(G,Q_2,\prec,g)$.

\item
\label{alg:two-parts-sort-step3}
Combining Pebble Presentation Algorithm for $H_1$ forward and for $H_2$ backward, move $Q_1$ to $V^{\rm first}_2(t_1)$ and move $Q_2$ to $V^{\rm first}_1(t_2)$.
Return the sequence of moves performed by the subroutine.

\end{alg-enumerate}
\end{alg}
%%%%%%%%%%%%%%%%%%%%%%%%%%%%%

Let $f(t)$ be the number of moves for executing  $\Alg.\mathrm{TwoPartsSort}$, where $t=|Q_0|$.
The main part of $\Alg.\mathrm{TwoPartsSort}$ is to call itself two times with a half number of pebbles of the input size.
The number of moves required for the preparatory and post-processing tasks involved in calling itself is $\Ord(tN)$.
Hence, we have $f(t) \le 2f(\lceil t/2 \rceil) + \Ord(tN)$.
By the inequality, we have $f(t) = \Ord(Nt\log t)$.   

We now describe the processing of $P_1$.

By $\Alg.\mathrm{TwoPartsSeparation}(G,P,\prec,f)$ with the initial configuration $f$, $P_1$ is gathered on $V(H_1)$.
By combining Pebble Presentation Algorithm for $H_1$ forward and for $H_2$ backward, move a half of $P_1$ to $V^{\rm first}_1(t_1)$ and another half of $P_1$ to $V^{\rm first}_2(t_2)$ with $t_1 + t_2 = |P_1|$.
Call $\Alg.\mathrm{TwoPartsSort}(G,P_1,\prec,f)$ with a current configuration $f$.
Now, $P_1$ is sorted according to $\prec$.
Again by combining Pebble Presentation Algorithm for $H_1$ backward and for $H_2$ forward, move sorted $P_1$ back to $V(H_1)$.
Restore $P_1$ to $V^{\rm last}_2(|P_1|)$ by combining Pebble Presentation Algorithm for $H_1$ forward and for $H_2$ backward.

Next we describe the processing of $P_i$ for $2 \le i \le \alpha$.

Let $D = P_1 \cup P_2 \cup \cdots \cup P_{i-1}$.
$D$ is the set of pebbles already sorted according to $\prec$ and restored in $V^{\rm last}_2(|D|)$.
If $i < \alpha$, by $\Alg.\mathrm{TwoPartsSeparation}(G,P\setminus D,\prec,f)$ with a current configuration $f$, $P_i$ is gathered on $V(H_1)$.
By combining Pebble Presentation Algorithm for $H_1$ forward and for $H_2$ backward, move a half of $P_i$ to $V^{\rm first}_1(t_1)$ and another half of $P_i$ to $V^{\rm first}_2(t_2)$ with $t_1 + t_2 = |P_i|$.
Call $\Alg.\mathrm{TwoPartsSort}(G,P_i,\prec,f)$ with a current configuration $f$.
Now, $P_i$ is sorted according to $\prec$.
Again by combining Pebble Presentation Algorithm for $H_1$ backward and for $H_2$ forward, move sorted $P_i$ back to $V(H_1)$.

If $i<\alpha$, restore $P_i$ to
$V_2^{\rm last}(|P_i|)$ by combining the Pebble Presentation Algorithm
for $H_1$ forward and that for $H_2$ backward.

For $i=\alpha$, no such restoration is necessary.  At this stage
$D=P_1\cup\cdots\cup P_{\alpha-1}$ already occupies, in its prescribed
order, the presentation positions corresponding to
$S_1,\ldots,S_{\alpha-1}$.  Hence the remaining presentation positions
are exactly those corresponding to $S_\alpha$, and the pebbles
occupying them are exactly $P_\alpha=P\setminus D$.  Sorting these
remaining pebbles according to the target presentation order therefore
produces the prescribed target configuration on all of $G$.

Let us estimate the number of moves of a sequence for a reconfiguration.

For the processing of $P_i$ with $1 \le i \le \alpha$, the number of moves to gather $P_i$ to $V(H_i)$ is $\Ord(N^2)$, since it applies $\Alg.\mathrm{TwoPartsSeparation}(G, P\setminus D,\prec, f)$.
Restoring $P_i$ to $V^{\rm first}_1(|P_i|/2)$ and $V^{\rm first}_2(|P_i|/2)$ contributes $\Ord(|P_i| N) = \Ord(N_1 N)$.
$\Alg.\mathrm{TwoPartsSort}(G, P_i, \prec, f)$ needs $\Ord(N |P_i| \log |P_i| )= \Ord(NN_1 \log N_1)$ moves.
Restoring $P_i$ to $V^{\rm last}_2(|P_i|)$ contributes $\Ord(NN_1)$.
Therefore, the processing of $P_i$ needs $\Ord(N^2 + NN_1 \log N_1)$ moves. 
Since we have $\alpha = \Ord(N/N_1)$ blocks, the total number of moves of the sequence is $\Ord(N/N_1) \Ord(N^2 + NN_1 \log N_1) = \Ord(N^3/N_1 + N^2 \log N_1)$.
\end{Proofof}

We have an immediate consequence from Theorem \ref{thm-two-cycles-o-n3}. 

\begin{cor}\label{cor-two-cycles-o-n3}
Let $H_1$, $H_2$ be cycles or board graphs of feasible puzzles with $k$ blanks.
Let $N_i = |V(H_i)|$ with $N_1 \le N_2$. 
Let $G$ be a graph constructed by identifying $u_1 \in V(H_1)$ and $u_2 \in V(H_2)$.
If $N_1 = \omega(1)$, then $(G,P)$ is an $\ord (N^3)$-puzzle.
\end{cor}

The next theorem shows that upper bound of reconfiguration cost in Theorem \ref{thm-two-cycles-o-n3} is approximately best if a board graph contains two cycles with $\Theta(N)$ vertices meeting in a single vertex.

\begin{thm}\label{thm-two-cycles-Theta-n2logn}
Let $H_1$, $H_2$ be cycles such that $N_i = |V(H_i)|$. 
Let $G$ be a graph constructed by identifying $u_1 \in V(H_1)$ and $u_2 \in V(H_2)$.
If $N_i = \Theta(N)$ for $i = 1,2$, then the worst-case shortest solution length of the PMP on $G$ is $\Theta(N^2 \log N)$.
\end{thm}

\begin{Proofof}{Theorem\tume\ref{thm-two-cycles-Theta-n2logn}}
Let us denote the vertex $u_1 = u_2$ as $u$.
$G$ has two corridors $Q_i$ such that $V(Q_i) = V(H_i) \setminus \{ u \}$ for $i=1,2$.
Let $Q_{i,0} \subset Q_i$ be a path obtained by omitting the $k$ vertices on either side of $Q_i$. 
The length of $Q_{i, 0}$ is $\Theta(N)$. 

For a pair of edges $(e_1, e_2) \in E(Q_{1, 0}) \times E(Q_{2, 0})$,
consider the $k$-isthmus-avoiding covering of $G \setminus \{ e_1, e_2 \}$ and a pebble block tree $\tilde{T}(G \setminus \{ e_1, e_2 \})$.
Then there exists a single nontrivial pebble block corresponding to a set of $4k-3 = \Ord(1)$ vertices around $u$ with $k$ blanks, and all other pebble blocks are trivial.

Fix an initial configuration $f_0$.
Define ${\cal M}_s$ as the family of all sequences of at most $s$ moves starting from $f_0$, and define ${\cal F}_s$ as the family of all configurations generated by ${\cal M}_s$.

We will show that if $s = \ord(N^2 \log N)$, there exists a configuration $g$ such that $g \not\in {\cal F}_s$.

For a pair of edges $(e_1, e_2) \in E(Q_{1, 0}) \times E(Q_{2, 0})$, define ${\cal M}(e_1, e_2; t)$ as the family of sequences of moves $M$ starting from $f_0$ such that $e_1$ and $e_2$ are traversed at most $t$ times in total in $M$.

For a move $m$ using neither $e_1$ nor $e_2$, all the pebble blocks of $\tilde{T}(G \setminus \{e_1, e_2\})$ are unchanged immediately before and after $m$.
On the other hand, for a move $m$ using $e_i$, each pebble may cyclically go from a current pebble block to an adjacent pebble block along the cycle $H_i$ for $i=1,2$. 
Since $\tilde{T}(G \setminus \{e_1, e_2\})$ contains only one nontrivial pebble block of order $\Ord(1)$, the number of all possible families of pebble blocks generated by ${\cal M}(e_1, e_2; 1)$, denoted by $c$, is $\Ord(1)$.
Iterating the process at most $t$ times, the number of all possible families of pebble blocks generated by ${\cal M}(e_1, e_2; t)$ is $\Ord(c^t)$.

Each family of pebble blocks of $\tilde{T}(G \setminus \{e_1, e_2\})$ corresponds to at most $\Ord(N^k)$ configurations, since only $\Ord(1)$ pebbles are permuted in the unique nontrivial subpuzzle and there exist $k$ blanks movable on $V(G)$.
Hence, the number of all possible configurations generated by ${\cal M}(e_1, e_2; t)$ is $\Ord ( c^t N^k )$.

Since $N_i=\Theta(N)$ for $i=1,2$ and $k=\Ord(1)$, there exists
a constant $\alpha>0$ such that
\[
    |E(Q_{i,0})|\geq \alpha N
    \qquad (i=1,2)
\]
for all sufficiently large $N$.
For each $M\in\mathcal{M}_s$ and each $i\in\{1,2\}$, the total
number of traversals of the edges of $Q_{i,0}$ is at most $s$.
Hence, by averaging, there exists an edge
$e_i\in E(Q_{i,0})$ that is traversed at most
$s/(\alpha N)$ times in $M$.
Therefore
\[
    M\in
    \mathcal{M}\!\left(
        e_1,e_2;
        \left\lceil\frac{2s}{\alpha N}\right\rceil
    \right).
\]
Hence, we have
\begin{eqnarray*}
\log |{\cal F}_s|
& \le &
\log \sum_{(e_1,e_2)\in E(Q_{1,0})\times E(Q_{2,0})}
\Ord\!\left(
c^{\left\lceil \frac{2s}{\alpha N}\right\rceil}N^k
\right)
\\
& \le &
\log\!\left(
N^2
\Ord\!\left(
c^{\left\lceil \frac{2s}{\alpha N}\right\rceil}N^k
\right)
\right)
\\
& = &
\Ord\!\left(\frac{s}{N}+\log N\right).
\end{eqnarray*}
Since $s=\ord(N^2\log N)$, we have
\[
    \log |{\cal F}_s|=\ord(N\log N).
\]

On the other hand, the number of all configurations, denoted by
${\cal F}$, is $(N!)/(k!)$.
Since every edge of $G$ lies on one of the two cycles, $G$ has no
bridge and hence no $k$-isthmus. Moreover, $G$ is not a cycle.
Therefore, since $k\geq 2$, Theorem~A implies that every pair of
configurations is mutually reachable. Hence
\[
    \log |{\cal F}|
    =
    \log\frac{N!}{k!}
    =
    \Theta(N\log N).
\]
Since $|{\cal F}_s|<|{\cal F}|$, there exists a target configuration
$g$ that cannot be reached from $f_0$ by any sequence
of at most $s=\ord(N^2\log N)$ moves.

Thus, the worst-case shortest solution length is
$\Omega(N^2\log N)$.
Together with the $\Ord(N^2 \log N)$ upper bound given by Theorem\tume\ref{thm-two-cycles-o-n3}, this proves the theorem.
\end{Proofof} 

\begin{thm}\label{thm-Theta-n3-unique-large-component}
Let $\{(G_N,P_N)\}_{N\to\infty}$ be a family of feasible
puzzles whose worst-case shortest solution length is $\Theta(N^3)$.
For all sufficiently large $N$, use the notation of the asymptotic
convention above and write $G=G_N$ and $n=|P_N|$.  Then
$\mathcal{B}(G)$ has the following properties:
\begin{enumerate}
\item{There exists a biconnected component $X$ such that $|V(X)|=\Theta(N)$ 
and the reconfiguration cost of $X$ is $\Theta(N^3)$. }\label{giant-component}
\item{$X$ contains no cycle $C$ with $|V(C)|=\omega(1)\cap \ord(N)$.}
\item{If $C_1,C_2\subseteq X$ are cycles of order $\Theta(N)$, then
$|V(C_1)\cap V(C_2)|\ge N-n$ or $V(C_1)\cap V(C_2) = \emptyset$.}
\item{For all components $Y$ except $X$, we have $|V(Y)| = \Ord(1)$.}
\item{For every cut vertex $v\in V(X)$, the union of all components of
$\mathcal{B}(G)$ attached to $X$ through $v$ has order $\Ord(1)$.}
\end{enumerate}
\end{thm}

\begin{Proofof}{Theorem~\ref{thm-Theta-n3-unique-large-component}}

Suppress the subscript $N$ according to the asymptotic convention above.
Put
\[
    b:=N-n.
\]
By Theorem~\ref{Restriction}, we have $b=\Ord(1)$, and $G$
contains a cycle of order $\Theta(N)$.  Such a cycle is contained in
a biconnected component of $\mathcal{B}(G)$.  Hence
$\mathcal{B}(G)$ contains at least one component $X$ satisfying
\[
    |V(X)|=\Theta(N).
\]
By Corollary~\ref{cor:bridge-tree-components-bounded}, every
bridge-tree component of $\mathcal{B}(G)$ has order $\Ord(1)$.
Thus $X$ is biconnected.  Moreover,
Theorem~\ref{thm:linear-component-cubic} gives
\[
    \text{the reconfiguration cost of }X=\Theta(N^3).
\]
This proves the existence asserted in~\ref{giant-component}.

We next prove that $X$ contains no cycle of order
$\omega(1)\cap\ord(N)$.
Suppose, to the contrary, that $X$ contains such a cycle $C$.
Since $X$ is biconnected, the ear-decomposition theorem allows us to
start an ear decomposition of $X$ from $C$.  Thus $X$ is obtained
from $C$ by successively adjoining ears.

Since
\[
    |V(C)|=\omega(1)\cap\ord(N),
\]
Theorem~\ref{thm:linear-ear-closure} applies with $G_0=C$ and gives
\[
    L(X)=\ord(N^3)
\]
for the puzzle on $X$ itself.

Under the convention of
Subsection~\ref{subsec:biconnected-bridge-tree-decomposition},
the constant number of blanks can first be brought into the working
region of $X$, and any temporarily used vertices outside $X$ can be
restored after the reconfiguration, with an additional
$\Ord(N^2)$ moves.
Hence the reconfiguration cost of $X$ is at most
\[
    L(X)+\Ord(N^2)=\ord(N^3).
\]
This contradicts Theorem~\ref{thm:linear-component-cubic}.  Hence
\[
    X\text{ contains no cycle }C
    \text{ with }|V(C)|=\omega(1)\cap\ord(N).
\]

We now prove the assertion concerning two linear-size cycles.
Suppose that $C_1,C_2\subseteq X$ satisfy
\[
    |V(C_1)|=\Theta(N),
    \qquad
    |V(C_2)|=\Theta(N),
\]
and put
\[
    s:=|V(C_1)\cap V(C_2)|.
\]
Assume, to the contrary, that
\[
   1 \leq s<b.
\]
Since there are $b$ blanks, all vertices of
$V(C_1)\cap V(C_2)$ may be made blank, with at least one further blank
remaining available.  Whenever control is transferred from one cycle
to the other, restore the blanks on all vertices of
$V(C_1)\cap V(C_2)$ before starting the next presentation operation.
Thus the common vertices serve only as blank interfaces.  With this
convention, the same two-part presentation-and-sorting construction as
in the proof of Theorem~\ref{thm-two-cycles-o-n3} applies to
$C_1\cup C_2$.
Since both cycles have order $\Theta(N)$, the
resulting bound is
\[
    \Ord(N^2\log N)=\ord(N^3).
\]
Thus $C_1\cup C_2$ is a linear-size subcubic core.

Again, because $X$ is biconnected and contains $C_1\cup C_2$, the
ear-decomposition theorem extends this core to $X$ by successively
adjoining ears.  Theorem~\ref{thm:linear-ear-closure} therefore implies
that the reconfiguration cost of $X$ is $\ord(N^3)$, contradicting
Theorem~\ref{thm:linear-component-cubic}.  
Consequently,
\[
    V(C_1)\cap V(C_2)=\emptyset
    \qquad\text{or}\qquad
    |V(C_1)\cap V(C_2)|\ge b=N-n.
\]
%Consequently,
%\[
%    |V(C_1)\cap V(C_2)|\ge b=N-n.
%\]

We next prove that $X$ is the unique component of
$\mathcal{B}(G)$ of order $\Theta(N)$.
By Theorem~\ref{thm:no-intermediate-decomposition-component}, every
component of $\mathcal{B}(G)$ has order either $\Ord(1)$ or
$\Theta(N)$.
Suppose that there is another component $Y\ne X$ with
\[
    |V(Y)|=\Theta(N).
\]
Corollary~\ref{cor:bridge-tree-components-bounded} shows that $Y$
is also biconnected.

Consider the component--attachment incidence tree of
$\mathcal{B}(G)$, and the unique path in this tree from $X$ to
$Y$.  At an attachment vertex $u$ on this path, the graph $G$
is separated into two connected sides $G_1$ and $G_2$, one
containing $X$ and the other containing $Y$.  Hence
\[
    |V(G_1)|=\Theta(N),
    \qquad
    |V(G_2)|=\Theta(N).
\]
Applying Theorem~\ref{thm-two-cycles-o-n3} to these two sides gives
\[
    \Ord(N^2\log N)=\ord(N^3),
\]
contradicting the assumption that the worst-case shortest solution
length on $G$ is $\Theta(N^3)$.
Therefore $X$ is the unique component of order $\Theta(N)$.

Combining this uniqueness with
Theorem~\ref{thm:no-intermediate-decomposition-component}, we obtain
\[
    |V(Y)|=\Ord(1)
\]
for every component
$Y\in\mathcal{B}(G)\setminus\{X\}$.

Finally, fix a cut vertex $v\in V(X)$ through which components of
$\mathcal{B}(G)$ are attached to $X$.  Let $A(v)$ denote the
connected subgraph consisting of $v$ together with all components
on the side away from the interior of $X$ that are attached through
$v$.  The complementary side contains $X$, and hence has order
$\Theta(N)$.

If
\[
    |V(A(v))|=\omega(1)\cap\ord(N),
\]
then Theorem~\ref{thm:intermediate-cut-component} applied at $v$
gives an $\ord(N^3)$ solution of the whole PMP, a contradiction.

If instead
\[
    |V(A(v))|=\Theta(N),
\]
then the two sides of the cut vertex $v$ both have order
$\Theta(N)$.  Theorem~\ref{thm-two-cycles-o-n3} again gives an
$\ord(N^3)$ solution, also a contradiction.

Therefore
\[
    |V(A(v))|=\Ord(1).
\]
In particular, the entire branch of
$\mathcal{B}(G)$ attached to $X$ through $v$ has bounded order.

This proves all the assertions.
\end{Proofof}

%%%%%%%%%%%%%%%%%%%%%%%%%%%%%%%%%%%%%%%%%%%%%%%%%%%%%%%%%%%%%%%%%%%%%%%%%%%%%%%%
%%%%%%%%%%%%%%%%%%%%%%%%%%%%%%%%%%%%%%%%%%%%%%%%%%%%%%%%%%%%%%%%%%%%%%%%%%%%%%%%
%%%%%%%%%%%%%%%%%%%%%%%%%%%%%%%%%%%%%%%%%%%%%%%%%%%%%%%%%%%%%%%%%%%%%%%%%%%%%%%%

\section{Concluding Remarks}\label{final}

In the concluding section ``Conclusion and Open Problems'' of
\cite{715921}, the authors remark that:
\begin{quote}
It would be useful to at least have an efficient algorithm which
approximates the number of moves required. For it seems that only a
small fraction of the graph puzzles actually require $\Ord(N^3)$
moves.
\end{quote}
To illustrate this phenomenon, the authors note that when the board
graph is a square grid with $N$ vertices, the generalized $15$-puzzle
can be solved in $\Ord(N^{3/2})$ moves.

Our Theorems \ref{refinement1},
\ref{refinement2-general-constant-blanc}, \ref{refinement2}, and
\ref{refinement3}, established in Section~\ref{refinement}, provide a
partial answer to this question by extending the square-grid example
to general board graphs and by presenting algorithms that achieve
efficient asymptotic bounds both on the length of the output move
sequence and on the computation time.
Indeed, in the above square-grid example we have
$\CL(G)=4=\Ord(1)$ and $\D(G)\in\Ord(\sqrt{N})$.
If the number of unoccupied vertices satisfies
\[
    2\leq N-n\in\Ord(1),
\]
then
a solution sequence of length
$\Ord(n\CL(G)\D(G))=\Ord(N^{3/2})$ can be computed in time linear in
its length.

Our results can be summarized as follows.
\begin{enumerate}
\item
For every feasible instance on an $N$-vertex tree, the length of the
shortest solution sequence is bounded by
$\Ord(N^2\log N)$.
This result not only rectifies an invalid proof given in
\cite{715921}, but also constitutes a genuine improvement over the
state-of-the-art bound presented in the most recent work
\cite{WOS:001177663400002}.

\item
For general board graphs, our algorithmic bounds and the structural
analysis of Sections~\ref{sec:pebble-presentation} and
\ref{sec:structural-restrictions} progressively restrict the possible
source of cubic complexity.
The general-graph bounds of Section~\ref{refinement} first imply that
$\Theta(N^3)$ behavior is possible only when the number of blanks is
bounded by a constant.
The Pebble Presentation Theorem,
Theorem\tume\ref{thm:pebble-presentation}, then provides an
$\Ord(N^2)$ mechanism for presenting all pebbles at an interface
vertex, and the resulting presentation order can be used as an
address system and transcribed across cut vertices.
Using this mechanism,
Theorem\tume\ref{thm:intermediate-cut-component} first excludes
intermediate-scale structures separated through a cut vertex.
Theorem\tume\ref{thm:no-intermediate-decomposition-component},
Theorem\tume\ref{thm:linear-component-cubic}, and
Corollary\tume\ref{cor:bridge-tree-components-bounded} then restrict the
possible scales and types of the components in the
biconnected--bridge-tree decomposition.
Finally, Theorem\tume\ref{thm:linear-ear-closure} shows that subcubic
reconfiguration complexity of the specified cores is preserved under
an arbitrary sequence of ear additions resulting in a biconnected graph,
with no restriction on the number or the individual lengths of the ears.

These results culminate in
Theorem\tume\ref{thm-Theta-n3-unique-large-component}.
If the worst-case shortest solution length on an $N$-vertex board
graph is $\Theta(N^3)$, then its biconnected--bridge-tree
decomposition has a unique component $X$ of order $\Theta(N)$, and
$X$ is biconnected and has reconfiguration cost $\Theta(N^3)$.
Every other component has order $\Ord(1)$, and the total structure
attached to $X$ through any one of its cut vertices also has bounded
order.
Moreover, $X$ contains no cycle of order
$\omega(1)\cap\ord(N)$, and any two intersecting cycles of order
$\Theta(N)$ in $X$ have at least $N-n$ common vertices.
Thus all unbounded structure responsible for cubic complexity is
forced into a single giant biconnected core whose internal cycle
structure is itself severely restricted.

\item
The results of Section~\ref{lower bound} show that these restrictions
are close to sharp.
Theorem\tume\ref{thm-Theta-n3} shows that a cycle containing all but
$\Ord(1)$ vertices of the graph can support the full
$\Theta(N^3)$ complexity.
In contrast, Theorem\tume\ref{thm-cycle-o-n3} shows that, within the
corresponding edge-amalgamation family, allowing an unbounded total
amount of structure outside the giant cycle while keeping every
individual attachment bounded forces the complexity below cubic.
Furthermore,
Theorem\tume\ref{thm-two-cycles-o-n3} and
Corollary\tume\ref{cor-two-cycles-o-n3} show that joining two
unbounded board graphs of feasible puzzles at a single vertex is
subcubic, while
Theorem\tume\ref{thm-two-cycles-Theta-n2logn} determines the exact
order $\Theta(N^2\log N)$ when the two parts are cycles of order
$\Theta(N)$.
These results motivate
Conjecture\tume\ref{conj-Theta-n3}, which asserts that the essentially
cyclic construction of Theorem\tume\ref{thm-Theta-n3} characterizes
the cubic regime.
\end{enumerate}

More broadly, these results offer a clearer structural understanding
of the \textsc{Pebble Motion Problem} by demonstrating that the
previously known $\Ord(N^3)$ worst-case bound can only arise under
highly restricted and specific structural configurations of the
underlying graph.
In particular, the combination of
Theorem\tume\ref{thm-Theta-n3-unique-large-component} with the
sharpness results of Section~\ref{lower bound} indicates that cubic
complexity is not merely associated with the presence of large
cycles: it requires a much more rigid concentration of the unbounded
structure inside a single biconnected core.

Next, we discuss the relationship between our results on the
\textsc{Pebble Motion Problem} (PMP) for general graphs and the
PMP-based techniques proposed in the context of
Multi-Agent Path Finding (MAPF).
Our results suggest a different perspective on algorithmic approaches
to PMP in MAPF.
Rather than treating special graph classes separately, our analysis
indicates that the tree case can serve as a natural base structure for
more general instances.
In particular, computations can be organized around a spanning
breadth-first search tree of the graph, while the remaining structural
features of the graph are addressed only when needed to resolve local
infeasibility in the tree subproblems.
This viewpoint leads to a more unified treatment of PMP across
different graph structures.
In particular, our improved understanding of the complexity of the
problem on trees---showing that it admits an
$\Ord(N^2\log N)$-time algorithm---suggests that tree-based
formulations may serve as a computationally viable foundational
framework.
We believe that this perspective may provide useful guidance for the
design of PMP solvers within MAPF systems.

Finally, we provide some intuition behind Conjecture~\ref{conjecture}.
It is already known that the shortest solution sequences for PMT may
require $\Omega(N^2)$ moves.
Moreover, there exist instances in which moving a single pebble by a
distance $\Theta(N)$ along the board tree already requires
$\Theta(N^2)$ moves.
In view of this phenomenon, establishing an upper bound of
$\tilde{\Ord}(N^2)$ on the length of the shortest solution sequence
appears difficult without relying on some form of divide-and-conquer
strategy.

To illustrate the underlying intuition, we consider the following
model case (which is also discussed in
Remark~\ref{counterexample} in the Appendix).
Let $P$ be a sufficiently long path graph of even length, and attach
a pendant edge $uv$ to the midpoint vertex $u$ of $P$.
The resulting board graph is $G:=P+uv$.
Let $P_L$ and $P_R$ denote the left and right subpaths of $P$
obtained by removing the vertex $u$.
Let $\ell:=|V(P_L)|=|V(P_R)|$.
In this board tree, the maximum isthmus size equals $\ell$, and hence
at most $\ell+1$ pebbles can be placed on the board.

Consider the following reconfiguration task.
Place $\ell/2$ pebbles on $P_L$ and $\ell/2$ pebbles on $P_R$.
The goal is to transform this configuration into its reflection
across the axis determined by the edge $uv$.
More precisely, if the pebbles on $P_L$ are arranged along the path as
$p_1,p_2,\ldots,p_k$, they must be relocated to the symmetric
positions on $P_R$ in the order $p_k,\ldots,p_2,p_1$.
At the same time, the pebbles originally on $P_R$ must be moved to
the corresponding symmetric positions on $P_L$.
To analyze this process, partition the $\ell/2$ pebbles on each of
$P_L$ and $P_R$ into $m$ groups of equal size.
In order to realize the above reconfiguration, each group on $P_L$
must exchange positions near the edge $uv$ with every group on
$P_R$.
Consequently, interactions may arise between essentially every pair
of these groups, giving rise to $m^2$ induced subproblems.

The main difficulty lies in the cost of preparing these subproblems
on the board tree.
Even if the subproblems are solved in a carefully chosen order,
our current intuition suggests that, at least for this example,
preparing each subproblem of size $\Theta(\ell/m)$ would require
$\Theta((\ell/m)^2)$ moves even when $m=\Ord(1)$.
Moreover, if $m=\omega(1)$, the cost of this preparation step
would likely grow significantly beyond the quadratic order.
Since the original instance has size $\Theta(\ell)$, this intuition
suggests that each level of recursive refinement incurs a preparation
cost of
\[
\Theta((\ell/m)^2)\times\Theta(m^2)=\Theta(\ell^2)
\]
moves.
Reducing the subproblem size down to $\Ord(1)$ requires
$\Theta(\log\ell)$ levels of such refinement.
Consequently, the total number of moves spent in these preparatory
steps accumulates to $\Theta(\ell^2\log\ell)$.
Since $\ell=\Theta(N)$ in this construction, this corresponds to
$\Theta(N^2\log N)$ moves.

This example suggests that long degree-two corridors constitute a
common structural source of difficulty in PMP: on trees, they appear
as large isthmuses that induce substantial coordination costs and,
when sufficiently long relative to the number of empty vertices, may
even render the instance infeasible, while in general graphs, where
cycles provide alternative routes, the same feature can lead, in
extreme cases, to $\Theta(N^3)$-length solution sequences.

The principal questions left open by our results concern the two main
directions developed in this paper.
For trees, it remains to determine whether the logarithmic factor in
Theorem\tume\ref{main} is necessary, as suggested by
Conjecture\tume\ref{conjecture}.
For general graphs, Conjecture\tume\ref{conj-Theta-n3} asks whether
every $\Theta(N^3)$-puzzle is, up to only $\Ord(1)$ additional
vertices, organized around a single giant cycle.
Resolving these two questions would substantially sharpen our
understanding of both the extremal complexity on trees and the
structural origin of cubic complexity on general graphs.

%%%%%%%%%%%%%%%%%%%%%%%%%%%%%%%%%%%%%%%%%%%%%%%%%%%%%%%%%%%%%%%%%%%%%%%%%%%%%%%%
%%%%%%%%%%%%%%%%%%%%%%%%%%%%%%%%%%%%%%%%%%%%%%%%%%%%%%%%%%%%%%%%%%%%%%%%%%%%%%%%
%%%%%%%%%%%%%%%%%%%%%%%%%%%%%%%%%%%%%%%%%%%%%%%%%%%%%%%%%%%%%%%%%%%%%%%%%%%%%%%%

\section*{Appendix}
\addcontentsline{toc}{section}{Appendix}

\begin{rem}\label{counterexample}
The following statement on page 247, left column, second paragraph from the bottom of Kornhauser, Miller, and Spirakis \cite{715921}:
\begin{quote}
\noindent
{\bf b. G Separable and transitive}

\vspace{5pt}

``If G is a tree, then the proof of transitivity implies that at most $\Ord(n)$ moves are needed to move a pebble anywhere;
so the proof of Theorem 1, Case 2, implies an upper bound of $\Ord(n^2)$ to move all the $k < n$ pebbles.''
\end{quote}
(where $k$ denotes the number of pebbles and $n$ is the number of vertices of the board graph $G$) does not hold in general.
% The initial part of the above statement should be corrected as follows:
% \begin{quote}
%     ``If G is a tree, then the proof of transitivity implies that at most $\Ord(n^2)$ moves are needed
%     to move a pebble anywhere;''
% \end{quote}
There exist counterexamples to the initial part of the above citation, namely,
``If G is a tree, then the proof of transitivity implies that at most $\Ord(n)$ moves are needed to move a pebble anywhere;''.
This is because, when pebbles are densely packed along a path, it may be impossible
to relocate a particular pebble unless all pebbles are moved simultaneously.
For example, suppose we are given a sufficiently long path graph $P$ of even length, 
where a pendant edge $uv$ is attached to the midpoint vertex $u$ of $P$.
Let the board graph be $G:=P+uv$.
Let $P_L$ and $P_R$ denote the left and right subpaths of $P$ obtained by removing $u$, respectively, so that $P = P_L - u - P_R$.
Note that $n=|V(G)|=|V(P)|+2$ and $|V(P_L)|=|V(P_R)|=(n-2)/2$ hold. 
Suppose that all vertices of $P_L$ are occupied by pebbles without any gaps, and no other vertex contains a pebble.
In this case, the number of pebbles is $k = (n - 2)/2 = n/2 - 1$.
%For simplicity (though this is not essential to the argument), suppose that $n - 2$ is divisible by $4$.
In this configuration, by Theorem\tume{A}, this pebble motion puzzle on the board tree $G$ is feasible and thus transitive.
Consequently, all the assumptions required for Case 2 of Theorem 1 in \cite{715921} are satisfied.
In this setup, consider the task of moving the pebble $X$ initially located at the leftmost vertex of $P_L$ to the rightmost vertex of $P_R$.
This requires $\Omega(n^2)$ moves, regardless of the strategy employed.
The pebble $X$ is blocked on its right by $n/2 - 2$ pebbles, packed tightly along the path $P_L$.
%Since direct swaps of adjacent pebbles are not allowed, $X$ must be temporarily evacuated to the leaf vertex $v$ on the pendant edge $uv$ at least once during the process.
%However, to make room for $X$ to reach $v$, all the $n/2 - 2$ pebbles to its right must first be moved to the subpath $P_R$, which has $n/2 - 1$ unoccupied vertices. 
Since direct swaps of adjacent pebbles are not allowed, $X$ must be temporarily evacuated to the base vertex $u$ on the pendant edge $uv$ at least once during the process.
However, to make room for $X$ to reach $u$, all of the $n/2 - 2$ pebbles packed immediately to the right of $X$ must at some point pass through the base vertex $u$.
Indeed, this sequence of operations already demands at least 
$$
\sum_{i=1}^{n/2 - 1} i = \frac{n(n-2)}{8} \in \Omega(n^2)
$$
moves of pebbles. 
%Specifically, by the end of this process, the right half of $P_R$ will contain at least $(n - 2)/4 - 1$ pebbles.
%Since all of these pebbles originally resided on $P_L$, each must have made at least $(n - 2)/4+2$ moves to reach their positions.
%Hence, the total number of moves incurred in this stage is at least
%$$
%\left( \frac{n - 2}{4} - 1 \right) \cdot \left(\frac{n - 2}{4} +2\right)\in \Omega(n^2).
%$$
Moreover, under the setting described above, it is straightforward to verify that, starting from a configuration $f$ in which $k(=n/2 - 1)$ pebbles
labeled $a_1$ through $a_k$ are tightly packed from left to right $(a_1,\ldots, a_k)$ beginning at the leftmost end-vertex of the left subpath $P_L$
of the board graph $G$, reconfiguring them into a target configuration $g$ where the same $k$ pebbles are tightly packed from right to left
$(a_k,\ldots,a_1)$ beginning at the rightmost end-vertex of the right subpath $P_R$ of $G$, according to the algorithm proposed in \cite{715921}
requires a sequence of $\Omega(n^3)$ moves.
Although we present here only a single counterexample, % showing a lower bound  $\Omega(n^3)$ of the algorithm in \cite{715921},
it is evident that similar situations necessarily arise whenever there exists an isthmus of size $\Omega(n)$
in the $n$-vertex board tree with $\Omega(n)$ pebbles. 
% Conversely, the fact that moving a pebble to an arbitrary position can be achieved by a sequence of $\Ord(n^2)$ moves
% requires a separate proof.
% That the procedure of the algorithm in \cite{715921} indeed achieves this has been independently verified through
% replication studies in several subsequent papers \cite{WOS:001177663400002, WOS:000077614700003}. 
% Therefore, a correct revision of their argument should read:
% \begin{quote}
% \noindent
% ``If G is a tree, then the proof of transitivity implies that at most $\Ord(n^2)$ moves are needed
% to move a pebble anywhere; so the proof of Theorem 1, Case 2, implies an upper bound of
% $\Ord(n^3)$ to move all the $k < n$ pebbles.''
% \end{quote}
% That said, immediately after the above-mentioned statement in \cite{715921},
% the authors write: "The existence of biconnected subgraphs, however, can force us to an upper bound of $\Ord(n^3)$..."
% This presents a clear contrast between the case where the board graph is a tree—suggesting a complexity of $\Ord(n^2)$—and
% the case where it contains biconnected subgraphs, in which the best attainable upper bound is only $\Ord(n^3)$.
% This contrast reveals that the error identified in \cite{715921}  is not merely a typographical mistake or a momentary oversight,
% but rather stems from a fundamental misunderstanding. 
Furthermore, the facts presented here clearly indicate that,
in order to reduce the computational complexity of the problem to strictly below $\Ord(n^3)$,
it is necessary to design algorithms in a clever way that allows amortized analysis to function effectively.
\end{rem}

\begin{Proofof}{Lemma\tume\ref{sbs0}}
We will show an algorithm of the reconstruction problem on trees for unlabeled pebbles.
%It is worth noting that in the case of unlabeled pebbles, we can echange a role of a pebble and that of an empty space.
%Hence, it is sufficient to show that there exists a sequence of moves of length $\Ord (n\D(G))$ which transfers all the pebbles from $S_1$ to $S_2$, and it can be computed in time $\Ord (N + n\D(G))$.

Let $T$ be a directed rooted tree with the root $r$.
For $v \in V(T)$, the closed descendants of $v$ is the set of vertices that are either $v$ itself or can be reached by a directed path starting from $v$.
Let us denote the closed descendants of $v$ of $T$ by ${\textrm{Des}}_{T}(v)$. 

%\begin{alg}[$\Alg.\Pac(T, \BT, x)$]\hfill
\begin{alg}[$\Alg.sbs0(G, S_1, S_2)$]\hfill
\begin{namelist}{OUTPUT}
\item[\sl Input]{A tree $G$ with $N$ vertices. Subsets of vertices $S_1, S_2  \subset V(G)$ with $|S_1| = |S_2| = n$. }
\item[\sl Output]{A sequence of moves which transfers all the pebbles from $S_1$ to $S_2$.}
\end{namelist}
\begin{alg-enumerate}
\item Choose a vertex $r$ of $G$, and build a directed rooted tree $T$ with the root $r$ such that $V(T)= V(G)$ and $E(T)=E(G)$.\label{sbs0_start0}
\item For $i=1,2$, build a directed rooted tree $T_i \subset T$ with the root $r$ such that $V(T_i) = \{ v \in V(T) \,:\, {\textrm{Des}}_{T}(v) \cap S_i \ne \emptyset \}$.
Let $\textrm{Stk}$ be a stack for delayed sequences of moves. 
\item Make a list $L$ of leaves of $T$. \label{sbs0_start1}
\item If $L = \emptyset$, then Goto step \maru{\ref{sbs0_pop}}.\label{sbs0_main}
Take a leaf $v \in L$.\\
If $v \in S_1 \cap S_2$, then:\\
$\ $ Set $S_1 = S_1 \setminus \{ v \}$, $S_2 = S_2 \setminus \{ v \}$, and update $T_1$ and $T_2$ appropriately.\\
else if $v \in S_1 \setminus S_2$, then:\\
$\ $ Let $x$ be the first vertex of $V(T_2)$ such that $x$ is on the $vr$-path.
Then there exists a vertex $y \in {\textrm{Des}}_{T_2}(x) \cap S_2$.
Let $z$ be the first vertex of $S_2$ on the $vxy$-path. 
Push the $vz$-path to $\textrm{Stk}$.\\
Set $S_1 = S_1 \setminus \{ v \}$, $S_2 = S_2 \setminus \{ z \}$, and update $T_1$ and $T_2$ appropriately.\\
else if $v \in S_2 \setminus S_1$, then:\\
$\ $ Let $x$ be the first vertex of $V(T_1)$ such that $x$ is on the $vr$-path.
Then there exists a vertex $y \in {\textrm{Des}}_{T_1}(x) \cap S_1$.
Let $z$ be the first vertex of $S_1$ on the $vxy$-path. 
Move a pebble from $z$ to $v$ in $G$.
Set $S_1 = S_1 \setminus \{ z \}$, $S_2 = S_2 \setminus \{ v \}$, and update $T_1$ and $T_2$ appropriately.\\
\item Set $T = T - v$.\label{sbs0_cut_v}
Let $u$ be the parent of $v$ in $T$.
Remove $v$ from $L$.
If $u$ is a leaf of $T$, then add $u$ to $L$.
Goto step \maru{\ref{sbs0_main}}. 
\item \label{sbs0_pop}
Pop an element, a $vz$-path, from $\textrm{Stk}$, and move a pebble from $v$ to $z$ in $G$  repeatedly, until $\textrm{Stk}$ is empty.
\end{alg-enumerate}
\end{alg}

In step \maru{\ref{sbs0_cut_v}}, $|V(T)|$ always decreases.
Hence, we have finally $L = \emptyset$.

In the algorithm, each pebble moves to a destination one by one through a path.
Since there is no other pebble on the path, the algorithm works correctly.
Moreover, since the length of each path is at most $\D(G)$, the total length of the resulting sequence of moves is $\Ord(n\D(G))$.

For steps from \maru{\ref{sbs0_start0}} to \maru{\ref{sbs0_start1}}, it is sufficient to apply Depth First Search of $G$, whose running time is $\Ord (N)$.
The total running time of the algorithm is dominated by step \maru{\ref{sbs0_main}} and step  \maru{\ref{sbs0_pop}}.
In case $v \not\in S_1 \cup S_2$, it takes a constant time.
Hence, the total running time of the case is $\Ord (N)$. 
In case $v \in S_1 \cup S_2$, it takes $\Ord (\D(G))$ time to update $T_1$ and $T_2$, and it takes $\Ord (\D(G))$ steps to move a pebble.
Since this case occurs at most $n$ times throughout the entire algorithm, the total running time of the case is $\Ord (n\D(G))$.   

Therefore, the total running time of the algorithm is $\Ord (N + n\D(G))$. 
\end{Proofof}

%%%%%%%%%%%%%%%%%%%%%%%%%%%%%%%%%%%%%%%%%%%%%%%%%%%%%%%%%%%%%%%%%%%%%%%%%%%%%%%%
%%%%%%%%%%%%%%%%%%%%%%%%%%%%%%%%%%%%%%%%%%%%%%%%%%%%%%%%%%%%%%%%%%%%%%%%%%%%%%%%
%%%%%%%%%%%%%%%%%%%%%%%%%%%%%%%%%%%%%%%%%%%%%%%%%%%%%%%%%%%%%%%%%%%%%%%%%%%%%%%%

\begin{Proofof}{Lemma\tume\ref{sbs1}}
We propose two algorithms for the proof.
The first algorithm is as follows. 

\begin{alg}[$\Alg.sbs1A(G, Q_1, Q_2)$](Partition with a rivet)\hfill
\begin{namelist}{OUTPUT}
\item[\sl Input]{A tree $G$ satisfying the assumption of Lemma\tume\ref{sbs1}.\\
A configuration $\pi \in \cF(G, P)$.\\
A partition of pebbles $P = Q_1 \cup Q_2$ such that $Q_1 \cap Q_2 = \emptyset$ and $|Q_i| = \lfloor n/2 \rfloor$ or $\lceil n/2 \rceil$ for $i=1, 2$.  
}
\item[\sl Output]{A sequence of moves which transfers $Q_i$ to $V(G_i)$ for $i = 1, 2$.}
\end{namelist}
\begin{alg-enumerate}
\item Find subsets $S_i$ of $V(G)$ for $1 \le i \le 4$, such that $|S_1| = |S_4| = \lfloor n/2 \rfloor$, $|S_2| = |S_3| = \lceil n/2 \rceil + 1$, $S_i \cap S_j \ne \emptyset$ for $1 \le i < j \le 4$, $S_1, S_2 \subset V(G_1)$, $S_3, S_4 \subset V(G_2)$, $S_2 \cap S_3 = \{ o \}$, and $S_i \cup S_{i+1}$ induces a subtree $G(i, i+1)$ of $G$ for $1 \le i \le 3$. 
\item Applying $\Alg.sbs0$, move a quarter of $P$ to $S_i$ for $1 \le i \le 4$.
Put $P_i$ as a set of pebbles on $S_i$ for $1 \le i \le 4$.
\item Consider the puzzle $(G(2,3), P_2 \cup P_3)$ and transfer half of $P_2 \cup P_3$ to $S_2$ and the other half to $S_3$ such that pebbles of $Q_1$ are moved to $S_2$ as many as possible.
Reset $P_2$ and $P_3$ as a set of pebbles on $S_2$ and $S_3$, respectively.
\item Remove $P_3 \cup P_4$ to $S_4$.
Load $P_1 \cup P_2$ to $S_2 \cup S_3$. 
Consider the puzzle $(G(2,3), P_1 \cup P_2)$ and transfer half of $P_1 \cup P_2$ to $S_2$ and the other half to $S_3$ such that pebbles of $Q_1$ are moved to $S_2$ as many as possible.
Reset $P_1$ and $P_2$ as a set of pebbles on $S_2$ and $S_3$, respectively.
\item Remove $P_1 \cup P_2$ to $S_1$.
Load $P_3 \cup P_4$ to $S_2 \cup S_3$. 
Consider the puzzle $(G(2,3), P_3 \cup P_4)$ and transfer half of $P_3 \cup P_4$ to $S_2$ and the other half to $S_3$ such that pebbles of $Q_1$ are moved to $S_2$ as many as possible.
Reset $P_3$ and $P_4$ as a set of pebbles on $S_2$ and $S_3$, respectively.
\item Remove $P_4$ to $S_4$.
Load $P_2 \cup P_3$ to $S_2 \cup S_3$. 
Consider the puzzle $(G(2,3), P_2 \cup P_3)$ and transfer half of $P_2 \cup P_3$ to $S_2$ and the other half to $S_3$ such that pebbles of $Q_1$ are moved to $S_2$ as many as possible.
Reset $P_2$ and $P_3$ as a set of pebbles on $S_2$ and $S_3$, respectively.
\end{alg-enumerate}
\end{alg}

We claim that $\Alg.sbs1A$ works correctly, it outputs a solution sequence of moves of length $\Ord (n\D(G) + n^2 \log n) $ and it can be implemented to run in time $\Ord (N + n\D(G) + n^2 \log n)$.
From steps \maru{3} to \maru{6}, $\Alg.sbs1A$ is recursively called four times with a half number of original pebbles.
The worst case is $Q_1 = P_3 \cup P_4$ initially.
In this case, after step \maru{3}, we have $Q_1 = P_2 \cup P_4$, and after step \maru{4}, we have $Q_1 = P_1 \cup P_4$, and after step \maru{5}, we have $Q_1 = P_1 \cup P_3$, and finally after step \maru{6}, we have $Q_1 = P_1 \cup P_2$.

Step \maru{2} can be implemented in $\Ord(N+n\D(G))$ time by Lemma\tume\ref{sbs0}.
Let us denote the running time of the algorithm from step \maru{3} to \maru{6} by $g_1(n)$.
We want to prove $g_1(n) \le c_1 n^2 \log_2 n$ with some constant $c_1$.
we proceed by induction on $n$.
The number of steps of the removal, loading, and restoration phase is at most $cn^2$ with some constant $c$, since $|S_i| = n/2 + \Ord (1)$ for $1 \le i \le 4$.
Then, by the inductive hypothesis, we have
\begin{eqnarray*}
g_1(n) & \le & 4 g_1(\frac{n}{2}) + c n^2 \\
 & \le & 4 c_1 (\frac{n}{2})^2 \log_2 (\frac{n}{2}) + c n^2 \\
 & = & c_1 n^2 \log_2 n - c_1 n^2 + c n^2.
\end{eqnarray*} 
By choosing $c_1 = c$, we have $g_1(n) \le c n^2 \log_2 n$, as claimed.

Next, we propose the second algorithm. 

\begin{alg}[$\Alg.sbs1B(G, \pi_1, \pi_2)$](Sort with a rivet)\hfill
\begin{namelist}{OUTPUT}
\item[\sl Input]{A tree $G$ satisfying the assumption of Lemma\tume\ref{sbs1}.\\
Configurations $\pi_1$, $\pi_2 \in \cF(G, P)$.
}
\item[\sl Output]{A sequence of moves which transforms $\pi_1$ to $\pi_2$.}
\end{namelist}
\begin{alg-enumerate}
\item Find subsets $S_i$ of $V(G)$ for $1 \le i \le 4$, such that $|S_1| = |S_4| = \lfloor n/2 \rfloor$, $|S_2| = |S_3| = \lceil n/2 \rceil + 1$, $S_i \cap S_j \ne \emptyset$ for $1 \le i < j \le 4$, $S_1, S_2 \subset V(G_1)$, $S_3, S_4 \subset V(G_2)$, $S_2 \cap S_3 = \{ o \}$, and $S_i \cup S_{i+1}$ induces a subtree $G(i, i+1)$ of $G$ for $1 \le i \le 3$. 
%Let $S = \cup_{1 \le i \le 4} S_i$. 
\item Choose a configuration $\pi_i^{\ast}$ for $i = 1, 2$ such that $|{\rm sup}(\pi_i^{\ast}) \cap S_j| = \lfloor n/4 \rfloor$ or $\lceil n/4 \rceil$ for $1 \le j \le 4$ and $\pi_i^{\ast}$ is strongly related to $\pi_i$ for $i = 1, 2$.
Transform $\pi_1$ to $\pi_1^{\ast}$.
\item Transfer $\pi_2^{\ast}(S_1 \cup S_2)$ to $S_1 \cup S_2$, and transfer $\pi_2^{\ast}(S_3 \cup S_4)$ to $S_3 \cup S_4$ by applying $\Alg.sbs1A$.
Reset $P_i$ as a set of pebbles on $S_i$ for $1 \le i \le 4$.
\item Remove $P_3 \cup P_4$ from $S_3 \cup S_4$ to $S_4$.
Load $P_1 \cup P_2$ to $S_2 \cup S_3$.
Consider the puzzle $(G(2,3), P_1 \cup P_2)$ and build a configuration $\pi_{2,1}^{\ast}$
such that $\pi_{2,1}^{\ast}$ and $\pi_2^{\ast}|_{S_1 \cup S_2}$ are strongly related.
\item Remove $P_1 \cup P_2$ from $S_2 \cup S_3$ to $S_1$.
Load $P_3 \cup P_4$ to $S_2 \cup S_3$.
Consider the puzzle $(G(2,3), P_3 \cup P_4)$ and build a configuration $\pi_{2,2}^{\ast}$ 
such that $\pi_{2,2}^{\ast}$ and $\pi_2^{\ast}|_{S_3 \cup S_4}$ are strongly related.
\item Make $\pi_2^{\ast}$.
\item Transform $\pi_2^{\ast}$ to $\pi_2$.
\end{alg-enumerate}
\end{alg}

We claim that $\Alg.sbs1B$ works correctly, it outputs a solution sequence of moves of length $\Ord (n\D(G) + n^2 \log n) $ and it can be implemented to run in time $\Ord (N + n\D(G) + n^2 \log n)$.

In step \maru{2} and step \maru{7}, since $\pi_i$ and $\pi_i^{\ast}$ are strongly related for $i = 1, 2$, 
the number of steps to transform $\pi_1$ to $\pi_1^{\ast}$ and to transform $\pi_2^{\ast}$ to $\pi_2$ are at most $\Ord (n\D(G))$, and it can be computed in time $\Ord (N + n\D(G))$ by applying $\Alg.sbs0$.
The number of steps of the removal, loading, and restoration phase is at most $cn^2$ with some constant $c$.
In step \maru{3}, by applying $\Alg.sbs1A$, we have at most $c_1 n^2 \log n$ time.
Let $f_1(n)$ be the algorithm's running time from step \maru{3} to step \maru{6}.
We want to prove $f_1(n) \le d_1 n^2 \log n$ with some constant $d_1$.
In step \maru{4} and step \maru{5}, by the inductive hypothesis, we have
\begin{eqnarray*}
f_1(n) & \le & 2 f_1(\frac{n}{2}) + c_1 n^2 \log n + c n^2 \\
 & \le & 2 d_1 (\frac{n}{2})^2 \log_2 (\frac{n}{2}) + c_1 n^2 \log n + c n^2 \\
 & = & (c_1 + \frac{d_1}{2} ) n^2 \log_2 n + (c - \frac{d_1}{2} ) n^2.
\end{eqnarray*} 
By choosing $d_1$ such that $d_1 \ge 2\max\{ c, c_1 \}$, we have $f_1(n) \le d_1 n^2 \log_2 n$.        
\end{Proofof} 

%%%%%%%%%%%%%%%%%%%%%%%%%%%%%%%%%%%%%%%%%%%%%%%%%%%%%%%%%%%%%%%%%%%%%%%%%%%%%%%%
%%%%%%%%%%%%%%%%%%%%%%%%%%%%%%%%%%%%%%%%%%%%%%%%%%%%%%%%%%%%%%%%%%%%%%%%%%%%%%%%
%%%%%%%%%%%%%%%%%%%%%%%%%%%%%%%%%%%%%%%%%%%%%%%%%%%%%%%%%%%%%%%%%%%%%%%%%%%%%%%%

\begin{Proofof}{Lemma\tume\ref{sbs2}}
We propose two algorithms for the proof.
The first algorithm is as follows. 

\begin{alg}[$\Alg.sbs2A(G, Q_1, Q_2)$](Partition with a long isthmus)\hfill
\begin{namelist}{OUTPUT}
\item[\sl Input]{A tree $G$ satisfying the assumption of Lemma\tume\ref{sbs2}.\\
A configuration $\pi \in \cF(G, P)$.\\
A partition of pebbles $P = Q_1 \cup Q_2$ such that $Q_1 \cap Q_2 = \emptyset$ and $|Q_i| = \lfloor n/2 \rfloor$ or $\lceil n/2 \rceil$ for $i=1, 2$.  
}
\item[\sl Output]{A sequence of moves which transfers $Q_i$ to $V(G_i)$ for $i = 1, 2$.}
\end{namelist}
\begin{alg-enumerate}
\item Find subsets $S_i$ of $V(G)$ for $1 \le i \le 4$, such that $|S_i| = \lfloor n/4 \rfloor$ or $\lceil n/4 \rceil$ for $1 \le i \le 4$, $S_i \cap S_j \ne \emptyset$ for $1 \le i < j \le 4$, $S_1, S_2 \subset V(G_1)$, $S_3, S_4 \subset V(G_2)$, and $S_i \cup S_{i+1} \cup V(I)$ induces a subtree $G(i, i+1)$ for $1 \le i \le 3$.
%Let $v_1$ and $v_2$ be the endvertices of $I$ such that $v_i \in V(G_i)$ for $i=1, 2$. 
\item Applying $\Alg.sbs0$, move a quarter of $P$ to $S_i$ for $1 \le i \le 4$.
Put $P_i$ as a set of pebbles on $S_i$ for $1 \le i \le 4$.
\item Consider the puzzle $(G(2,3), P_2 \cup P_3)$ and transfer half of $P_2 \cup P_3$ to $S_2$ and the other half to $S_3$ such that pebbles of $Q_1$ are moved to $S_2$ as many as possible.
Reset $P_2$ and $P_3$ as a set of pebbles on $S_2$ and $S_3$, respectively.
\item Consider the puzzle $(G(1,2), P_1 \cup P_2)$ and transfer half of $P_1 \cup P_2$ to $S_1$ and the other half to $S_2$ such that pebbles of $Q_1$ are moved to $S_1$ as many as possible.
Reset $P_1$ and $P_2$ as a set of pebbles on $S_1$ and $S_2$, respectively.
\item Consider the puzzle $(G(3,4), P_3 \cup P_4)$ and transfer half of $P_3 \cup P_4$ to $S_3$ and the other half to $S_4$ such that pebbles of $Q_1$ are moved to $S_3$ as many as possible.
Reset $P_3$ and $P_4$ as a set of pebbles on $S_3$ and $S_4$, respectively.
\item Consider the puzzle $(G(2,3), P_2 \cup P_3)$ and transfer half of $P_2 \cup P_3$ to $S_2$ and the other half to $S_3$ such that pebbles of $Q_1$ are moved to $S_2$ as many as possible.
Reset $P_2$ and $P_3$ as a set of pebbles on $S_2$ and $S_3$, respectively.
\end{alg-enumerate}
\end{alg}

We claim that $\Alg.sbs2A$ works correctly, it outputs a solution sequence of moves of length $\Ord (n\D(G) + n^2 \log n) $ and it can be implemented to run in time $\Ord (N + n\D(G) + n^2 \log n)$.

In step \maru{4} and step \maru{5}, $\Alg.sbs1A$ is called two times with a half number of original pebbles, and 
in step \maru{3} and step \maru{6}, $\Alg.sbs2A$ is recursively called two times with a half number of original pebbles.
The worst case is $Q_1 = P_3 \cup P_4$ initially.
In this case, after step \maru{3}, we have $Q_1 = P_2 \cup P_4$, and after step \maru{4}, we have $Q_1 = P_1 \cup P_4$, and after step \maru{5}, we have $Q_1 = P_1 \cup P_3$, and finally after step \maru{6}, we have $Q_1 = P_1 \cup P_2$.

The running time $g_2(n)$ of the algorithm is dominated by step \maru{3} - step \maru{6}.
We want to prove $g_2(n) \le n\D(G) + c_2 n^2 \log_2 n$ with some constant $c_2$.
we proceed by induction on $n$.
The number of steps of removal, loading, and restoration phase is at most $cn^2$ with some constant $c$, since $|S_1 \cup S_2| = n/2 + \Ord (1)$ and $|S_3 \cup S_4| = n/2 + \Ord (1)$.
The cost of step \maru{4} and step \maru{5} is at most $c_1 (n/2)^2 \log_2 (n/2)$ by Lemma\tume\ref{sbs1}.

Then, by the inductive hypothesis, we have
\begin{eqnarray*}
g_2(n) & \le & 2 g_2(\frac{n}{2}) + 2c_1 (\frac{n}{2})^2 \log_2 (\frac{n}{2}) + c n^2 \\
  & \le & 2 \left( \frac{n}{2}\D(G) + c_2 \left(\frac{n}{2}\right)^2 \log_2 \left(\frac{n}{2}\right) \right)
 + 2c_1 \left(\frac{n}{2}\right)^2 \log_2 \left(\frac{n}{2}\right) + c n^2 \\
 & \le & n\D(G) + (\frac{c_2}{2} + \frac{c_1}{2}) n^2 \log_2 n + (c - \frac{c_2}{2} - \frac{c_1}{2} )n^2.
\end{eqnarray*} 
By choosing $c_2$ such that $c_2 \ge \max \{ c_1, 2c - c_1 \}$, we have $g_2(n) \le n\D(G) + c_2 n^2 \log_2 n$, as claimed.              

Next, we propose the second algorithm. 

\begin{alg}[$\Alg.sbs2B(G, \pi_1, \pi_2)$](Sort with a long isthmus)\hfill
\begin{namelist}{OUTPUT}
\item[\sl Input]{A tree $G$ satisfying the assumption of Lemma\tume\ref{sbs2}.\\
Configurations $\pi_1$, $\pi_2 \in \cF(G, P)$.
}
\item[\sl Output]{A sequence of moves which transforms $\pi_1$ to $\pi_2$.}
\end{namelist}
\begin{alg-enumerate}
\item Find subsets $S_i$ of $V(G)$ for $1 \le i \le 4$, such that $|S_i| = \lfloor n/4 \rfloor$ or $\lceil n/4 \rceil$ for $1 \le i \le 4$, $S_i \cap S_j \ne \emptyset$ for $1 \le i < j \le 4$, $S_1, S_2 \subset V(G_1)$, $S_3, S_4 \subset V(G_2)$, and $S_i \cup S_{i+1} \cup V(I)$ induces a subtree $G(i, i+1)$ for $1 \le i \le 3$.
%Let $v_1$ and $v_2$ be the endvertices of $I$ such that $v_i \in V(G_i)$ for $i=1, 2$. 
\item Choose $\pi_i^{\ast}$ for $i = 1, 2$ such that $|{\rm sup}(\pi_i^{\ast}) \cap S_j| = \lfloor n/4 \rfloor$ or $\lceil n/4 \rceil$ for $1 \le j \le 4$ and $\pi_i^{\ast}$ is strongly related to $\pi_i$ for $i = 1, 2$.
Transform $\pi_1$ to $\pi_1^{\ast}$.
\item Consider the puzzle $(G, P)$ and transfer $\pi_2^{\ast}(S_1 \cup S_2)$ to $S_1 \cup S_2$ and transfer $\pi_2^{\ast}(S_3 \cup S_4)$ to $S_3 \cup S_4$ by applying $\Alg.sbs2A$.
Reset $P_i$ as a set of pebbles on $S_i$ for $1 \le i \le 4$.
\item 
Consider the puzzle $(G(1,2), P_1 \cup P_2)$ and build a configuration $\pi_{2,1}^{\ast}$ 
such that $\pi_{2,1}^{\ast}$ and $\pi_2^{\ast}|_{S_1 \cup S_2}$ are strongly related.
\item
Consider the puzzle $(G(3,4), P_3 \cup P_4)$ and build a configuration $\pi_{2,2}^{\ast}$ 
such that $\pi_{2,2}^{\ast}$ and $\pi_2^{\ast}|_{S_3 \cup S_4}$ are strongly related.
\item Make $\pi_2^{\ast}$.
\item Transform $\pi_2^{\ast}$ to $\pi_2$.
\end{alg-enumerate}
\end{alg}

We claim that $\Alg.sbs2B$ works correctly, it outputs a solution sequence of moves of length $\Ord (n\D(G) + n^2 \log n) $ and it can be implemented to run in time $\Ord (N + n\D(G) + n^2 \log n)$.

In step \maru{2} and step \maru{7}, since $\pi_i$ and $\pi_i^{\ast}$ are strongly related for $i = 1, 2$, the number of steps to transfer $\pi_1$ to $\pi_1^{\ast}$ and to transfer $\pi_2^{\ast}$ to $\pi_2$ are at most $\Ord (n\D(G))$ by applying $\Alg.sbs0$.
In step \maru{3}, $\Alg.sbs2A$ is called, and 
in step \maru{4} and step \maru{5}, $\Alg.sbs1B$ is called two times with a half number of original pebbles.

The running time $f_2(n)$ of the algorithm is dominated by step \maru{3} - step \maru{6}.
We want to prove $f_2(n) \le n\D(G) + d_2 n^2 \log_2 n$ with some constant $d_2$.
The number of steps of removal, loading, and restoration phase is at most $cn^2$ with some constant $c$.
The cost of  step \maru{3} is at most $Nn + c_2 n^2 \log_2 n$ by $\Alg.sbs2A$.
The cost of  steps \maru{4} and \maru{5} is at most $d_1 (n/2)^2 \log_2 (n/2)$ by $\Alg.sbs1B$.
Then, we have
\begin{eqnarray*}
f_2(n) & \le & n\D(G) + c_2 n^2 \log_2 n + 2 d_1 (\frac{n}{2})^2 \log_2 (\frac{n}{2}) + c n^2 \\
 & \le & n\D(G) + ( c_2 + \frac{d_1}{2}) n^2 \log_2 n + (c - \frac{d_1}{2})n^2.
\end{eqnarray*} 
By choosing $d_2$ such that $d_2 \ge c_2 + c$, we have $f_2(n) \le n\D(G) + d_2 n^2 \log_2 n$.              
\end{Proofof}

%%%%%%%%%%%%%%%%%%%%%%%%%%%%%%%%%%%%%%%%%%%%%%%%%%%%%%%%%%%%%%%%%%%%%%%%%%%%%%%%
%%%%%%%%%%%%%%%%%%%%%%%%%%%%%%%%%%%%%%%%%%%%%%%%%%%%%%%%%%%%%%%%%%%%%%%%%%%%%%%%
%%%%%%%%%%%%%%%%%%%%%%%%%%%%%%%%%%%%%%%%%%%%%%%%%%%%%%%%%%%%%%%%%%%%%%%%%%%%%%%%
\begin{Proofof}{Theorem\tume\ref{sbs3}}
Since $G$ is a tree, $G$ has a centroid vertex or a centroid edge.
We consider two cases.
\medskip\\
\underline{Case 1.} There exists a centroid $o$ of degree at least $3$.
\medskip\\
In this case, since $N \ge 3n$, $G$ has subtrees $G_1$ and $G_2$ such that $V(G_1) \cap V(G_2) = \{ o \}$ and $|V(G_i)| \ge n + 1$ for $i=1, 2$.
Hence, by Lemma\tume\ref{sbs1}, the proof is finished.
\medskip\\
\underline{Case 2.} There exists no centroid of degree at least $3$.
\medskip\\  
In this case, the centroid is contained in an isthmus $I$.
Let $v_1$ and $v_2$ be endvertices of $I$, and let $V_0$ be the set of inner vertices of $I$.
We have two subtrees $G_1$ and $G_2$ induced by $V(G) \setminus V_0$ such that $v_i \in V(G_i)$ for $i = 1, 2$.
Put $N_i = |V(G_i)|$ for $i = 1, 2$.
By the feasibility of the puzzle, we have $n \le N_1 + N_2 - 3$.
We may assume $N_1 \le N_2$. 
Let $k= |V_0|$. 
\medskip\\
\underline{Subcase 2. 1.} $k < n$.
\medskip\\  
In this case, we have $N_2 \ge |V(G) \setminus V_0|/2 = (N-k)/2 > (N - n)/2 \ge n$.
On the other hand, since $I$ contains the centroid, we have $|V(G_1) \cup V(I)| \ge N/2 \ge (3/2)n$.
Hence, by taking $v_2$ as the rivet in the assumption of Lemma\tume\ref{sbs1}, we can apply Lemma\tume\ref{sbs1} to finish the proof.
\medskip\\
\underline{Subcase 2. 2.} $k \ge n$.
\medskip\\
Let $H$ be a subtree of $G$ induced by $V(I) \cup V(G_2)$, and 
let $K$ be a subtree of $G$ induced by $V(G_1) \cup V(I) \cup V(G'_2)$, 
where $G'_2$ is a subtree of $G_2$ containing $v_2$ with $N_1$ vertices.
Let us take positive integers $n_1$, $n_2$ such that $n_1 + n_2 = n$ and $n_i \le N_i - 1$ for $i=1,2$.  
We will use $\Alg.sbs1B$ and $\Alg.sbs2B$ to solve the puzzle $(G,P)$ as follows;
\begin{alg-enumerate}
\item{Apply $\Alg.sbs1B$ for the puzzle $(H, n_2)$ with a rivet $v_2$.}
\item{Apply $\Alg.sbs2B$ for the puzzle $(K, 2n_1)$ with an isthmus $I$.}
\item{Again, apply $\Alg.sbs1B$ for the puzzle $(H, n_2)$ with a rivet $v_2$.}
\end{alg-enumerate}
The running time of steps \maru{1} - \maru{3} is $\Ord(N + n\D(G) + n^2\log n)$ by Lemma\tume\ref{sbs1} and Lemma\tume\ref{sbs2}.
This completes the proof.
\end{Proofof}

\begin{Proofof}{Lemma\tume\ref{sbs4}}
Since $G$ is a tree, $G$ has a centroid vertex or a centroid edge.
We consider two cases.
\medskip\\
\underline{Case 1.} There exists a centroid $o$ of degree at least $3$.
\medskip\\
In this case, since $N \ge 3n$, $G$ has subtrees $G_1$ and $G_2$ such that $V(G_1) \cap V(G_2) = \{ o \}$ and $|V(G_i)| \ge n + 1$ for $i=1, 2$.
Since $\deg \, o \ge 3$, we may assume that there exists a set of $3$ neighbors $v_1, v_2, v_3$ of $o$ such that $v_1 \in G_1$ and $v_2, v_3 \in G_2$.
Set $S = \{ o, v_1, v_2, v_3 \}$.

We will use the following algorithm;
\begin{alg-enumerate}
\item Starting from $\pi$, gather all the pebbles in $V(G_1)$.
\item Move pebbles to $V(G_2)$ through $o$.
If $p$ or $q$ passes through $o$, stop them at $v_2$ or $v_3$.
Let all pebbles except $p$ and $q$ pass through $S$ and put them in $V(G_2) \setminus \{ v_2, v _3 \}$, until both $p$ and $q$ arrive at $S$.
\item Interchange $p$ and $q$ in $S$.
\item Follow all of the steps before interchanging $p$ and $q$  in reverse order to achieve the desired configuration $\pi(p\leftrightarrow q)$.
\end{alg-enumerate}
\noindent
\underline{Case 2.} There exists no centroid of degree at least $3$.
\medskip\\  
In this case, the centroid is contained in an isthmus $I$.
Let $v_1$ and $v_2$ be endvertices of $I$, and let $V_0$ be the set of inner vertices of $I$.
We have two subtrees $G_1$ and $G_2$ induced by $V(G) \setminus V_0$ such that $v_i \in V(G_i)$ for $i = 1, 2$.
Put $N_i = |V(G_i)|$ for $i = 1, 2$.
By the feasibility of the puzzle, we have $n \le N_1 + N_2 - 3$.
We may assume $N_1 \le N_2$. 
Let $k= |V_0|$. 
\medskip\\
\underline{Subcase 2. 1.} $k < n$.
\medskip\\  
In this case, we have $N_2 \ge |V(G) \setminus V_0|/2 = (N-k)/2 > (N - n)/2 \ge n$.
On the other hand, since $I$ contains the centroid, we have $|V(G_1) \cup V(I)| \ge N/2 \ge (3/2)n$.
Hence, by taking $v_2$ as the vertex $o$ in Case 2.1, we can apply the same algorithm in Case 2.1 to finish the proof.
\medskip\\
\underline{Subcase 2. 2.} $k \ge n$.
\medskip\\
Let $H$ be a subtree of $G$ induced by $V(I) \cup V(G_2)$, and 
let $K$ be a subtree of $G$ induced by $V(G_1) \cup V(I) \cup V(G'_2)$, 
where $G'_2$ is a subtree of $G_2$ containing $v_2$ with $N_1$ vertices.
Let us take positive integers $n_1$, $n_2$ such that $n_1 + n_2 = n$ and $n_i \le N_i - 1$ for $i=1,2$.  

Since $\deg \, v_i \ge 3$ for $i=1,2$, we may assume that there exists a set of vertices $u_{1,1}, u_{1,2}, u_{2,1}, u_{2,2}$ such that $u_{i,1}, u_{i,2}$ are neighbors of $v_i$ in $G_i$ for $i=1,2$.
Set $S = V(I) \cup \{ u_{1,1}, u_{1,2}, u_{2,1}, u_{2,2} \}$.

We will use the following algorithm;
\begin{alg-enumerate}
\item Starting from $\pi$, gather all the pebbles in $V_0$.
\item Move $n_i$ pebbles to $V(G_i)$ through $v_i$ for $i=1,2$.
If $p$ or $q$ passes through $v_i$, stop them at $u_{i,1}$ or $u_{i,2}$.
Let all pebbles except $p$ and $q$ pass through $v_1$ or $v_2$ and put them in $(V(G_1) \cup V(G_2))\setminus \{ u_{1,1}, u_{1,2}, u_{2,1}, u_{2,2} \}$.
\item If both $p$ and $q$ arrive at common $v_i$ for $i=1$ or $2$, then interchange $p$ and $q$ around $v_i$.  
Otherwise, we may assume $p$ arrives at $v_1$ and $q$ arrives at $v_2$.
Let all pebbles except $p$ and $q$ pass through $v_1$ or $v_2$ and put them in $(V(G_1) \cup V(G_2))\setminus \{ u_{1,1}, u_{1,2}, u_{2,1}, u_{2,2} \}$, then interchange $p$ and $q$ in $S$.
\item Follow all of the steps before interchanging $p$ and $q$ in reverse order to achieve the desired configuration $\pi(p\leftrightarrow q)$.
\end{alg-enumerate}
In all cases, the total running time of steps is $\Ord(N + n\D(G))$, since each pebble moves $\Ord(\D(G))$ times in the algorithms.
This completes the proof.
\end{Proofof}

%%%%%%%%%%%%%%%%%%%%%%%%%%%%%%%%%%%%%%%%%%%%%%%%%%%%%%%%%%%%%%%%%%%%%%%%%%%%%%%%
%%%%%%%%%%%%%%%%%%%%%%%%%%%%%%%%%%%%%%%%%%%%%%%%%%%%%%%%%%%%%%%%%%%%%%%%%%%%%%%%
%%%%%%%%%%%%%%%%%%%%%%%%%%%%%%%%%%%%%%%%%%%%%%%%%%%%%%%%%%%%%%%%%%%%%%%%%%%%%%%%

\begin{Proofof}{Lemma\tume\ref{centroid-splitting}}
If $T$ has a centroid edge, then by definition the statement of this lemma clearly holds.
Hence let us assume that $T$ has a centroid vertex $o$.
Let us consider the set of components $\cS$ of $T-o$. 
Without loss of generality, let $X$ be a graph with the fewest order, 
and $Y$ be a graph with the second-fewest order in the set $\cS$. 
Of course it is possible that $|V(X)|=|V(Y)|$. 
If $m \geq 3$, since $|V(X)|+|V(Y)| \leq \frac{2(N-1)}{3}$ holds, let us denote  
$Z:=X \cup Y$ and update $\cS:=(\cS \setminus \{X,Y\})\cup \{Z\}$. 
Repeating this operation, the size of the set S becomes $2$. 
At that stage, let us denote $S=\{G,H\}$, and let us assume $|V(G)|\leq |V(H)|$.
In this case, we can take $T_1$ to be the subgraph of $T$ induced by 
$V(G) \cup \{o\}$, and similarly $T_2$ to be the subgraph of $T$ induced by 
$V(H) \cup \{o\}$. 
Since the size of each element of the set $\cS$ is between $1$ and $n$, 
we can sort and manage each element of $\cS$ in order of its size using 
a $1\times n$-array $A(\cS)$ whose each cell is a list. 
Since the minimum of the orders of the graphs in $\cS$ is 
monotonically non-decreasing, the task of identifying a graph with 
the smallest order and a graph with the second smallest order in $A(\cS)$
takes only $\Ord(N)$ effort in total.
\end{Proofof}

%%%%%%%%%%%%%%%%%%%%%%%%%%%%%%%%%%%%%%%%%%%%%%%%%%%%%%%%%%%%%%%%%%%%%%%%%%%%%%%%
%%%%%%%%%%%%%%%%%%%%%%%%%%%%%%%%%%%%%%%%%%%%%%%%%%%%%%%%%%%%%%%%%%%%%%%%%%%%%%%%
%%%%%%%%%%%%%%%%%%%%%%%%%%%%%%%%%%%%%%%%%%%%%%%%%%%%%%%%%%%%%%%%%%%%%%%%%%%%%%%%

\begin{Proofof}{Lemma\tume\ref{feasible-subtree}}
%$T^{\prime}$が長さ$1$のパス$(K_2)$である場合には、$T^{\prime}$のいずれかの点を$v$と置き、$H:=T$とする。それ以外の場合には、$T^{\prime}$のleaf以外のすべての点を$T$上で$1$頂点$v$にまで縮約した木を$H$と置く。さらに$T^{\prime}$に対応する$H$の部分グラフを$X$と置く。
%$X$に属していない$v$の隣接頂点$u$が$H$に見つかる限り、その頂点$u$と辺$vu$を$X$に追加（$X:=X;vu$に更新）し続ける。
%$X$に属していない$v$の隣接頂点$u$が$H$に存在しなくなったら、$X$のどこかのleaf $\ell$に隣接する$X$外の頂点$w$と辺$\ell{w}$を$X$に追加する。その結果$X$に新しく出来た長さ$2$のarm $v-\ell-w$の$v$に接続する辺$v\ell$を縮約して$\ell$を$v$に同一視したものを新たに$X$に置き直し、$H$上で$v\ell$を縮約した結果得られる新たなグラフを$H$に置き直して、$v$の隣接頂点を$X$に追加するフェイズに戻る。
Let us consider the following algorithm:
%下記の一連の作業を$X$に対応する$T$の点部分グラフの位数が$d$になったところで止めて、復元して出力すればよい。
\begin{alg-enumerate}
\item{If $T^{\prime}$ consists of a single edge $e$, let $v$ be an end-vertex of $e$ and set $H:=T$. 
Otherwise, let $H$ be the graph obtained by contracting all parts of $T^{\prime}$ except the pendant edges 
to a single vertex $v$. Let $X$ be the subgraph of $H$ corresponding to $T^{\prime}$. }
%\item{$T^{\prime}$が長さ$1$のパス$(K_2)$である場合には、$T^{\prime}$のいずれかの点を$v$と置き、$H:=T$とする。それ以外の場合には、$T^{\prime}$のleaf以外のすべての点を$T$上で$1$頂点$v$にまで縮約した木を$H$と置く。さらに$T^{\prime}$に対応する$H$の部分グラフを$X$と置く。}
\item{Continue to add vertex $u$ and edge $vu$ to $X$ (i.e., update to $X:=X+vu$) as long as 
the rank of the subgraph $T(X)$ of $T$ corresponding to $X$ is strictly less than $d$, and 
a neighboring vertex $u$ of $v$ that does not belong to $X$ is found in $H$.}\label{leaf-add}
%\item{$X$に属していない$v$の隣接頂点$u$が$H$に見つかる限り、その頂点$u$と辺$vu$を$X$に追加（$X:=X+vu$に更新）し続ける。}\label{leaf-add}
\item{If the rank of the subgraph $T(X)$ of $T$ corresponding to $X$ is strictly less than $d$, do:
\begin{alg-enumerate}
\item{Find a vertex $w$ outside $X$ adjacent to some leaf $\ell$ in $X$ and add that vertex $w$ and edge $\ell{w}$ to $X$. 
The resulting $X$ has the length $2$ arm $v-\ell-w$.}
\item{Update $X$ (resp. $H$) as the graph obtained by contracting the edge $v\ell$ of $X$ (resp. $H$) and identifying $\ell$ with $v$. }
\item{Go to Step\tume\maru{\ref{leaf-add}}}
\end{alg-enumerate}
}\label{spider}
%\item{$X$のどこかのleaf $\ell$に隣接する$X$外の頂点$w$を見つけて、その頂点$w$と辺$\ell{w}$を$X$に追加する。その結果$X$に生じた長さ$2$のarm $v-\ell-w$上の辺$v\ell$を縮約して$\ell$を$v$と同一視させたグラフを$X$に取り直し、$H$の辺$v\ell$を縮約したグラフを$H$に取り直して、ステップ\maru{\ref{leaf-add}}に戻る。}\label{spider}
\item{Return the subgraph $T(X)$ of $T$ corresponding to $X$.}
\end{alg-enumerate}

Let us show the correctness of the above algorithm. 
%上記のアルゴリズムの正当性について述べる。
To do so, we need only check the fact that the maximal isthmus size of the subgraph $T(X)$ of $T$ corresponding to $X$ 
does not exceed the maximal isthmus size of $T$ by the operation of Step\tume{\maru{\ref{leaf-add}}} or Step\tume{\maru{\ref{spider}}}.
%それには、$X$に対応する$T$の部分グラフ$T(X)$の最大isthmus長が、ステップ{\maru{\ref{leaf-add}}}やステップ{\maru{\ref{spider}}}の操作によって
%$T$の最大isthmus長を上回ることがない、という事実を確認すればよい。
First, for the operation of Step\tume\maru{\ref{leaf-add}}, the maximum isthmus size of $T(X)$ is not increased, 
since the degree of the non-leaf vertices in $T(X)$ is increased.
%まず、ステップ\maru{\ref{leaf-add}}の操作は、$T(X)$におけるleaf以外の頂点の次数しか増やさないので、最大isthmus長は減ることはあっても増えることはない。
Next, the vertex $w$ that is added to $T$ in Step\tume\maru{\ref{spider}} operation is not a cut vertex because it is a leaf of $T(X)$, 
so it cannot be part of any isthmus in $T(X)$.
%次に、ステップ\maru{\ref{spider}}の操作で$T(X)$に付け加える頂点$w$は$T(X)$のleafなのでcut vertexではないから、$T(X)$のisthmusを
%構成する頂点ではあり得ない。
Therefore, the new isthmus arising in $T(X)$ at this step has one end vertex at $\ell$. 
%したがって、この段階で$T(X)$に生じる新たなisthmusは一方の端点を$\ell$に持つ。
On the other hand, just before entering the operation of Step\tume\maru{\ref{spider}}, 
the degree on $T(X)$ and the degree on $T$ are the same for any vertex other than leaves 
in $T(X)$, so this newly generated isthmus in $T(X)$ is just an isthmus of $T$ itself. 
%一方、ステップ\maru{\ref{spider}}の操作に入る直前の状態においては、$T(X)$におけるleaf以外のどの頂点についても、$T(X)$上の次数と$T$上の次数は一致しているので、$T(X)$に新たに生じたこのisthmusは$T$のisthmusに過ぎない。
It is clear that the time complexity of this algorithm is linear. 
%なお、このアルゴリズムの計算量が線形時間であることは明らかである。
\end{Proofof}

%%%%%%%%%%%%%%%%%%%%%%%%%%%%%%%%%%%%%%%%%%%%%%%%%%%%%%%%%%%%%%%%%%%%%%%%%%%%%%%%
%%%%%%%%%%%%%%%%%%%%%%%%%%%%%%%%%%%%%%%%%%%%%%%%%%%%%%%%%%%%%%%%%%%%%%%%%%%%%%%%
%%%%%%%%%%%%%%%%%%%%%%%%%%%%%%%%%%%%%%%%%%%%%%%%%%%%%%%%%%%%%%%%%%%%%%%%%%%%%%%%

\begin{Proofof}{Lemma\tume\ref{bunkatsu}}
Let us consider the operation $X$ mentioned in Lemma\tume\ref{centroid-splitting} of dividing a given tree into two subtrees. 
%$T$のcentroid vertexもしくはcentroid edgeを境にして、小さい方の部分木の位数が$n/3$以上となるような2つの部分木$T_1, T_2$（一般性を失わず、$|V(T_1)|=<|V(T_2)|$を仮定）に分けるという操作$X$について考えます。
Let us recursively apply the operation $X$ to each subtree of $T$ to the point such that the order of each subtree is greater than or equal to $k$ and that another $X$ on those subtrees will reduce the order of the resulting subtrees to strictly less than $k$.
%上記の「centroidを境に2つの部分木に分ける操作$X$」を、$T$の各部分木が「位数$k$を下回らないギリギリのサイズ」になるまで再帰的に適用します。

Consequently, if $r$ is the number of subtrees obtained so far, then $2n/k > (n-1)/(k-1) \geq r$.
%その結果、ここまでに得られた部分木の個数を$r$とすると、$2n/k > (n-1)/(k-1) \geq r$が成り立ちます。
In fact, one $X$ operation increases the number of subtrees by one.
%実際、このアルゴリズムにおいて、操作$X$を1回施すと、部分木の個数は1つ増えます。
Therefore, the total overlap of the vertices of the subtrees is at most $r-1$.
%したがって、ここまでに得られた部分木の個数を$r$とすると、部分木の点の重複度は高々$r-1$に過ぎません。
From this it follows that $((r-1)+n)/k \geq r$, and solving this inequality yields $n-1 \geq (k-1)r$.
%このことから、$((r-1)+n)/k \geq  r$ であり、この不等式を解くと、$n-1 \geq (k-1)r$ です。
Now, obviously, if $k=1$, the degree of overlap is zero, hence we can assume $k \geq 2$, so we have $2n/k > (n-1)/(k-1) \geq r$.
%今、明らかに、$k=1$ならば重複度はゼロですから$k$は$2$以上としてよいので、$2n/k > (n-1)/(k-1) \geq  r$ です。
Therefore, even if we apply one more $X$ operation to each subtree here, the number of subtrees can be kept to $4n/k$.
%よって、ここで各subtreeを高々もう一回ずつcentroidで分解したとしても、その個数は$4n/k$個で抑えられます。
Each subtree produced by this final splitting step has order at most
$k$.  Moreover, its parent before the split has order at least $k$, and
Lemma~\ref{centroid-splitting} implies that each of its two children has
order greater than $k/3$.  Thus every resulting subtree $H$ satisfies
\[
\frac{k}{3}<|V(H)|\le k.
\]
Therefore $T$ is covered by at most $4n/k$ subtrees satisfying all the
conditions of the lemma.
%従って、このアルゴリズムにより、$T$は、高々$4n/k$個の位数$k$以下の$T$の部分木達により被覆されます。
\end{Proofof}

%%%%%%%%%%%%%%%%%%%%%%%%%%%%%%%%%%%%%%%%%%%%%%%%%%%%%%%%%%%%%%%%%%%%%%%%%%%%%%%%
%%%%%%%%%%%%%%%%%%%%%%%%%%%%%%%%%%%%%%%%%%%%%%%%%%%%%%%%%%%%%%%%%%%%%%%%%%%%%%%%
%%%%%%%%%%%%%%%%%%%%%%%%%%%%%%%%%%%%%%%%%%%%%%%%%%%%%%%%%%%%%%%%%%%%%%%%%%%%%%%%

\begin{Proofof}{Lemma\tume\ref{children}}
According to the definition of the subtree stored in the second child of each node in this binary tree $\BT$ (cf. Step\tume\ref{No2Child}), for the sibling node  $c_2(p(\ell(i)))$ of $\ell(i)$, the graph $X(r)-(X(\ell(i))-\JoiV(X(\ell(i))))$ includes the graph $X(c_2(p(\ell(i))))$. 
The statement of this lemma can be directly derived from this fact. 
%第2子ノードに格納する木の選び方（ステップ{\ref{No2Child}}）から、$\ell(i)$の兄弟(Sibling)である第2子leafノード$c_2(p(\ell(i)))$について、
%$X(c_2(p(\ell(i))))$は$X(r)-(X(\ell(i))-\JoiV(X(\ell(i))))$の中に含まれる。当該命題の言明はこの事実から直ちに導かれる。
\end{Proofof}

%%%%%%%%%%%%%%%%%%%%%%%%%%%%%%%%%%%%%%%%%%%%%%%%%%%%%%%%%%%%%%%%%%%%%%%%%%%%%%%%
%%%%%%%%%%%%%%%%%%%%%%%%%%%%%%%%%%%%%%%%%%%%%%%%%%%%%%%%%%%%%%%%%%%%%%%%%%%%%%%%
%%%%%%%%%%%%%%%%%%%%%%%%%%%%%%%%%%%%%%%%%%%%%%%%%%%%%%%%%%%%%%%%%%%%%%%%%%%%%%%%

\begin{Proofof}{Lemma\tume\ref{BT}}
%Lemma~\ref{bunkatsu}の証明で用いられるアルゴリズムを用いる。
Let $\NS(\ell)$ be the set of nodes that are at distance $\ell$ from the root node of the binary tree $\BT$.
%二分木$\BT$のルートノードから距離$\ell$にあるノードの集合を$\NS(\ell)$と置く。
Since $\forall i, \sum_{x \in\ NS(i)}|V(X(x))|\leq 4n$ holds, 
according to Lemma~\ref{centroid}, for any distance $\ell$, the effort required to correctly store all the objects that should be placed in all nodes of $\NS(\ell)$ is at most $\Od(n)$. 
%$\forall i, \sum_{x \in \NS(i)}|V(X(x))|\leq 4n$であるから、Lemma~\ref{centroid}を勘案すれば、任意の距離$\ell$に対して、$\NS(\ell)$の構成には、高々$\Od(n)$程度の計算量しか掛からない。
The height of $\BT$ is $\Od(\log n)$, so the total effort required to complete $\BT$ is $\Od(n\log n)$. 
%この二分木$\BT$の高さは高々$\Od(\log n)$であるから、この二分木$\BT$を完成させるのに必要な手間は、$\Od(n\log n )$である。
\end{Proofof}

%%%%%%%%%%%%%%%%%%%%%%%%%%%%%%%%%%%%%%%%%%%%%%%%%%%%%%%%%%%%%%%%%%%%%%%%%%%%%%%%
%%%%%%%%%%%%%%%%%%%%%%%%%%%%%%%%%%%%%%%%%%%%%%%%%%%%%%%%%%%%%%%%%%%%%%%%%%%%%%%%
%%%%%%%%%%%%%%%%%%%%%%%%%%%%%%%%%%%%%%%%%%%%%%%%%%%%%%%%%%%%%%%%%%%%%%%%%%%%%%%%

\begin{Proofof}{Lemma\tume\ref{constant}}
In this case, by applying Lemma\tume\ref{feasible-subtree} and consecutively solving at most 
$\Ord(\D(G))$ feasible sub-puzzles --- each on a board tree with at most $2M+2$ vertices --- 
along the path connecting any two pebbles, it is possible to swap their positions 
without altering the positions of the other pebbles. Since the size of each board tree is 
bounded by the constant $2M+2$ in these sub-puzzles, the effort and number of moves 
to rearrange any configuration on such a board tree is also in constant order $\Ord(1)$. 
Note that, in general, after performing a suitable preprocessing of $\Ord(n\D(G))$ steps, 
the entire configuration can be completed by repeating the operation of 
'swapping the positions of two pebbles' at most $n-1$ times. 
%%%%%%%%%%%%%%%%%%%%%%%%%%%%%%%%%%%%%%%%%%%%%%%%%%%%%%%%%%%%%%%%%%%%%%%%
%%%%%%%%%%%%%%%%%%%%%%%%%%%%%%%%%%%%%%%%%%%%%%%%%%%%%%%%%%%%%%%%%%%%%%%%
%%%%%%%%%%%%%%%%%%%%%%%%%%%%%%%%%%%%%%%%%%%%%%%%%%%%%%%%%%%%%%%%%%%%%%%%
%%%%%%%%%%%%%%%%%%%%%%%%%%%%%%%%%%%%%%%%%%%%%%%%%%%%%%%%%%%%%%%%%%%%%%%%
%%%%%%%%%%%%%%%%%%%%%%%%%%%%%%%%%%%%%%%%%%%%%%%%%%%%%%%%%%%%%%%%%%%%%%%%
%%%%%%%%%%%%%%%%%%%%%%%%%%%%%%%%%%%%%%%%%%%%%%%%%%%%%%%%%%%%%%%%%%%%%%%%
%%%%%%%%%%%%%%%%%%%%%%%%%%%%%%%%%%%%%%%%%%%%%%%%%%%%%%%%%%%%%%%%%%%%%%%%
Moreover, in general, for any vertex $u$ of $G$ occupied by a pebble $A$, 
and any other vertex $v$ occupied either by a different pebble $B$ or 
by an unoccupied blank space, the positions of pebble $A$ and pebble $B$ 
(or the blank space) can be exchanged—without altering the positions of 
the other pebbles—by performing at most $\Ord\left(\frac{\D(G)}{k}\right)$
sub-puzzles of size $\Ord(k)$ along a shortest path in $G$ connecting $u$ 
and $v$, whose length is at most $\D(G)$.
%%%%%%%%%%%%%%%%%%%%%%%%%%%%%%%%%%%%%%%%%%%%%%%%%%%%%%%%%%%%%%%%%%%%%%%%
%%%%%%%%%%%%%%%%%%%%%%%%%%%%%%%%%%%%%%%%%%%%%%%%%%%%%%%%%%%%%%%%%%%%%%%%
%%%%%%%%%%%%%%%%%%%%%%%%%%%%%%%%%%%%%%%%%%%%%%%%%%%%%%%%%%%%%%%%%%%%%%%%
%%%%%%%%%%%%%%%%%%%%%%%%%%%%%%%%%%%%%%%%%%%%%%%%%%%%%%%%%%%%%%%%%%%%%%%%
%%%%%%%%%%%%%%%%%%%%%%%%%%%%%%%%%%%%%%%%%%%%%%%%%%%%%%%%%%%%%%%%%%%%%%%%
%%%%%%%%%%%%%%%%%%%%%%%%%%%%%%%%%%%%%%%%%%%%%%%%%%%%%%%%%%%%%%%%%%%%%%%%
%%%%%%%%%%%%%%%%%%%%%%%%%%%%%%%%%%%%%%%%%%%%%%%%%%%%%%%%%%%%%%%%%%%%%%%%
From these observations, it immediately follows that, in this case, even if scanning 
the entire structure of the input tree $G$ requires up to $\Ord(N)$ time, 
the \textsc{Pebble Motion Problem on Trees} can still be computed in time  
$\Ord(N+n\D(G))$, and the length of the resulting 
solution sequence is bounded by $\Ord(n\D(G)).$
\end{Proofof}

%%%%%%%%%%%%%%%%%%%%%%%%%%%%%%%%%%%%%%%%%%%%%%%%%%%%%%%%%%%%%%%%%%%%%%%%%%%%%%%%
%%%%%%%%%%%%%%%%%%%%%%%%%%%%%%%%%%%%%%%%%%%%%%%%%%%%%%%%%%%%%%%%%%%%%%%%%%%%%%%%
%%%%%%%%%%%%%%%%%%%%%%%%%%%%%%%%%%%%%%%%%%%%%%%%%%%%%%%%%%%%%%%%%%%%%%%%%%%%%%%%

\begin{Proofof}{Lemma\tume\ref{subtrees-containing-v}}
The operations in
Step\tume\maru{\ref{hige}} of
$\Alg.\CutOff(G,T,\BT,f,g^{\prime})$
and
Step\tume\maru{\ref{hige2}} of
$\Alg.\Extraction(T,\BT,g^{\prime},t)$
are identical for the purpose of the complexity estimate, and the
same is true of the deletions in
Steps\tume\maru{\ref{del-t}} and
\maru{\ref{del-ell}}.
We therefore discuss only the corresponding operations in
$\Alg.\CutOff(G,T,\BT,f,g^{\prime})$.

When the vertices of
$X(t)-\JoiV(X(t))$ are removed from the current black subtree, the
possible change in the maximum isthmus size is confined to an
$\Ord(k)$-vertex region around $\Joi(X(t))$.
Hence whether $T_{\mathrm{black}}$ contains an isthmus of size at
least $k+1$ can be determined in $\Ord(k)$ time.

During the initial construction of $\BT$, for each vertex of $T$,
store the list of leaf nodes $w$ of $\BT$ for which the vertex
belongs to $X(w)$, together with the position of each leaf in the
corresponding lists.
These data can be constructed in total $\Ord(n+k)$ time.
They allow a leaf node $w$ satisfying $v\in V(X(w))$ to be found in
$\Ord(1)$ time and allow all list entries belonging to a deleted
leaf to be removed in $\Ord(k)$ time.

Consequently, each of
Steps\tume\maru{\ref{hige}},
\maru{\ref{del-t}},
\maru{\ref{hige2}}, and
\maru{\ref{del-ell}}
can be carried out in $\Ord(k)$ time.
\end{Proofof}

%%%%%%%%%%%%%%%%%%%%%%%%%%%%%%%%%%%%%%%%%%%%%%%%%%%%%%%%%%%%%%%%%%%%%%%%%%%%%%%%
%%%%%%%%%%%%%%%%%%%%%%%%%%%%%%%%%%%%%%%%%%%%%%%%%%%%%%%%%%%%%%%%%%%%%%%%%%%%%%%%
%%%%%%%%%%%%%%%%%%%%%%%%%%%%%%%%%%%%%%%%%%%%%%%%%%%%%%%%%%%%%%%%%%%%%%%%%%%%%%%%

\begin{Proofof}{Lemma\tume\ref{mini-jyunbi}}
Let $\mathcal{L}$ be the set of current leaf nodes of $\BT2$.
Let $H$ be the minimal subtree of the tree induced by the current black
vertices that contains, for every $\ell\in\mathcal{L}$, the joint
$\Joi(X(\ell))$, where a joint may be either a joint vertex or a joint
edge.

Throughout the repetitions of
Step\tume\maru{\ref{mini-puzzle}} of
$\Alg.\Extraction(T,\BT,g^{\prime},t)$,
the temporary subtree $X(\ell)$ to be processed is attached to the
remaining black part of $T$ through its joint.
Hence the connected set of $k+1$ blanks can be moved from one such joint
to the next along $H$.
Equivalently, the whole movement of these blanks can be organized by a
depth-first traversal of $H$.

Since $H$ is a subtree of the current black part of $T$, we have
\[
|V(H)|\leq |V(T)|=n+k+1.
\]
Thus one blank travels only $\Ord(n+k)$ edges during the entire traversal.
Moving the connected set of $k+1$ blanks therefore requires
\[
\Ord\bigl(k(n+k)\bigr)
=
\Ord\bigl(k(k+n)\bigr)
\]
moves in total.
\end{Proofof}

\section*{AI methodology.} 
All theorems, ideas, proof strategies, and mathematical content of this paper are due solely to 
the authors. We used ChatGPT (OpenAI) and Claude (Anthropic) to assist in writing and revising 
the text of this paper. Neither tool contributed to the research itself or to its mathematical 
content. The authors verified the correctness and originality of all content, including references.

%%%%%%%%%%%%%%%%%%%%%%%%%%%%%%%%%%%%%%%%%%%%%%%%%%%%%%%%%%%%%%%%%%%%%%%%%%%%%%%%
%%%%%%%%%%%%%%%%%%%%%%%%%%%%%%%%%%%%%%%%%%%%%%%%%%%%%%%%%%%%%%%%%%%%%%%%%%%%%%%%
%%%%%%%%%%%%%%%%%%%%%%%%%%%%%%%%%%%%%%%%%%%%%%%%%%%%%%%%%%%%%%%%%%%%%%%%%%%%%%%%

\bibliographystyle{alphaurl}
\bibliography{pebble-motion-problem-on-trees}
\end{document}